\documentclass{article}
\usepackage{amsmath}
\usepackage{amsthm}

\newtheorem{theorem}{Theorem}[section]

\newtheorem{lem}[theorem]{Lemma}
\newtheorem{prop}[theorem]{Proposition}
\newtheorem{conjecture}[theorem]{Conjecture}
\theoremstyle{definition}

\newtheorem{rem}[theorem]{Remark}
\newtheorem{example}[theorem]{Example}
\numberwithin{equation}{section}
\usepackage{amssymb}
\usepackage{amsfonts}
\usepackage{epsfig}
\usepackage{psfrag}
\newcommand{\A}{\mathbb{A}}
\newcommand{\B}{\mathbb{B}}
\newcommand{\C}{\mathbb{C}}
\newcommand{\N}{\mathbb{N}}
\newcommand{\R}{\mathbb{R}}
\newcommand{\I}{\mathbb{I}}
\newcommand{\Z}{\mathbb{Z}}
\newcommand{\OO}{\mathbb{O}}
\newcommand{\Q}{\mathbb{Q}}
\newcommand{\F}{\mathbb{F}}
\newcommand{\J}{\mathbb{J}}
\newcommand{\DD}{\mathbb{D}}
\newcommand{\cC}{\mathcal{C}}
\newcommand{\cF}{\mathcal{F}}
\newcommand{\rr}{{\sf Re}\, }
\newcommand{\ii}{{\sf Im}\, }
\renewcommand{\delta}{b}
\begin{document}

\title{On planar curves with position-dependent curvature}

 \author{Arno Berger \\Mathematical and Statistical Sciences\\
University of Alberta\\Edmonton, Alberta, {\sc Canada}\\[14mm]
{\em Dedicated to Theodore P.\ Hill on the occasion of many a birthday}\\[-22mm]}


\maketitle

\vspace*{10mm}

\begin{abstract}

\noindent
Motivated by homothetic solutions in curvature-driven flows of planar
curves, as well as their many physical applications, this article carries out
a systematic study of oriented smooth curves whose curvature $\kappa$ 
is a given function of position or direction. The analysis is informed
by a dynamical systems point of view. Though focussed on situations 
where the prescribed curvature depends only on the distance $r$ from 
one distinguished point, the basic dynamical concepts are seen to be 
applicable in other situations as well. As an application, a complete 
classification of all closed solutions of $\kappa = ar^{\delta}$, with 
arbitrary real constants $a,\delta$, is established.
\end{abstract}

\noindent
\hspace*{8.3mm}{\small {\bf Keywords.} Curve shortening flow, planar curve,
  curvature, net winding, Jordan solution.}

\noindent
\hspace*{8.3mm}{\small {\bf MSC2010.} 34C05, 34C025, 34C35, 53A04,
  53C44, 58F05, 70H06.}


\section{Introduction}\label{sec1}

Modulo rotations and translations, an oriented smooth planar curve
is completely determined by its {\em curvature\/} \cite{BG, doCarmo,
  klingenberg}. Naturally, therefore, curvature plays a central role in
the study of planar shapes. The evolution and characterization of
planar shapes has been studied extensively and in a great variety of
contexts, including curve flows \cite{AL, a1, andrews, ayala, urbas}, 
growth and abrasion processes \cite{DG, DSV, FDK, Hill}, optimization 
problems \cite{dubins, KC}, among many others. In these
contexts, curvature typically is but one aspect of a more complicated
process or model. Leaving aside all additional layers of complexity,
the present article aims at classifying those planar
curves whose curvature simply is a (given) function of position or
direction. Formally, given any smooth function $k:U\times S^1 \to
\R$, where $U\subset \C$ is open, it aims at studying all smooth curve
solutions $Z$ of
\begin{equation}\label{eq1a_1}
\kappa = k (Z,N) \, ;
\end{equation}
here $\kappa\in \R$ and $N\in S^1$ denote the curvature and unit
normal vector at $Z(s) \in \C$, respectively. In particular, the
article asks whether or not (\ref{eq1a_1}) allows for solutions
that have further desirable properties such as being, for instance,
simple, closed, or convex. The main goal is to address these questions
in a systematic way, informed by a dynamical systems point of view. 
Although the pertinent dynamical ideas apply in greater generality,
for concreteness most of the analysis is focussed on situations
where $k$ in (\ref{eq1a_1}) depends only on $|Z|$ and is
independent of $N$, that is, on the special case 
\begin{equation}\label{eq1_5}
\kappa = f(|Z|) \, ,
\end{equation}
where $f:\R^+\to \R$ is a given smooth function. Thus, whereas solving
(\ref{eq1a_1}) amounts to finding planar curves with prescribed
position- and direction-dependent curvature, in (\ref{eq1_5}) the
prescribed curvature only depends on the distance from one distinguished
point (namely, the origin). 

\begin{figure}[ht]
\psfrag{tx1}[]{\small $\rr z$}
\psfrag{txx1}[]{\small $\rr Z$}
\psfrag{ty1}[]{\small $1$}
\psfrag{ty2}[]{\small $2$}
\psfrag{ty3}[]{\small $3$}
\psfrag{ty4}[]{\small $4$}
\psfrag{ty5}[]{\small $5$}
\psfrag{ty6}[]{\small $6$}
\psfrag{tx2}[]{\small $\ii z$}
\psfrag{txx2}[]{\small $\ii Z$}
\psfrag{tee1}[]{\small $1$}
\begin{center}
\includegraphics{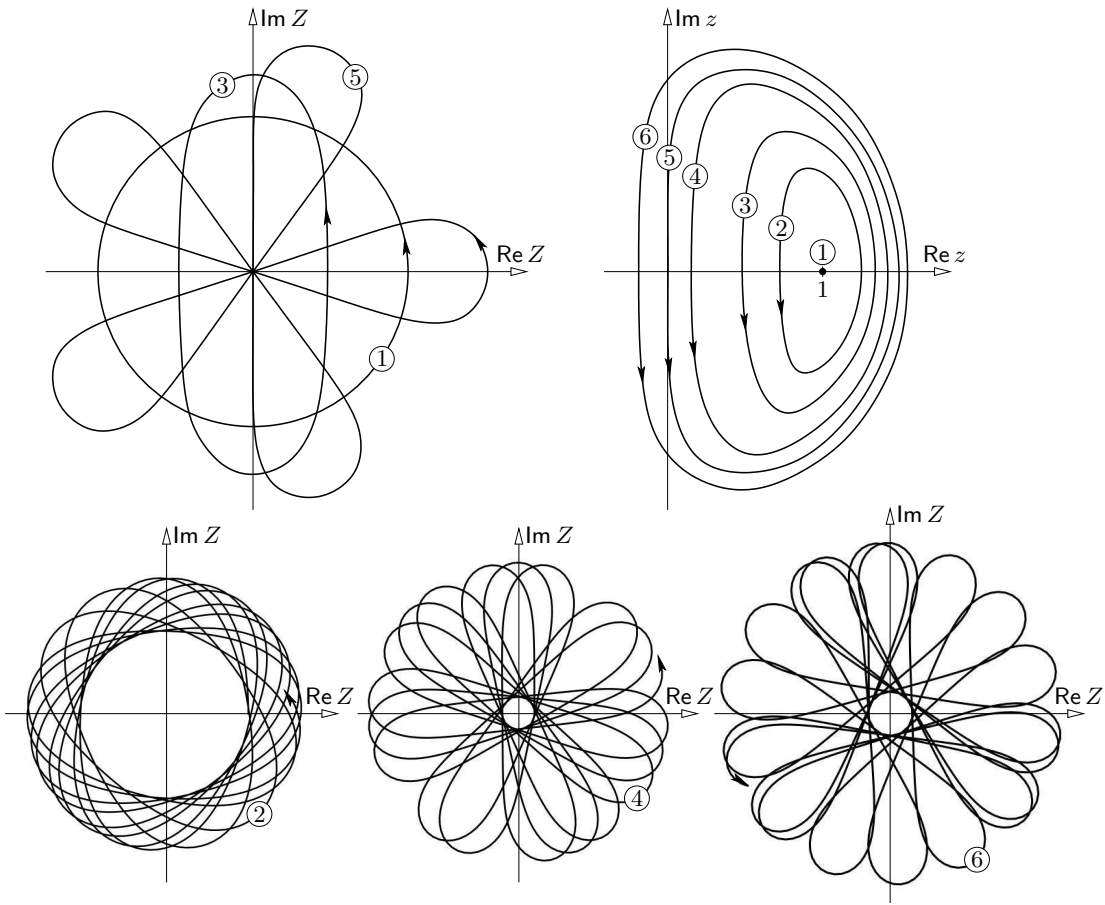}
\end{center}
\caption{
Among all smooth-curve solutions of $\kappa =
|Z|^4$, modulo rotations only precisely two are closed and simple: the
unit circle and one non-circular oval, labelled $1$ and $3$,
respectively (top left). The structure of {\em all\/} solutions is best
understood by means of the associated planar flow generated by
(\ref{eq1_6}), notably its phase portrait (top right).}\label{fig0}
\end{figure}

One key tool in this article is a planar (topological) flow
$\Phi_f$ associated with (\ref{eq1_5}). Properties of $\Phi_f$
translate into properties of solutions of (\ref{eq1_5}) that
may be hard to recognize directly. For a simple illustration,
take for instance $f(r) = r^4$. Figure \ref{fig0} displays a few
solutions of (\ref{eq1_5}) in this case; see also the illustrations in
\cite{urbas}. As it turns out, the family
of {\em all\/} (maximal) solutions of $\kappa = |Z|^4$ is most easily understood by
considering the flow on $\C$ generated by
\begin{equation}\label{eq1_6}
\dot z = i z |z|^4 - i \, .
\end{equation}
Notice that (\ref{eq1_6}) is Hamiltonian, with $H(z) = \rr z  - \frac16|z|^6$ being strictly
convex and having a unique, non-degenerate global
minimum at $z=1$. Moreover, one of the main results of the present article,
Theorem \ref{thm57} below, implies that modulo rotations, $\kappa = |Z|^4$
has precisely {\em two\/} solutions that are simple closed
(counter-clockwise oriented) curves: the unit circle and one non-circular oval;
see Figure \ref{fig0}.

One prominent natural source for (\ref{eq1a_1}) are curvature-driven
curve flows which continue to be studied for their deep mathematical
properties as well as their broad physical applications. Consider for
example the {\em Andrews--Bloore flow\/} generated by 
\begin{equation}\label{eq1_2}
\frac{\partial Z_t}{\partial t} = - (a + \kappa^b)  N \, ;
\end{equation}
here $(Z_t)_{0\le t < T}$ with the appropriate $T>0$
describes a parametrized family of oriented simple closed smooth
curves, and $a,b\ge 0$ are real constants. Curve flows
such as Andrews--Bloore have been suggested as simple models for a
variety of physical processes, ranging from the growth of crystal
surfaces to the abrasion of pebbles. As detailled
in \cite{FDK}, some aspects of the dynamics of (\ref{eq1_2}) in
general remain a challenge, both mathematically and computationally. However,
important special cases, including $b=0$ (also referred to as the {\em eikonal
  flow}), $a=0$, $b=1$ (the {\em curve-shortening flow}), and $a=0$, $b=
\frac13$ (the {\em affine curve-shortening flow}), are now well
understood; see, e.g., \cite{AL, a1, andrews, urbas} and the many references therein.
Since they may represent limiting or equilibrium shapes, {\em homothetic\/}
solutions, i.e., solutions $Z_t$ that are mere $t$-dependent
rescalings of one fixed curve $Z$, play a key role in any
analysis of (\ref{eq1_2}). Assume for instance that $Z_t(s) = \varphi(t) Z(s)$, with
$\varphi:[0,T[\: \to \R^+$, is a solution of
(\ref{eq1_2}) with $a=0$. Then $\dot \varphi = \lambda \varphi^{-b}$
for some real constant $\lambda \ne 0$, and $Z$ is a solution of
\begin{equation}\label{eq1_3}
\kappa = \big| \lambda \, \rr (Z \overline{N})\big|^{1/b} \, .
\end{equation}
While the (counter-clockwise oriented) circle with radius
$|\lambda|^{-1/(1+b)}$ centered at the origin obviously solves
(\ref{eq1_3}), it is much less obvious whether or not any other (simple closed) solutions
exist. In essence, this question is answered in
\cite{andrews} by means of an ad-hoc analysis, with the answer
depending on $b$ in a non-trivial way. Note that (\ref{eq1_3}) has the form (\ref{eq1a_1}). 
Towards the end of the present article, it will become clear that an analysis similar to, but simpler than the one
presented here can be carried out, for instance, when $k$ in
(\ref{eq1a_1}) depends only on $\rr (Z \overline{ N})$, a scenario that includes (\ref{eq1_3}).

The unique non-circular oval solution of $\kappa = |Z|^4$ is readily seen to {\em not\/} be an
ellipse. This observation nicely contrasts the ellipticity of all limiting
shapes for the affine curve-shortening flow \cite{andrews,
  FDK}. From a physical point of view, therefore, if (\ref{eq1_5}) were to be
interpreted as describing (rescaled) limiting shapes of a
curvature-driven abrasion process on pebbles, say, then such an
interpretation, however physically questionable it may be in other respects, would at
least be consistent with the well-documented empirical observation that shapes
of worn stones are not exactly elliptical either, but rather appear to
be a bit bulkier \cite{DG, DSV, F, Hill}.

\subsubsection*{Organization and notation}

This article is organized as follows. Section \ref{sec2}
introduces the planar flow $\Phi_f$ associated with (\ref{eq1_5})
and discusses a few of its basic properties. Section \ref{sec3}
establishes the crucial correspondance between the dynamics of
$\Phi_f$ on the one hand and solutions of (\ref{eq1_5}) on the
other. For a reasonably wide class of smooth functions $f$, Section \ref{sec4} characterizes
simple closed solutions of (\ref{eq1_5}) and provides tools to find
all such solutions, or else to prove that none exist. Section \ref{sec5} employs
the machinery developed in earlier sections for an in-depth study of the monomial
family $f(r) = r^{\delta}$, where $\delta$ is an arbitrary real constant. A
concluding supplemental section illustrates how the analysis, though
quite specific and delicate, nonetheless is representative
of arguments and techniques that can be applied to other
classes of functions in (\ref{eq1_5}), as well as variants thereof; as
such, it is quite similar in spirit to analyses in, e.g., \cite{andrews, benguria, MY, y}.

The following, mostly standard symbols, notation, and terminology are used throughout. The sets
of all positive integers, non-negative integers, integers, rational,
positive real, real, and complex numbers are
denoted $\N$, $\N_0$, $\Z$, $\Q$, $\R^+$, $\R$, and $\C$, respectively, and
$\R_{\infty}:= \R \cup \{-\infty, \infty\}$ is the extended real line with its familiar
order, topology, and arithmetic \cite[Sec.1.22]{Rudin}. As usual,
$\varnothing$ is the empty set, with $\inf \varnothing := \infty$ and
$\sup \varnothing := -\infty$. Limits of
real-valued objects (such as sequences, functions, or integrals) are understood in $\R_{\infty}$ unless
stated otherwise. The terms {\em increasing\/} (respectively, {\em decreasing\/})
for sequences $(a_n)$ in $\R_{\infty}$ are interpreted strictly, i.e.,
$a_n>a_m$ (respectively, $a_n < a_m$) whenever $n>m$, and similarly for
functions. Usage of an inequality such as, e.g., $a>b$ with $a,b\in
\C$ is understood to automatically imply that $a,b\in \R_{\infty}$. 
Numerical values of real numbers are displayed to
four correct significant decimal digits. The real part, imaginary
part, complex conjugate, and Euclidean norm of $z\in \C$ are $\rr z$,
$\ii z$, $\overline{z}$, and $|z|$, respectively. For convenience, let
$\C_{\times} := \C \setminus \{0\}$, $\DD = \{z\in \C : |z| < 1\}$, and
$S^1 = \partial \DD = \{z\in \C : |z| = 1\}$. Also, for every $p,w\in
\C$ and $\A \subset \C$ let $p+w\A = \{p+wz : z\in \A\}$, as well as
$\overline{\A} = \{\overline{z} : z \in \A\}$ and $\mbox{\rm dist}
(p,\A) = \inf_{z\in \A} |p-z|$, and denote the cardinality of $\A$ by
$\# \A$. Moreover, $[p,w]$ is the closed line segment with end-points $p,w$, i.e., $[p,w] =
\{(1-t) p+tw : 0\le t \le 1\}$, and similarly for the open line
segment $]p,w[$ etc. Given any function $f:\A\to \R_{\infty}$, write the set $\{ z \in \A : f( z) =
a\}$ simply as $\{f=a\}$, and its complement in $\A$ as $\{ f \ne
a\}$. As usual, $\Gamma$ denotes the Euler Gamma function.
In a slight abuse of familiar notation, let
${\rm O}(2)$ be the group of all isometries of $\C$ that fix
$0$. Recall that for every $Q\in {\rm O}(2)$ there exist $\vartheta
\in \R$ and $\epsilon_Q \in \{-1,1\}$ such that $Q(z) = e^{\textstyle i\vartheta}
(\rr z + i \epsilon_Q \ii z)$ for all $z\in \C$. Say that
$\A, \B \subset \C$ are ${\rm O}(2)$-{\bf congruent} if $Q (\A)  = \B$ for some
$Q\in {\rm O}(2)$.

\section{An auxiliary planar flow}\label{sec2}

Throughout, let $f:\R^+ \to \R$ be a smooth, that is, $C^{\infty}$-function, with additional
properties specified explicitly whenever needed. Given $f$, fix an $F : \R^+ \to \R$ with $F'(s) = s f(s) - 1$ for all
$s\in \R^+$. (The particular choice of $F$ is not going to matter
prior to Proposition \ref{prop46} below.) Also, let $\F_f = \bigl\{ z\in
\C_{\times}  : z f(|z|) = 1 \bigr\}$, a (possibly empty) subset of the
real axis. On $\C_{\times}$, consider the ODE for $z=z(t)$,
\begin{equation}\label{eq21}
\dot z = i z f(|z|) - i  \, .
\end{equation}
Recall that (\ref{eq21}) generates a local flow $\Phi_f$ on
$\C_{\times}$, in that for every $p \in \C_{\times}$ there exists
an open interval $\J \subset \R$ with $0\in \J$ such that $z(t) = \Phi_f
(t,p)$ for all $t\in \J$ yields the unique, non-extendable solution of
(\ref{eq21}) with $z(0) = p$; see, e.g., \cite{arnold, BG, perko}. Since 
$\big| {\rm d}|z|/{\rm d}t \big| =   \big|-\ii z / |z|
\big| \le 1$, it is impossible that $z(t)\to \infty$ for
finite $t$. By contrast, it is possible that $z(t) \to 0$ for finite
$t$; as recorded in Proposition \ref{lem21} below, however, such behaviour is
unproblematic. More precisely, $\Phi_f$ can be extended to a
unique (global) flow on $\C$, henceforth
also denoted $\Phi_f$. (Here and throughout, the term {\em flow\/} is
understood to mean {\em topological flow}, so $\Phi_f$ corresponds to
a one-parameter group of homeomorphisms of $\C$; see
\cite[Sec.1.I]{Irv} and Remark
\ref{rem23}(v) below.) If $\limsup_{s\to 0} |f'(s)|<\infty$ then this
immediately follows from standard facts \cite{amann, arnold, perko, walter}, but
otherwise a more tailor-made argument is required. With a
view towards the subsequent analysis of (\ref{eq21}), one such
argument utilizes the smooth function $H_f : \C_{\times}\to \R$ given by
$$
H_f (z) = \rr z  - F(|z|) - |z| \quad \forall z \in
\C_{\times} \, .
$$
Evidently, $H = H_f$ has the property that $H(z)
- \rr z$ is constant along circles centered at $0$, i.e., 
\begin{equation}\label{eq23}
\mbox{\rm for every $s\in \R^+$ there exists $a\in \R$ such that }
H(z) - \rr z  = a \enspace \mbox{\rm for all $|z| = s$} \, .
\end{equation} 
Every smooth function $H: \C_{\times}\to \R$ satisfying (\ref{eq23})
equals $H_f$ for an appropriate $f$ --- simply take $f(s) = (1- {\rm
  d}H(s)/{\rm d}s )/ s$. Most importantly, $H_f$ is a first integral of (\ref{eq21}),
and hence the study of the latter ODE for the most part reduces to an analysis of the
level sets of $H_f$. Notice that $H_f (\overline{z}) = H_f(z)$ for all
$z\in \C_{\times}$, so all level sets are symmetric w.r.t.\ the real axis. Correspondingly,
(\ref{eq21}) has a basic symmetry as well: If $z(\cdot)$ is a solution then so
is its {\bf conjugate-reversal} $\overline{z (- \, \cdot \, )}$. Utilizing $H_f$, it is a routine exercise to
establish the basic fact alluded to earlier.

\begin{prop}\label{lem21}
Let $f:\R^+ \to \R$ be smooth. Then {\rm (\ref{eq21})} generates a
unique flow $\Phi_f$ on $\C$, and
\begin{equation}\label{eq24}
\overline{\Phi_f (t,z)} = \Phi_f (-t, \overline{z}) \quad \forall (t,z) \in \R \times \C
\, .
\end{equation}
\end{prop}

\noindent
For the flow $\Phi_f$, the time-$t$ map $\Phi_f (t, \, \cdot \, )$ is
a homeomorphism of $\C$ for every $t\in \R$, and $\Phi_f (\R, z) :=
\{\Phi_f (t, z) : t \in \R\}$ is the {\bf orbit} of $z\in \C$. To
exploit the basic symmetry (\ref{eq24}), say that $z,p\in \C$ are
$\Phi_f$-{\bf conjugate} if $\Phi_f(t,z) \in \{ p, \overline{p} \}$
for some $t\in \R$. Plainly, $\Phi_f$-conjugacy is an equivalence relation, and
$z,p$ are $\Phi_f$-conjugate if and only if $\Phi_f (\R,z)$ equals
either $\Phi_f(\R, p)$ or $\overline{\Phi_f (\R, p)}$. For every $z\in \C$,
let $T_f(z) = \inf\{ t \in \R^+ : \Phi_f(t,z) = z\}$. Thus $z\in \C$ is a
periodic (respectively, fixed) point of $\Phi_f$, or $\Phi_f (\R, z)$
is a periodic orbit, in symbols $z\in
\mbox{\rm Per}\, \Phi_f$ (respectively, $z \in \mbox{\rm Fix}\,
\Phi_f$), if and only if $T_f(z)<\infty$ (respectively, $T_f(z) =
0$). Notice that if $z\in \mbox{\rm Per}\, \Phi_f$ then $z,w$ are
$\Phi_f$-conjugate precisely if $\Phi_f (t,z)=w$ for some $0\le t < T_f(z)$.
Also, $\mbox{\rm Fix}\, \Phi_f \setminus \{ 0  \}= \F_f$,
whereas $0$ may or may not be a fixed point; see Examples \ref{exa28}
to \ref{exa210} below. To clarify the nature of the fixed points of
$\Phi_f$, recall the notions of a (topological) saddle and a
center, e.g., from \cite[Sec.2.10]{perko}. Specifically, $z\in \C$
is a {\bf center} of $\Phi_f$ if a punctured neighbourhood of $z$ is
the disjoint union of periodic orbits, each of which has $z$ in its
interior. For instance, the linearization of (\ref{eq21}) at $z\in \F_f$ is 
\begin{equation}\label{eq28}
\dot p =  i F''(|z|) \, \rr p  - f(|z|) \, \ii p  \, ,
\end{equation}
suggesting that $z$ is a center when $z F''(|z|)>0$, and a saddle
when $z F''(|z|) < 0$; see Lemma \ref{lem24} below. The origin may be a fixed point 
as well; if it is, it cannot be a saddle due to (\ref{eq23}), but may
be a center, a situation that is straightforward to characterize.

\begin{prop}\label{prop22}
Let $f:\R^+ \to \R$ be smooth. Then $0$ is a center of
$\Phi_f$ if and only if both of the following conditions hold:
\begin{enumerate}
\item $\mbox{\rm dist} (0, \F_f) >0$;
\item for every $a\in \R$ there exists a decreasing sequence $(s_n)$
  with $\lim_{n\to \infty} s_n =0$ such that $|F(s_n) + s_n + a| >
  s_n$ for all $n$.
\end{enumerate}
\end{prop}

\begin{rem}\label{rem23}
(i) By Proposition \ref{prop22}, the point $0$ is a center of
$\Phi_f$ precisely if it is not an accumulation point of $\F_f$, and
$F(0+) := \lim_{s\to 0} F(s)$ either does not exist in $\R$, or
else $\bigl( F(s)  - F(0+)\bigr)/s \not \in [-2,0]$ for
arbitrarily small $s\in \R^+$. Letting $a = \liminf_{s\to 0} s
|f(s)|$, therefore, $0$ is a center when $a>1$, but is not a
center when $a<1$; if $a=1$ then $0$ may or may not be a
center as the examples $f(s) = 1 + 1/s$ and $f(s) = 1/s$ show.

(ii) Though this is of no direct consequence for the present article, note
that (\ref{eq21}) actually is Hamiltonian since
$$
\frac{{\rm d} \rr z}{{\rm d}t} = \frac{\partial H_f}{\partial \ii z}
\, ,  \quad
\frac{{\rm d} \ii z}{{\rm d}t} = - \frac{\partial H_f}{\partial \rr z}
\, .
$$
Being a 1-DOF Hamiltonian flow severely constrains the dynamical
complexity of $\Phi_f$, notably the nature of its fixed and periodic points
(e.g., no sources, sinks, or limit cycles). 

(iii) In the setting of Proposition \ref{lem21}, note that
$$
\Phi_{-f} (t,z) = - \Phi_f (-t , -z) \quad \forall (t,z) \in \R \times
\C \, .
$$
As far as the dynamics of $\Phi_f$ is concerned, therefore, $f$ may be
replaced by $-f$ whenever convenient, e.g., if $f(s) \ne 0$ or
$f'(s) \ne 0$ for all $s\in \R^+$ then it may be assumed that $f>0$ or that $f$ is
increasing, respectively.

(iv) With $\Phi_f(t,\infty):= \infty$ for all $t\in \R$, the flow
$\Phi_f$ may be considered a flow on the compactified complex plane
$\C \cup \{\infty\}$. In this setting, the fixed point $\infty$ cannot
be a saddle, but may be a centre. In perfect analogy to
Proposition \ref{prop22}, it is straightforward to show that $\infty$
is a center of $\Phi_f$ if and only if $\F_f$ is bounded, and for
every $a\in \R$ there exists an increasing sequence $(s_n)$ with
$\lim_{n\to \infty} s_n = \infty$ such that $|F(s_n) + s_n + a|>s_n$
for all $n$; see Proposition \ref{prop26} below.

(v) Smoothness of $f:\R^+ \to
\R$ is assumed throughout for convenience only. All results remain valid under appropriate finite differentiability
assumptions; in most instances it suffices to assume the function $f$
to be $C^1$. Also note that $\Phi_f$ is not in general a smooth
flow, due to $f(|z|)$ being non-smooth or indeed undefined at
$z=0$. However, if for instance $f$ is an even polynomial, as it is, e.g., in
(\ref{eq1_6}), then clearly $\Phi_f$ is smooth (on $\R \times \C$).
\end{rem}

Given $z\in \mbox{\rm Per}\, \Phi_f$, say that the orbit $\Phi_f (\R,
z)$ is {\bf untwisted} if $0$ lies in the exterior of the
closed path $\Phi_f (\, \cdot \, , z)$; otherwise $\Phi_f ( \R, z)$ is {\bf twisted}. Note that $\Phi_f (\R, 0)$,
if at all periodic, is twisted. By contrast, $\Phi_f (\R, z) = \{ z \}$ is
untwisted for every $z\in \F_f$. Given $z\in \mbox{\rm Per}\, \Phi_f \setminus \Phi_f
(\R, 0 )$, define the {\bf net winding} of $z$ as
\begin{equation}\label{eq29}
\omega_f (z) = \pm \frac1{2\pi} \int_0^{T_f(z)} \!\! f (|\Phi_f (t,z)|) \,
{\rm d}t \, ,
\end{equation}
where the plus sign (respectively, minus sign) applies in (\ref{eq29})
when $\Phi_f (\, \cdot \, , z)$ is oriented counter-clockwise
(respectively, clockwise). Net winding plays a key role in later sections. Here only a few basic
properties are recorded. Clearly, $\omega_f$ is constant along
orbits. If $z\in \Phi_f (\R, 0 )$
then the integral in (\ref{eq29}) may or may not exist (in $\R_{\infty}$); 
it does exist, for instance, if $f\ge 0$ or $f' \ge 0$. Notice,
however, that strictly speaking $\omega_f(z)$ is defined only for
$z\in \mbox{\rm Per}\, \Phi_f \setminus \Phi_f (\R, 0 )$. The
following is an immediate consequence of (\ref{eq21}) and (\ref{eq29}).

\begin{prop}\label{prop23a}
Let $f:\R^+ \to \R$ be smooth, and $z\in \mbox{\rm Per}\,
\Phi_f \setminus   \Phi_f (\R, 0 )$. Then, with the same signs
as in {\rm (\ref{eq29})},
\begin{equation}\label{eq29aa}
\omega_f (z) - k_z = 
\pm \frac1{2\pi} \int_0^{T_f(z)}    \frac{{\rm d}t}{\Phi_f (t,z)} =
\pm \frac1{2\pi} \int_0^{T_f(z)}   \frac{\rr \Phi_f (t,z)}{|\Phi_f
  (t,z)|^2} \, {\rm d}t  \, ,
\end{equation}
where $k_z = 0$ or $k_z = 1$ when $\Phi_f (\R, z)$ is untwisted
or twisted, respectively.
\end{prop}

On the (possibly empty or disconnected) set $\mbox{\rm Per}\, \Phi_f \setminus ( \Phi_f(\R,
0) \cup \mbox{\rm Fix}\, \Phi_f)$, the function $\omega_f$ is
continuous, but it is not in general continuous at $z\in \mbox{\rm Fix}\,
\Phi_f$ as, for instance, $\omega_f(z) = 0\ne \lim_{p\to z} \omega_f
(p)$ for $z\in \F_f$, provided that $z$ is a non-degenerate center.
(Recall that $\F_f \subset \R$.)

\begin{lem}\label{lem24}
Let $f:\R^+ \to \R$ be smooth. If $z\in \F_f$ and $z
F''(|z|)>0$ then $z$ is a center of $\Phi_f$, and
\begin{equation}\label{eq241}
\lim\nolimits_{p\to z} \omega_f (p) = \frac1{\sqrt{z F''(|z|)}}
\, .
\end{equation}
\end{lem}

\begin{proof}
It is readily seen that $z\in \F_f$ is a non-degenerate maximum or
minimum of $H_f$ if and only if $f(|z|) F''(|z|) = z
F''(|z|)/|z|^2 > 0$. In this case, $z$ is a center, and $f(|\Phi_f
(t,p)|) \to f(|z|)$ uniformly in $t$ as $p\to z$, whereas
$T_f (p) \to 2\pi/\sqrt{ f(|z|) F'' (|z|)}$, the minimal period
of (\ref{eq28}). Consequently,
$$
\lim\nolimits_{p\to z} \omega_f (p) = \pm \frac1{2\pi} f(|z|) \cdot
\frac{2\pi}{ \sqrt{ f (|z|) F''(|z|)}} = \pm \frac{|z|}{z
  \sqrt{z F''(|z|)}} \, ,
$$
and since $\Phi_f (\, \cdot \, , z)$ is oriented counter-clockwise
(respectively, clockwise) when $z>0$ (respectively, $z<0$), this
proves (\ref{eq241}).
\end{proof}

Recall from Remark \ref{rem23}(i) that $0$ is a center of
$\Phi_f$ whenever $\liminf_{s\to 0} s |f(s)|>1$. Under a slightly
stronger assumption, the behaviour of $\omega_f$ near $0$ is
as follows.

\begin{lem}\label{lem25}
Let $f:\R^+ \to \R$ be smooth. If $\lim_{s\to 0} s |f(s)| =
a > 1$ then $0$ is a center of $\Phi_f$, and
\begin{equation}\label{eq251}
\lim\nolimits_{z\to 0} \omega_f(z) = \frac1{\sqrt{1 - 1/a^{2}}} \, .
\end{equation}
\end{lem}

\begin{proof}
By Proposition \ref{prop22}, $0$ is a center, and clearly $s_0 \DD
\setminus \{ 0\} \subset \mbox{\rm Per}\, \Phi_f \setminus
\mbox{\rm Fix}\, \Phi_f$ for some $s_0 \in \R^+$. To establish
(\ref{eq251}), assume first that $a= \infty$, and in fact $\lim_{s\to
\infty} s f(s) = \infty$. Given any $b\in \R^+$, it can be assumed
that $|\Phi_f (t,z)| f (|\Phi_f (t,z)|) \ge b + 1$ for all $0<|z|<s_0$
and all $t\in \R$. Noting that $\Phi_f (\, \cdot \, , z)$ winds around
$0$ counter-clockwise, write $\Phi_f (t,z) = \rho e^{\textstyle i \varphi}$, with 
smooth functions $\rho = \rho(t)>0$ and $\varphi =
\varphi (t)$. With this, (\ref{eq21}) reads
$$
\dot \rho = - \sin \varphi \, , \quad \dot \varphi = f (\rho) -
\frac{\cos \varphi}{\rho} \, .
$$
Note that $\dot \varphi \ge (\rho f (\rho) - 1)/\rho \ge b / \rho >
0$. It follows that
$$
\omega_f (z) = \frac1{2\pi} \int_0^{T_f(z)} \!\! f (\rho) \, {\rm d}t =
\frac1{2\pi} \int_0^{T_f (z)} \left( \dot \varphi + \frac{\cos
    \varphi}{\rho} \right)  {\rm d}t = 1 + \frac1{2\pi}
\int_0^{T_f(z)} \frac{\cos \varphi}{\rho} \, {\rm d}t \, ,
$$
and consequently
$$
|\omega_f (z) - 1| \le \frac1{2\pi} \int_0^{T_f (z)} \frac{{\rm
    d}t}{\rho} \le \frac1{2\pi} \int_0^{T_f (z)} \frac{\dot
  \varphi}{b} \, {\rm d}t = \frac1{b} \, .
$$
Since $b\in \R^+$ has been arbitrary, $\lim_{z\to 0} \omega_f
(z) = 1$. The argument in case $\lim_{s\to 0} s f(s) = - \infty$ is
completely analogous, and hence (\ref{eq251}) is correct when
$a=\infty$.

It remains to consider the case $1<a<\infty$. Assume first that $\lim_{s\to 0} s f(s) = a$. Note that $\lim_{s\to 0}
F(s)$ exists in $\R$, and so does $\lim_{z \to 0} H_f
(z)$. Since $0$ is a center, it suffices to consider $\omega_f
(s )$ for sufficiently small $s\in \R^+$. For every such $s$, there
exists a unique $0< s^* < s$ with $H_f (- s^* ) =
H_f (s ) = - F (s)$, or equivalently $s^* + s =
\int_{s^*}^s u f(u) \, {\rm d}u$. From the latter, it is
easily deduced that $\lim_{s\to 0} s^*/s =(a-1)/(a+1)$. By
the symmetry of (\ref{eq21}) and ${\rm d}|z|/{\rm d}t = - \ii z/|z|$,
$$
\omega_f (s ) = \frac1{\pi} \int_0^{\frac12 T_f (s  )} \!\! f (|\Phi_f
(t, s )|) \, {\rm d}t = \frac1{\pi} \int_{s^*}^s
\frac{uf(u)}{y(u)} \, {\rm d}u \, ,
$$
where $y = y(u)\ge 0$ is determined uniquely by $|x+ i y|=u$ and $H_f
(x + i y) = - F (s)$, that is,
\begin{align*}
y^2 = u^2 - x^2 & = u^2 - \bigl( u + F (u) - F (s)\bigr)^2  = \bigl( F (s) - F (u)\bigr) \bigl( 2u + F (u) - F (s)
\bigr) \\
& = (s-u) \left( \frac1{s-u} \int_u^s v f(v) \, {\rm d}v - 1\right) (u
- s^*) \left( \frac1{u - s^*} \int_{s^*}^u v
  f(v) \, {\rm d}v + 1 \right)\, ,
\end{align*}
and consequently
$$
\omega_f (s ) = \frac1{\pi} \int_{s^*/s}^1
\frac{\widehat{f}_s (u)\, {\rm d}u}{\sqrt{(1-u) (u -
    s^*/s)}} \, ,
$$
where $\widehat{f_s}:[0,1] \to \R$ is the continuous function with
$$
\widehat{f}_s (u) = \frac{su f(su)}{\sqrt{\frac1{1 - u} \int_u^1 sv
    f(sv) \, {\rm d}v - 1}   \sqrt{\frac1{u - s^*/s}
    \int_{s^*/s}^u sv f(sv) \, {\rm d}v + 1}}  \quad \forall 0< u < 1,
u\ne \frac{s^*}{ s}  \, .
$$
Notice that $\lim_{s\to 0}\widehat{ f}_s(u) = a/\sqrt{a^2 - 1}$ uniformly
on $[0,1]$, and hence
$$
\lim\nolimits_{s\to 0} \omega_f (s ) = \frac1{\pi}
\int_{\frac{a-1}{a+1}}^1 \frac{a}{\sqrt{a^2 - 1}} \cdot \frac{{\rm
    d}u}{\sqrt{(1-u) \bigl( u - {\textstyle \frac{a-1}{a+1}} \bigr)}} =
\frac1{\sqrt{1 - 1/a^{2}}} \, .
$$
This establishes (\ref{eq251}) when $\lim_{s\to 0} s f(s) = a$, and
again the case $\lim_{s\to 0} s f(s) = -a$ is completely analogous.
\end{proof}

In order for $\omega_f (z)$ to be defined whenever $|z|$ is large, note that
$z\in \mbox{\rm Per}\, \Phi_f \setminus \mbox{\rm Fix}\, \Phi_f$ for
all sufficiently large $|z|$, provided that $\liminf_{s\to \infty} s
|f(s)|> 1$. In the terminology of Remark \ref{rem23}(iv), the fixed point
$\infty$ is a center of $\Phi_f$ in this case. The following, then, is an analogue of Lemma \ref{lem25}; its
very similar proof is left to the interested reader.

\begin{prop}\label{prop26}
Let $f:\R^+ \to \R$ be smooth. If $\lim_{s\to \infty} s
|f(s)| = a> 1$ then every solution of {\rm (\ref{eq21})} is bounded,
$\C \setminus s \DD \subset \mbox{\rm Per}\, \Phi_f \setminus
\mbox{\rm Fix}\, \Phi_f$ for some $s\in \R^+$, and
$$
\lim\nolimits_{z \to \infty} \omega_f (z) = \frac1{\sqrt{1 - 1/a^{2}}} \, .
$$
\end{prop}

\begin{rem}\label{rem26a}
By its very definition (\ref{eq29}), the function $\omega_f$ bears some resemblance to the minimal period
function $T_f$. The literature on minimal
periods in Hamiltonian systems, notably near non-degenerate centers,
is substantial; see, e.g., \cite{benguria, chicone, CD, CW, cima1, cima2,
  Rothe, villa, walter, ZP} and the many references
therein. The author does not know whether these fine studies
can fruitfully be applied for the purpose of the present article, and
in particular whether a multiple of $\omega_f$ can be interpreted as the {\em
  true\/} period in a $1$-DOF Hamiltonian flow. Usage of $\omega_f$ in later sections may also remind the
reader of the basic differential geometry notions of total curvature
and rotation index \cite{BG, doCarmo,klingenberg, kuehnel}. Unlike the latter, however,
the value of $\omega_f$ need not be an integer but can in fact be any
(extended) real number.
\end{rem}

For the analysis in later sections, it is crucial whether or not
$\omega_f$ attains certain particular values. To state a simple first
observation in this regard, note that if $f(s) \ne 0$ or $f'(s) \ne 0$
for all $s\in \R^+$ then $\omega_f (z)$ is well-defined (in
$\R_{\infty}$) for every $z\in \mbox{\rm Per}\, \Phi_f$, unless
$z = 0 \in \mbox{\rm Fix}\, \Phi_f$, in which case simply define
$\omega_f (0) = 0$. With this, $\omega_f (z) = 0$ for every
$z\in \mbox{\rm Fix}\, \Phi_f$. Moreover, the possible values of $\omega_f$
always are constrained as follows.

\begin{lem}\label{lem27}
Let $f:\R^+ \to \R$ be smooth, and $z\in \mbox{\rm Per}\,
\Phi_f\setminus \mbox{\rm Fix}\, \Phi_f$.
\begin{enumerate}
\item If $f(s) \ne 0$ for all $s\in \R^+$ then $\omega_f (z) >0$.
\item If $f'(s) \ne 0$ for all $s\in \R^+$ then $\omega_f (z) \ne 1$.
\item If $f(s) f'(s) \ne 0$ for all $s\in \R^+$ then $0<\omega_f (z) <
  1$ when $ff'>0$, and $\omega_f (z) >1$ when $ff' < 0$.
\end{enumerate}
\end{lem}

\begin{proof}
To see (i), simply note that $\Phi_f (\, \cdot \, , z)$ is positively
(respectively, negatively) oriented when $f>0$ (respectively, $f<0$),
and hence $\omega_f(z)>0$ in either case, by (\ref{eq29}).

To prove (ii), let $s_1 < s_2$ be the 
intersection points of $\Phi_f (\R, z)$ with the real axis. For
convenience, let $a = H_f(z) = H_f (s_1) = H_f
(s_2)$. By Remark \ref{rem23}(iii), it may be assumed that
$f'>0$. Let $s_0 = \inf \{ f>0 \}$. If $\max \{|s_1|,|s_2|\}\le s_0$
then, utilizing (\ref{eq21}) and its symmetry,
$$
\omega_f (z) = \omega_f (s_1 ) = - \frac1{\pi} \int_0^{\frac12 T_f
  (s_1)} \!\!  f (|\Phi_f (t, s_1)|) \, {\rm d}t  = \frac1{\pi}
\int_{s_1}^{s_2} \!\! \frac{{\rm d}s}{y (s)} \, ,
$$
where $y = y (s) > 0$ for $s_1 < s < s_2$ is given implicitly
by $H_f (s + i y  ) = a$. Since $f$ is increasing, it is readily seen
that the relevant component of the level set $\{ H_f = a\}$ intersects
the set $\bigl\{ z\in \C : (s_1 + s_2) \rr z  = |z|^2 + s_1 s_2 \bigr\}$, i.e.,
the circle with radius $\frac12 (s_2 - s_1)$ centered at $\frac12 (s_2
+ s_1)$, only in $s_1, s_2$, or, more algebraically, $y
(s) \ne \sqrt{(s_2 - s) (s-s_1)}$ for all $s_1 < s < s_2$. With this,
\begin{equation}\label{eq271}
\omega_f (z) - 1  = \frac1{\pi} \int_{s_1}^{s_2} \!\! \frac{{\rm d}s}{y
  (s)} - \frac1{\pi} \int_{s_1}^{s_2} \!\! \frac{{\rm d}s}{\sqrt{(s_2 - s)
    (s- s_1)}} 
= \frac1{\pi} \int_{s_1}^{s_2} \!\! \frac{\sqrt{(s_2 - s)
    (s- s_1)}  - y (s)}{y (s) \sqrt{(s_2 - s) (s- s_1)}} \, {\rm
  d}s \ne 0 \, .
\end{equation}
Since a virtually identical argument applies when $\min \{|s_1|,
|s_2|\}\ge s_0$, it only remains to consider the case $\min \{|s_1|, |s_2|\}< s_0 < \max \{
|s_1| , |s_2| \}$. Thus assume for
instance that $|s_1|> s_0 > |s_2|$. There exists a unique $s_2 < s_3 <
s_0$ such that $\rr \dot z =0$ when $\rr z  = s_3$, and consequently
$$
\omega_f (z)  = - \frac{1}{\pi} \!\! \int_{s_2}^{s_3}  \frac{{\rm
    d}s}{y^-(s)} -  \frac1{\pi} \int_{s_3}^{s_1} \!\!  \frac{{\rm
    d}s}{y^+ (s)} = \frac1{\pi} \int_{s_1}^{s_2} \!\! \frac{{\rm
    d}s}{y^+(s)} + \frac1{\pi} \int_{s_2}^{s_3} \!\! \left(
  \frac1{y^+(s)} - \frac1{y^-(s)}\right) {\rm d} s 
 < \frac1{\pi} \int_{s_1}^{s_2} \!\! \frac{{\rm
    d}s}{y^+(s)}  \, ,
$$
where $0 < y^-(s) < y^+(s)$ for $s_2 < s < s_3$ are the
two solutions of $H_f(s+ i y) = a$. Since $y^+ > \sqrt{(s_2 - s)
  (s - s_1)}$ for all $s_1 < s < s_2$, similarly to (\ref{eq271}),
$$
\omega_f (z) -1 < \frac1{\pi} \int_{s_1}^{s_2} \!\! \frac{{\rm d}s}{y^+
  (s)} - 1 = \frac1{\pi} \int_{s_1}^{s_2} \frac{\sqrt{(s_2 - s)
    (s- s_1)}  - y^+ (s)}{y^+ (s) \sqrt{(s_2 - s) (s- s_1)}} \, {\rm
  d}s < 0 \, .
$$
A completely analogous argument applies when $|s_1|< s_0 < |s_2|$.

To prove (iii), again assume w.l.o.g.\ that $f'>0$. Using the same
quantities as in the proof of (ii), it is readily checked that $f>0$
implies $y (s) > \sqrt{(s_2 - s) (s- s_1)}$ for all $s_1 < s < s_2$,
whereas this inequality is reversed when $f<0$. By (\ref{eq271}),
therefore, $\omega_f (z) < 1$ when $f>0$, and $\omega_f (z) >1$ when $f<0$.
\end{proof}

The following examples illustrate the notions
introduced in this section; in particular, they show how $0$ may
be non-periodic, a periodic but non-fixed point, or a fixed point of
$\Phi_f$, respectively.

\begin{figure}[ht]
\psfrag{tx1}[]{\small $\rr z$}
\psfrag{tx2}[]{\small $\ii z$}
\psfrag{th1}[]{\small $1$}
\psfrag{tee1}[]{\small $1$}
\psfrag{tv1}[]{\small $i$}
\psfrag{tvm1}[]{\small $-i$}
\psfrag{torig}[]{\small $0$}
\psfrag{tf1}[r]{\small $f(s)=1/(s+1)$}
\psfrag{tper1}[r]{\small $\mbox{\rm Per} \, \Phi_f = \varnothing$}
\psfrag{tf2}[r]{\small $f(s)=s$}
\psfrag{tper2}[r]{\small $\mbox{\rm Per} \, \Phi_f = \C$}
\psfrag{tfix2}[l]{\small $\mbox{\rm Fix} \, \Phi_f = \{1 \}$}
\psfrag{twis}[l]{\small twisted}
\psfrag{tuntwis}[l]{\small untwisted}
\psfrag{tom1}[]{\small $\omega_f \to \frac12 \sqrt{2}$}
\psfrag{tom2}[]{\small $\omega_f =\frac34 $}
\psfrag{tom3}[r]{\small $\omega_f \to 1$}
\begin{center}
\includegraphics{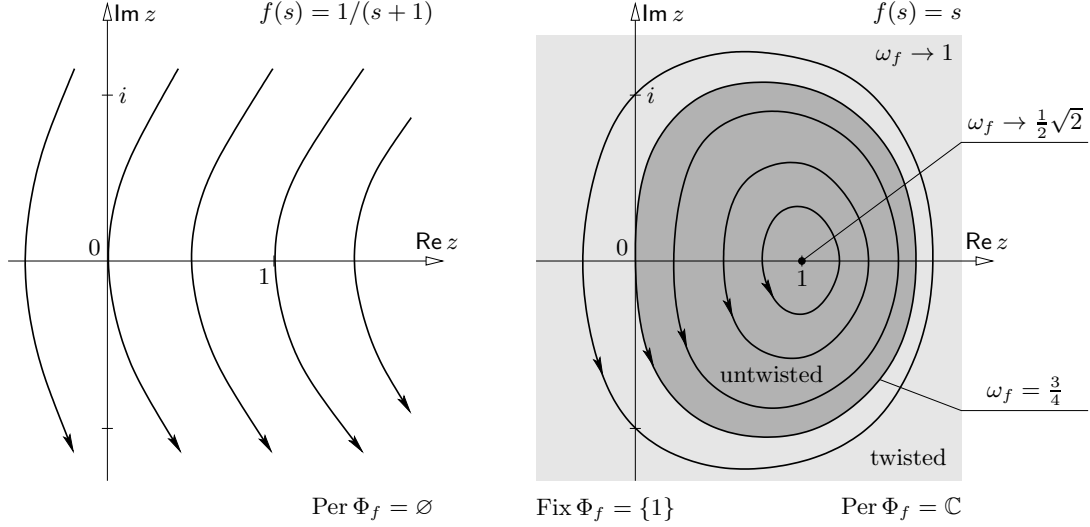}
\end{center}
\caption{For $f(s) = 1/(1+s)$, no point is periodic under $\Phi_f$
  (left; see Example \ref{exa28}), whereas for $f(s) = s$ every point is
periodic, with the center $1$ being the only fixed point (right; see
Example \ref{exa29}).}\label{fig1}
\end{figure}

\begin{example}\label{exa28}
Let $f(s) = 1/(1+s)$ for all $s\in \R^+$. Then $\mbox{\rm Per}\,
\Phi_f = \varnothing$, and every solution of (\ref{eq21}) is
unbounded. In particular, $0$ is non-periodic, its orbit given
implicitly by $\rr z = |z| - \log (1 + |z|)$;
see Figure \ref{fig1}.
\end{example}

\begin{example}\label{exa29}
Let $f(s) = s$ for all $s\in \R^+$. Then $\mbox{\rm Fix}\, \Phi_f =
\F_f = \{1 \}$, and $1$ is a center with $F''(1) =2$, so
$\lim_{z\to 1} \omega_f (z) = \frac12 \sqrt{2}$. Every 
orbit is periodic, and correspondingly every solution of
(\ref{eq21}) is bounded. In particular, $0$ is periodic with
$T_f (0) =  \frac14 \sqrt{6} \, \Gamma (\frac14)^2/\sqrt{\pi}  = 4.541$, its orbit
given implicitly by $3\, \rr z = |z|^3$, and $\omega_f (0) =
\frac34$. Proposition \ref{prop26}  and Lemma \ref{lem27} imply that
$\lim_{z\to \infty} \omega_f (z) = 1$, and 
$0< \omega_f(z) < 1$ for every $z \in \C \setminus \{1\}$; see Figure \ref{fig1}.
\end{example}

\begin{figure}[!ht]
\psfrag{tini}[l]{(i)}
\psfrag{tinii}[l]{(ii)}
\psfrag{tiniii}[]{(iii)}
\psfrag{tpp}[]{\small $p$}
\psfrag{tx1}[]{\small $\rr z$}
\psfrag{tx2}[]{\small $\ii z$}
\psfrag{th1}[]{\small $1$}
\psfrag{tee1}[]{\small $1$}
\psfrag{tv1}[]{\small $i$}
\psfrag{tvm1}[]{\small $-i$}
\psfrag{torig}[]{\small $0$}
\psfrag{tf1}[r]{\small $f(s)=1/s^2$}
\psfrag{tper1}[]{\small $\mbox{\rm Per} \, \Phi_f $}
\psfrag{tfix1}[l]{\small $\mbox{\rm Fix} \, \Phi_f = \{0, 1 \}$}
\psfrag{tper3}[r]{\small $\mbox{\rm Per} \, \Phi_f \! = \C \setminus
 \! \Phi_f (\R, p) $}
\psfrag{tfix3}[l]{\small $\mbox{\rm Fix} \, \Phi_f \! = \{0, 1 \}$}
\psfrag{tf2}[r]{\small $f(s)=1/s$}
\psfrag{tf3}[r]{\small $f(s) = (1+3s^4)/(s+3s^3)$}
\psfrag{tper2}[l]{\small $\mbox{\rm Fix} \, \Phi_f  = \mbox{\rm Per} \,
  \Phi_f = [ 0,\infty [ $}
\psfrag{twis}[l]{\small twisted}
\psfrag{tuntwis}[l]{\small untwisted}
\psfrag{tom1}[]{\small $\omega_f \to 1$}
\psfrag{tominf}[]{\small $\omega_f \to \infty$}
\psfrag{tom3}[r]{\small $\omega_f \to 1$}
\psfrag{tom1a}[]{\small $\omega_f \to \frac13 \sqrt{6}$}
\psfrag{tom2a}[]{\small $\omega_f = 1$}
\psfrag{tom3a}[]{\small $\omega_f \to \infty$}
\begin{center}
\includegraphics{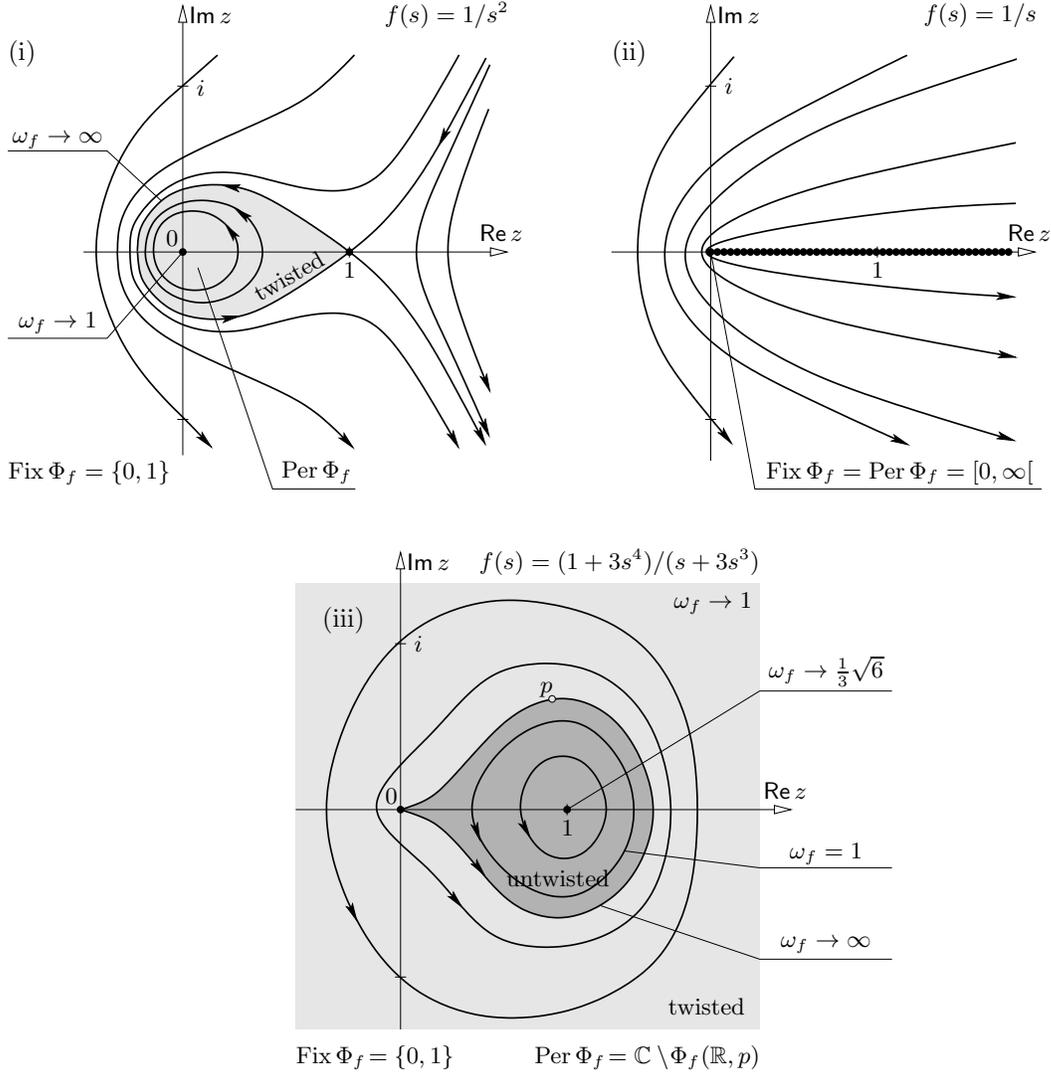}
\end{center}
\caption{When $0$ is a fixed point of $\Phi_f$, it may be a
  center (i), non-isolated (ii), or isolated (iii); see
  Example \ref{exa210}.}\label{fig2}
\end{figure}

\begin{example}\label{exa210}
This example illustrates several different ways how $0$ may be a fixed point of
$\Phi_f$; see Figure \ref{fig2}.

(i) Let $f(s) = 1/s^{2}$ for all $s\in \R^+$. Then $\F_f= \{1\}$, and
$F''(1) = -1$, so the fixed point $1$ is a saddle. By Proposition
\ref{prop22}, the fixed point $0$ is a center, and $\lim_{z \to 0 } \omega_f (z) =1$, 
by Lemma \ref{lem25}. All other periodic
orbits are twisted; they lie inside the homoclinic loop associated with
the saddle and given implicitly by $\rr z = 1 + \log |z|$. By Lemma \ref{lem27},
$\omega_f(z)> 1$ for every $z\in \mbox{\rm Per}\, \Phi_f \setminus \{
0, 1 \}$, and $\omega_f (z) \to \infty$ as $z$ approaches the
homoclinic loop.

(ii) Let $f(s) =1/ s$ for all $s\in \R^+$. Then $\mbox{\rm Fix}\,
\Phi_f = \mbox{\rm Per}\, \Phi_f = \{z\in \C : \rr z \ge 0, \ii z =
0\}$, hence the fixed point $0$ is not isolated. Every point not
on the non-negative real axis has an unbounded orbit which in fact is a parabola.

(iii) Let $f(s) = (1+3s^4)/(s+3s^3)$ for all $s\in
\R^+$. Again, $\F_f = \{1 \}$, but $F''(1) = \frac32$, so unlike
in (i), the fixed point $1$ now is a center. Every point not on the homoclinic
loop associated with the fixed point $0$ and given implicitly by
$3\, \rr z = |z|^3 - |z| + \frac43 \sqrt{3} \arctan (|z| \sqrt{3})$ is periodic. Note that $\omega_f (z)
\to \infty$ as $z$ approaches the homoclinic loop. Also, $\lim_{z\to
  1} \omega_f (z) = \frac13 \sqrt{6 }<1 $ and $\lim_{ z\to \infty}
\omega_f (z) = 1$. Since $f'(0.8294) = 0$, Lemma \ref{lem27}(ii,iii) do not
apply. By the intermediate value theorem there exists at least
one untwisted periodic orbit for which $\omega_f = 1$.
\end{example}

\begin{rem}\label{rem211}
(i) Note that $f(s) f'(s) \ne 0$ for all $s\in \R^+$, and hence Lemma
\ref{lem27}(iii) applies in most examples above, the only exception being Example \ref{exa210}(iii).

(ii) If $0\in \mbox{\rm Fix}\, \Phi_f$ then $0$ can be neither a (topological) saddle nor a source
or sink, due to (\ref{eq23}) and the fact that $\Phi_f$ preserves
Lebesgue measure on $\C$. Rather, the fixed point $0$ must be
a center or degenerate, as in Example \ref{exa210}(i) or (ii,iii), respectively.
\end{rem}

\section{Characterizing closed solutions of $\kappa = f(r)$}\label{sec3}

In all that follows, let $\J_c\subset \R$ be a non-empty open interval
with $0\in \J_c$ and $c: \J_c\to \C$ a  smooth path parametrized
by arc length, that is, $|\dot c (t)| =1$ for all $t\in \J_c$. For every $Q\in
{\rm O} (2)$ let $c_Q(t) = Q \circ c ( \epsilon_Q t )$ for all $t\in
\J_{c_Q} = \epsilon_Q  \J_c$. Thus $c_Q$ is either a rotated (if $\epsilon_Q = 1$) or a reflection-reversed (if
$\epsilon_Q= -1$) copy of $c$. 
Two smooth paths $c, \widehat{c}$ parametrized by arc length are {\bf
  equivalent} if $ \J_{\widehat{c}} = a + \J_c$ for some $a\in \R$ and
$\widehat{c}(t) = c(t-a)$ for all $t\in \J_{\widehat{c}}$. Refer to any equivalence
class as an {\bf oriented smooth curve} $\cC$, and let $[\cC]
= c(\J_c)$ where $c\in \cC$; thus $[\cC]\subset \C$ simply is the set
of points parametrized by some (and hence any) $c\in \cC$. Also, for every $Q \in {\rm O}(2)$ let $\cC_Q$
be the equivalence class of $c_Q$ for some $c\in
\cC$. Note that $\widehat{\cC} = \cC_Q$ implies 
$\bigl[\widehat{\cC}\bigr]= Q ([\cC ])$, but the converse is not true in general.

Given any $c:\J_c \to \C$, associate with it a smooth function
$\vartheta_c : \J_c \to \R$ such that $\dot c =
e^{\textstyle i\vartheta_c}$. Clearly, $\vartheta_c$ is determined by $c$ only up
to an additive integer multiple of $2\pi$. Recall that the
curvature of $c$ is $\kappa_c = \dot \vartheta_c$. If
$c,\widehat{c}$ are equivalent then $\kappa_{\widehat{c}}(t) =
\kappa_c (t-a)$ for all $t\in \J_{\widehat{c}}$. It makes sense,
therefore, so say that, given any smooth function $f:\R^+\to \R$, the
oriented smooth curve $\cC$ is a {\bf solution} of
\begin{equation}\label{eq33}
\kappa = f (r)
\end{equation}
if $\kappa_c (t) = f(|c(t)|)$ for all $t\in \J_c$ with $c(t) \ne 0$, where $c$
is some (and hence any) element of $\cC$. A solution $\cC$ of
(\ref{eq33}) is {\bf maximal} if the set $[\cC]$ cannot be enlarged
any further, that is, if $[\cC ] \subset \bigl[ \widehat{\cC} \bigr]$ for any solution
$\widehat{\cC}$ of (\ref{eq33}) necessarily implies that 
$[\cC] = \bigl[ \widehat{\cC} \bigr]$; see, e.g., \cite[Sec.1]{urbas}. Note that $\cC$ is a (maximal)
solution of (\ref{eq33}) if and only if $\cC_Q$ is a (maximal)
solution for every $Q \in {\rm O}(2)$. Also, given any $p \in \C_{\times}$
and $\vartheta \in \R$, there exists a (locally unique) maximal solution $\cC$
of (\ref{eq33}) such that $c(0) = p$ and $\vartheta_c (0) = \vartheta$ for some
$c\in \cC$.

As alluded to already in the Introduction, the main objective of this
article is to systematically study all (maximal) solutions of
(\ref{eq33}). This is accomplished by making these
solutions correspond to the orbits of the planar flow $\Phi_f$ introduced in the
previous section. To establish such a correspondance, notice first that for
a smooth path $c:\J_c \to \C$ parametrized by arc length, (\ref{eq33})
with $p = e^{\textstyle i\vartheta} \in S^1$ simply reads
$$
\dot c = p \, , \quad \dot p  = i  f(|c|) p \, ,
$$
provided that $c\ne 0$. While this may be read as an ODE on
$\C_{\times} \times S^1$, the dimension of the latter phase space can actually be reduced with
very little effort. Specifically, given $c$, associate
with it the smooth path $z_c : \J_c \to \C$ with
\begin{equation}\label{eq31}
z_c(t) = \overline{i c(t) e^{\textstyle -i\vartheta_c (t)}} = -i \,
\overline{c(t)} e^{\textstyle i\vartheta_c(t)}
\quad \forall t \in \J_c  \, .
\end{equation}
Note that $|z_c| = |c|$. Now, assume that the oriented smooth curve
$\cC$ is a solution of (\ref{eq33}), and pick any $c\in
\cC$. Differentiation of (\ref{eq31}) yields, for $t\in \J_c$,
$$
\dot z_c(t) = -i \, \overline{\dot c(t)} e^{\textstyle i\vartheta_c(t)} - i\, 
\overline{c(t)} e^{\textstyle i\vartheta_c(t)} i \kappa_c(t) = -i + i z_c(t) 
f(|z_c(t)|) \, ,
$$
provided that $c(t)\ne 0$. At least on the non-empty open set $\{t \in
\J_c : c(t) \ne 0\}$, therefore, $z_c$ is a solution of (\ref{eq21}),
an ODE on $\C_{\times}$. It is the purpose of this section to
demonstrate that the correspondance $c \leftrightarrow z_c$ indeed
enables the systematic study of (\ref{eq33}) by way of $\Phi_f$. A
first simple observation in this regard is that most solutions of
(\ref{eq33}) can be reconstructed from the corresponding orbit of
$\Phi_f$; the routine proof is left to the interested reader, as are
the proofs of several equally elementary observations below.

\begin{prop}\label{prop31a}
Let $f:\R^+ \to \R$ be smooth. For every $p \in \C \setminus
\Phi_f(\R, 0)$ and $\vartheta\in \R$, the smooth path $c:\R \to
\C_{\times}$ given by
$$
c(t) = |p| e^{\textstyle i\vartheta + i \int_0^t {\rm d}u/ \overline{\Phi_f (u,p)}}
\quad \forall t \in \R \, ,
$$
is parametrized by arc length, satisfies $\kappa_c = f(|c|)$, and
$z_c(t) = \Phi_f(t,p)$ for all $t \in \R$.
\end{prop}

\begin{rem}\label{rem31b}
Arguably, (\ref{eq31}) might be more natural still if
it were made the definition of $\overline{z_c}$ rather than of
$z_c$, a modification that would not affect the substance of the
subsequent analysis. However, the specific form of (\ref{eq31}) has
been chosen in order to ensure that, in all relevant situations, the paths $c$ and $z_c$ have the same
orientation; see, for instance, the proof of Lemma \ref{lem37} below.
\end{rem}

Another basic observation is that replacing $c$ with
$\widehat{c} _Q$, where $\widehat{c}$ is equivalent to $c$ and $Q \in
{\rm O}(2)$, affects $z_c$ only in a trivial way.

\begin{prop}\label{prop31}
Let $\cC$ be an oriented smooth curve, and $c,\widehat{c}\in
\cC$. For every $Q \in {\rm O}(2)$ there exists $a\in \R$ such that
$$
z_{\widehat{c}_Q} (t) = \left\{
\begin{array}{ll}
z_c(t-a) & \mbox{if } \epsilon_Q = 1 \, , \\[3mm]
\overline{z_c (a-t)} & \mbox{if } \epsilon_Q = -1  \, , \\
\end{array}
\right.
$$
for all $t\in \J_{\widehat{c}_Q} = a + \epsilon_Q \J$.
\end{prop}

With Propositions \ref{lem21} and
\ref{prop31}, every solution $\cC$ of (\ref{eq33})
corresponds to a uniquely determined $\Phi_f$-orbit, namely $\Phi_f
\bigl( \R, z_c (0) \bigr)$, and it is easily seen that every $\Phi_f$-orbit, with the
possible exception of $\Phi_f(\R, 0)$, can be
obtained that way. Moreover, by Proposition \ref{prop31} a
$\Phi_f$-conjugate orbit is obtained if $\cC$ is replaced by $\cC_Q$
for any $Q\in {\rm O}(2)$. Leaving aside trivial exceptions, a stronger statement can be made that has the additional
benefit of being reversible.

\begin{prop}\label{prop32}
Let $f:\R^+ \to \R$ be smooth with $\sup \{f \ne 0 \}= \infty$. For every two maximal solutions $\cC, \widehat
{\cC}$ of {\rm (\ref{eq33})} the following are equivalent:
\begin{enumerate}
\item $[\cC], \bigl[\widehat{\cC}\bigr]$ are ${\rm O}(2)$-congruent;
\item $z_c (0), z_{\widehat{c}}(0)$ are $\Phi_f$-conjugate for some
  (and hence every) $c\in \cC$, $\widehat{c} \in \widehat{\cC}$.
\end{enumerate}
\end{prop}

Whenever $\{f \ne 0\}$ is unbounded, therefore, Proposition
\ref{prop32} establishes a bijection between the maximal solutions of (\ref{eq33}) modulo rotations and
reflection-reversals on the one hand, and the orbits of $\Phi_f$
modulo $\Phi_f$-conjugacy on the other hand. As a consequence, the study of
maximal solutions of (\ref{eq33}) modulo ${\rm O}(2)$-congruence, the central theme throughout the
remainder of the present article, can proceed mostly via a careful
analysis of the orbits of $\Phi_f$.

\begin{rem}\label{rem33}
Maximality of $\cC, \widehat{\cC}$ is
essential for both implications in Proposition \ref{prop32}. Moreover,
(i)$\Rightarrow$(ii) may fail when $\{f \ne 0\}$ is bounded, whereas
(ii)$\Rightarrow$(i) remains correct in this case also.
\end{rem}

Say that an oriented smooth curve $\cC$ is {\bf closed} if $[\cC]$ is compact,
and {\bf simple closed} if $[\cC]$ is homeomorphic to the unit
circle. Though not without precursors in the literature, e.g., in \cite{kuehnel}, this terminology is
tailor-made for the present article and may seem 
unconventional. For maximal solutions of (\ref{eq33}), though, it is easily seen
to be equivalent to a more conventional notion \cite{BG, doCarmo,
  klingenberg}.

\begin{prop}\label{prop34}
Let $f:\R^+ \to \R$ be smooth. For every maximal solution
$\cC$ of {\rm (\ref{eq33})} the following are equivalent:
\begin{enumerate}
\item there exists an oriented smooth curve $\widehat{\cC}$ and $a\in \R^+$ with $\bigl[ \widehat{\cC}\bigr] = [\cC ]$ and
  $\widehat{c} (t+a) = \widehat{c} (t)$ for some $\widehat{c}
  \in \widehat{\cC}$ and all $t\in \J_{\widehat{c}} = \R$;
\item $\cC$ is closed.
\end{enumerate}
Moreover, $\cC$ is simple closed if and only if $\widehat{c}$ in
{\rm (i)} can be chosen to be one-to-one on $[0,a[$.
\end{prop}

\begin{rem}\label{rem35}
Plainly, (i)$\Rightarrow$(ii) in Proposition \ref{prop34} for every
oriented smooth curve $\cC$. However, (ii)$\Rightarrow$(i) may fail
when $\cC$ is not a maximal solution of (\ref{eq33}).
\end{rem}

Notice that if two closed solutions $\cC , \widehat{\cC}$ of
(\ref{eq33}) are ${\rm O}(2)$-congruent then in fact $[\cC] =
e^{\textstyle i \vartheta} [\widehat{\cC}]$ for some $\vartheta \in \R$.
The main goal of the remainder of this section is to characterize
maximal solutions that are closed, perhaps even simple closed curves. As the
reader may have suspected all along, the planar flow $\Phi_f$ and in
particular the 
net winding $\omega_f$ of its periodic points are instrumental in this characterization.

\begin{lem}\label{lem36}
Let $f:\R^+ \to \R$ be smooth. For every maximal solution
$\cC$ of {\rm (\ref{eq33})} with $\mbox{\rm dist} (0,  [\cC] ) > 0$ the following are equivalent:
\begin{enumerate}
\item $\cC$ is closed;
\item $z_c (0)\in \mbox{\rm Per}\, \Phi_f$ and $\omega_f \bigl( z_c
  (0)\bigr)\in \Q$ for some (and hence every) $c\in \cC$.
\end{enumerate}
\end{lem}

\begin{proof}
To prove (i)$\Rightarrow$(ii), let $\cC$ be a closed maximal solution
of (\ref{eq33}), and $c\in \cC$. By Propositions \ref{prop32} and \ref{prop34}, it can
be assumed that $c(t+a) = c(t)$ for some $a\in \R^+$ and all $t\in \J_c =
\R$. Differentiation yields $\vartheta_c (t+a) - \vartheta_c (t) =
2\pi k$ for some $k\in \Z$ and all $t\in \R$. From (\ref{eq31}), it is
clear that $z_c (t+a) = z_c(t)$, and hence $z_c(0) \in \mbox{\rm
  Per}\, \Phi_f$. Since obviously $\omega_f \bigl( z_c
(0)\bigr) = 0 \in \Q$ whenever $z_c (0) \in \mbox{\rm
  Fix}\, \Phi_f$, henceforth assume that $z_c(0) \in \mbox{\rm Per}\,
\Phi_f \setminus \mbox{\rm Fix}\, \Phi_f$, in which case $a =
m T_f \bigl( z_c (0)\bigr)$ for some $m\in \N$. Moreover, since $\kappa_c (t) =
f(|c(t)|)$ for almost all $t$,
\begin{equation}\label{eq36aa}
2\pi k  = \int_0^a \dot \vartheta_c (t) \, {\rm d}t = \int_0^a f(|z_c
(t)|) \, {\rm d}t   = m \int_{0}^{T_f ( z_c (0) ) } \!\!\! f\bigl( \big|\Phi_f\bigl(
t,z_c(0)\bigr) \big|\bigr) \, {\rm d}t = \pm 2 \pi m \omega_f \bigl(
z_c (0)\bigr) \, ,
\end{equation}
and so $\omega_f\bigl( z_c (0)\bigr)\in \Q$, as claimed.

To prove (ii)$\Rightarrow$(i), assume that $z_c (0) \in \mbox{\rm Per}\,
\Phi_f$ and $\omega_f\bigl( z_c (0)\bigr)\in \Q$. Note that $0 \not
\in \Phi_f \bigl( \R, z_c (0)\bigr)$ since otherwise $\mbox{\rm dist}
(0, [\cC]) = 0$. If $z_c (0) \in \mbox{\rm Fix}\,
\Phi_f$ then $[\cC]$ equals the circle with radius $|z_c (0)|>0$
centered at $0$, so clearly $\cC$ is closed. Assume from now on
that $z_c(0) \in \mbox{\rm Per}\, \Phi_f \setminus \mbox{\rm Fix}\,
\Phi_f$, with $b:= \frac12 T_f \bigl( z_c (0)\bigr)>0$ for
convenience. Also, let
$p= z_c (0)$, and pick $\vartheta \in \R$ such that $c(0) = |p|
e^{\textstyle i\vartheta}$. (This is possible because $|c(0)| = |p|$.) By
Proposition \ref{prop31a}, the smooth path $\widehat {c} : \R \to \C$
given by
$$
\widehat{c}(t) = |p| 
e^{\textstyle i\vartheta + i \int_0^t {\rm d}u/\overline{\Phi_f  (u,p)}} \quad \forall t \in \R \, ,
$$
is parametrized by arc length, and $\kappa_{\widehat{c}} = f
(|\widehat{c}|)$. 
Moreover, $\widehat{c} (0) = c(0)$, $ \dot{ \widehat{c} }(0) = \dot c (0)$, 
and consequently $c(\J_c) \subset \widehat{c} (\R)$. By maximality, $[\cC] = \widehat{c} (\R)$, so it
suffices to show that $\widehat{c} (\R)$ is compact. To this end,
simply notice that by (\ref{eq29aa}),
$$
\widehat{c} (t+2b) = \widehat{c}(t) e^{\textstyle i\int_0^{2b} {\rm
    d}u/\overline{\Phi_f (u,p)}} = \widehat{c} (t) e^{\textstyle \pm 2\pi i
  (\omega_f (p) - k_p)} = \widehat{c} (t) e^{\textstyle \pm 2\pi i \omega_f (p)}
\quad \forall t \in \R \, .
$$
Picking $n\in \N$ such that $n\omega_f (p) \in \Z$ yields
$\widehat{c} (t+2nb) = \widehat{c} (t)$ for all $t\in \R$. Thus
$\widehat{c} (\R) = \widehat{c} ([0,2nb])$ indeed is compact, and $\cC$
is closed.
\end{proof}

\begin{rem}\label{rem36a}
As can be seen from the above proof, the assumption $\mbox{\rm dist}
(0, [\cC]) >0$ is not needed for (i)$\Rightarrow$(ii). By
contrast, (ii)$\Rightarrow$(i) may fail without it. Notice, however,
that (i)$\Leftrightarrow$(ii) for {\em every\/} maximal solution $\cC$
of (\ref{eq33}) provided that $f$ can be extended smoothly to $s=0$,
e.g., if $f$ is a polynomial. In
this case, $\cC \leftrightarrow \Phi_f \bigl( \R, z_c (0)\bigr)$ establishes
a bijection between the closed maximal solutions of (\ref{eq33})
modulo rotations and those periodic orbits of
$\Phi_f$ whose net winding is a rational number.
\end{rem}

Let $\cC$ be a closed maximal solution of (\ref{eq33}). By Proposition
\ref{prop31} and Lemma \ref{lem36}, $z_c (0)$ is contained in the same
$\Phi_f$-orbit for every $c\in \cC$. It makes sense, therefore, to
refer to $\cC$ as being (un)twisted whenever $\Phi_f \bigl(\R, z_c (0)
\bigr)$ is (un)twisted, and to let $\omega_f (\cC)= \omega_f \bigl(
z_c (0)\bigr)$ for any $c\in \cC$. Not too surprisingly, if $\cC$ is 
simple closed then the possible values of $\omega_f (\cC)$ are severely
constrained. Henceforth the term {\bf Jordan solution} is used to
refer to any maximal solution of (\ref{eq33}) that is simple closed.
The following observation and its partial converse (Theorem
\ref{thm38} below) are the main results of this section.

\begin{lem}\label{lem37}
Let $f:\R^+ \to \R$ be smooth, and $\cC$ a Jordan solution of {\rm (\ref{eq33})}.
\begin{enumerate}
\item If $\cC$ is untwisted then $\omega_f (\cC) = 0$ or $\omega_f
  (\cC) = 1/n$ for some $n\in \N$.
\item If $\cC$ is twisted then $\omega_f (\cC) = 1$.
\end{enumerate}
\end{lem}

\begin{proof}
Let $\cC$ be a Jordan solution, and $c\in \cC$. By Proposition
\ref{prop34}, it may be assumed that $c(t+a) = c(t)$ for some $a\in \R^+$ and all $t\in \J_c
= \R$, with $c$ being one-to-one on $[0,a[$. If $z_c(0) \in \mbox{\rm
  Fix}\, \Phi_f$ then $\{z_c (0)\} \ne \{ 0 \}$ clearly is
untwisted, and $\omega_f (\cC) = 0$, so henceforth assume that $z_c(0)
\in \mbox{\rm Per}\, \Phi_f \setminus \mbox{\rm Fix}\, \Phi_f$.

It will first be shown that $\Phi_f \bigl(\, \cdot \, , z_c(0) \bigr)$
and $c$ have the same orientation: Either both are oriented
counter-clockwise, or both are oriented clockwise. To see this, let
$s_1 < s_2$ be the intersection points of
$\Phi_f\bigl(\R, z_c(0)\bigr)$ with the real axis. Then $H_f(s_1)
= H_f(s_2) $, and hence $|s_1| \ne |s_2|$. Assume for instance that
$|s_1|< |s_2|$. In this case, $s_2 = \max_{t\in \R} \big| \Phi_f \bigl( t,
z_c(0) \bigr)\big| = \max_{t\in \R} |c(t)|>0$. Since ${\rm d}|z|/{\rm
  d}t = - \ii z/|z|<0$ whenever $\ii z >0$, and $\dot z|_{z = s_2} = i
(s_2 f(s_2) - 1)$, necessarily $s_2f(s_2)>1$. Thus
$f(s_2)>1/s_2>0$, and $\Phi_f \bigl( \, \cdot \, , z_{c}(0) \bigr)$ is
oriented counter-clockwise. Pick $t_0$ with $|c(t_0)| = s_2$. Since
$|c|$ attains a maximum for $t=t_0$, and $\kappa_c (t_0) = f(s_2)>0$,
clearly $c$ is oriented counter-clockwise also. A completely analogous
argument shows that $\Phi_f \bigl(\, \cdot \, , z_c(0) \bigr)$
and $c$ both are oriented clockwise when $|s_1|>|s_2|$.

With $k\in \Z$ and $m\in \N$ as in the proof of Lemma \ref{lem36}, the
theorem of turning tangents (see, e.g., \cite[Thm.5.7.2]{doCarmo})
yields $k = \pm 1$, and (\ref{eq36aa}) simply reads $\pm 2 \pi = \pm 2
\pi m \omega_f (\cC)$, where the plus (respectively, minus) signs
apply when $\Phi_f \bigl(\, \cdot \, , z_c(0) \bigr)$
and $c$ are oriented counter-clockwise (respectively, clockwise). Thus
$\omega_f (\cC) = 1/m$ regardless of whether $\cC$ is untwisted or
twisted. Clearly, this proves (i).

To prove (ii), let $\cC$ be twisted, and assume first that $0
\in [\cC]$. Then $\omega_f (\cC) = \omega_f (0)$ and $0
\in \mbox{\rm Per}\, \Phi_f \setminus \mbox{\rm Fix}\, \Phi_f$. Assume
w.l.o.g.\ that $c(0) = 0$, so $z_c(0) =0$ as well. But $c \bigl( T_f
(0)\bigr)=0$ also, and hence $T_f(0) \ge a = m T_f(0) >0$, yielding
$m\le 1$, that is, $\omega_f (\cC) = 1$.

To complete the proof of (ii), let $\cC$ be twisted but assume that $0\not\in [\cC]$. With the
numbers $s_1 < s_2$ as above, this means that $s_1 < 0 < s_2$. Assume
for instance that $|s_1|< s_2$, and w.l.o.g.\ that $c(0) = s_2$ and
$\vartheta_c (0) = \frac12 \pi$. Then $z_c(0) = s_2$, and $\Phi_f (\,
\cdot \, , s_2 )$ is oriented
counter-clockwise. As in the proof of Lemma \ref{lem25}, write
$\Phi_f (t, s_2) = |\Phi_f(t,s_2)|e^{\textstyle i\varphi(t)} $, where $\varphi$ is smooth and $\varphi(0) = 0$. As seen there,
\begin{equation}\label{eq38a}
\dot \varphi(t)  = f(|z_c (t)|) - \frac{\cos \varphi
  (t)}{|z_c (t)|} =  f(|z_c (t)|) - \frac{\rr z_c(t)}{|z_c(t)|^2}
\quad \forall t \in \R \, .
\end{equation}
Deduce from (\ref{eq31}) that $c(t) = |c(t)|e^{\textstyle i\alpha (t)}$ for all
$t\in \R$, 
with the smooth function $\alpha = \vartheta_c  - \frac12
\pi - \varphi$. Notice that $\alpha (0) = 0$, and by (\ref{eq38a}),
\begin{equation}\label{eq38b}
\dot \alpha (t) = \dot \vartheta_c (t) - \dot \varphi (t)
= \frac{\rr z_c(t)}{|z_c(t)|^2} = \rr \frac1{z_c(t)}
\quad \forall t \in \R \, .
\end{equation}
In particular, $\dot \alpha (0)>0$ and $\dot \alpha (b)<0$, where $b :=
\frac12 T_f (s_2)\le \frac12 a$ for convenience. Also, from the basic
symmetry (\ref{eq24}) it is readily deduced that
\begin{equation}\label{eq38c}
c(-t) = \overline{c(t)} \qquad \mbox{\rm and} \qquad
c(2b-t) = e^{2i\alpha (b)} \overline{c(t)}
\quad \forall t
\in \R \, .
\end{equation}
Now, suppose that $\alpha (t_0) \ge \pi$ or $\alpha (t_0) \le 0$ for some
$0< t_0 < b$. Then $\overline{c(t_0)} = c(t_0)$, and so
$c (-t_0) = c  (t_0)$, by the left equality in
(\ref{eq38c}). Thus $c$ would not be one-to-one on $]-\! b,b[
\: \subset \: ]-\! \frac12 a, \frac12 a[$, contradicting the initial
assumption on $c$. Consequently, $0< \alpha
(t) < \pi$ for all $0<t< b$, and $0\le \alpha (b) \le \pi$. Similarly,
if $\alpha (b)>0$ then $\alpha (t_0) = \alpha (b)$ for some $0<t_0 <
b$, so $e^{\textstyle 2i\alpha (b)} \overline{c(t_0)} = c(t_0)$, and
hence $c (2b - t_0) = c (t_0)$, by the right
equality in (\ref{eq38c}), leading again to the contradictory
conclusion that $c$ is not one-to-one on $]0,2b[\: \subset
\: [0,a[$. In summary, therefore, $\alpha (b) = 0$. Utilizing
(\ref{eq29aa}) and (\ref{eq38b}) yields
$$
0 = \alpha(b) = \int_0^{\frac12 T_f (s_2 )} \dot \alpha (t) \, {\rm d}t
= \int_0^{\frac12 T_f (s_2 )} \frac{ {\rm d}t}{\Phi_f
    (t,s_2)} = \pi (\omega_f (s_2 ) - 1) \, ,
$$
and so $\omega_f (\cC) = \omega_f (s_2 ) = 1$. Since the case
$|s_1|>s_2$ is completely analogous, $\omega_f (\cC) = 1$ for every
twisted Jordan solution $\cC$.
\end{proof}

The following two examples illustrate how reversing the conclusion of Lemma \ref{lem37}
in general may be delicate: Whether or not an untwisted closed maximal
solution $\cC$ of (\ref{eq33}) with $\omega_f (\cC) = 1/n$ for some
$n\in \N$ actually is a Jordan solution may depend on properties of
$f$ not obvious from the outset.

\begin{figure}[ht]
\psfrag{tx1}[]{\small $\rr c$}
\psfrag{tx2}[]{\small $\ii c$}
\psfrag{t1ooeps}[]{\small $1$}
\psfrag{tangle}[]{\small $\pi/6$}
\psfrag{ts0}[]{\small $s_0$}
\psfrag{tvm1}[]{\small $-e_2$}
\psfrag{torig}[]{\small $0$}
\psfrag{tfin}[]{\small $f(s)=-1/a $}
\psfrag{tfout}[]{\small $f(s)= 1$}
\psfrag{tc}[]{\small $\cC$}
\psfrag{tra}[]{\small $a$}
\psfrag{tomc}[]{\small $\omega_f(\cC) = \frac16$}
\begin{center}
\includegraphics{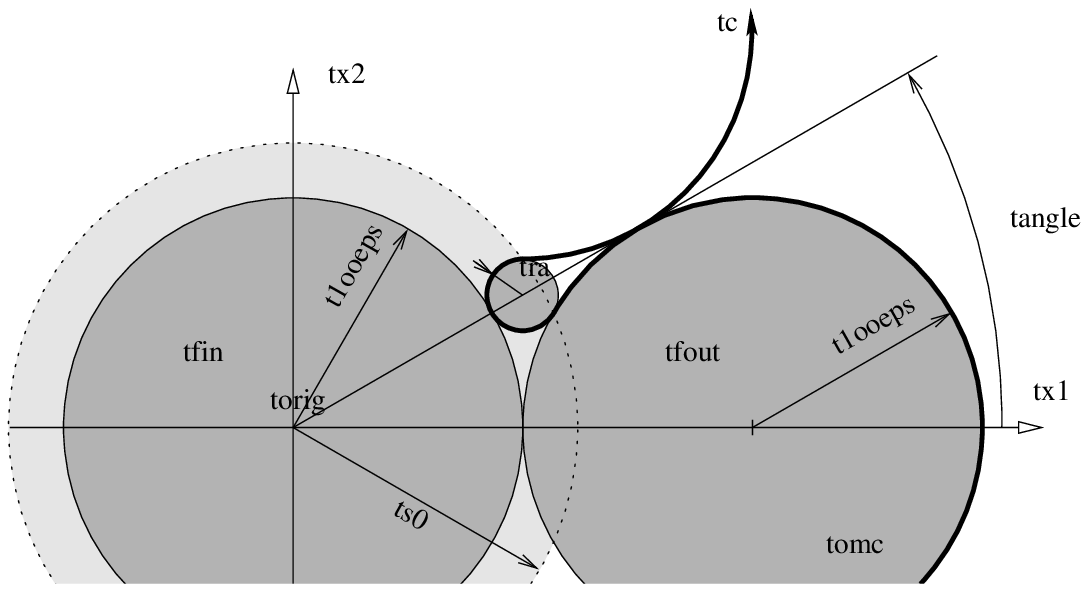}
\end{center}
\vspace*{-18mm}
\caption{A closed maximal solution $\cC$ of (\ref{eq33}) with
  $\omega_f (\cC) = \frac16$ may or may not be a Jordan solution; see
  Example \ref{ex38a}.}\label{fign1}
\end{figure}

\begin{example}\label{ex38a}
Consider three touching discs with radii
$1$, $1$, and $a$, respectively,
positioned as shown in Figure \ref{fign1} (dark grey), where simple trigonometry yields
$a = \frac23 \sqrt{3} - 1 = 0.1547$ and $s_0 = \sqrt{5 - 2\sqrt{3}}
= 1.239$. Let $f(s) = - 1/a$ when $0<s<s_0$ (light grey), and $f(s) = 1$ when
$s>s_0$. The closed solution $\cC$ indicated in Figure \ref{fign1} is
not simple, and yet $\omega_f (\cC) = \frac16$. Of course, $f$ is not
continuous, but given any $\varepsilon >0$, it is
straightforward to construct smooth functions
$f_{\varepsilon}, \widehat{f}_{\varepsilon} : \R^+ \to \R$ with
$\|f_{\varepsilon} - \widehat{f}_{\varepsilon}  \|_{\infty} + \|
f_{\varepsilon} ' - \widehat{f}_{\varepsilon}' \|_{\infty}  <  \varepsilon$, and
corresponding closed maximal solutions $\cC_{\varepsilon} ,
\widehat{\cC}_{\varepsilon}$ with $\omega_{f_{\varepsilon}} (\cC_{\varepsilon}) =
\omega_{\widehat{f}_{\varepsilon}}(\widehat{\cC}_{\varepsilon}) =
\frac16$ such that $\cC_{\varepsilon}$ is simple closed whereas $\widehat{\cC}_{\varepsilon}$ is not.
\end{example}

\begin{figure}[ht]
\psfrag{tp}[]{\small $p$}
\psfrag{tq}[]{\small $q$}
\psfrag{tx1}[]{\small $\rr z$}
\psfrag{tx2}[]{\small $\ii z$}
\psfrag{th1}[]{\small $1$}
\psfrag{tee1}[]{\small $1$}
\psfrag{tv1}[]{\small $i$}
\psfrag{tvm1}[]{\small $-i$}
\psfrag{torig}[]{\small $0$}
\psfrag{tf1}[r]{\small $f(s)=1/(s+1)$}
\psfrag{tf2}[r]{\small $f(s)=3-2/s$}
\psfrag{tper2}[r]{\small $\mbox{\rm Per} \, \Phi_f = \C\setminus
  \Phi_f(\R , \{p,q\})$}
\psfrag{tfix2}[l]{\small $\mbox{\rm Fix} \, \Phi_f = \{-\frac13 ,  0,
  1 \}$}
\psfrag{twis}[r]{\small twisted}
\psfrag{tuntwis}[l]{\small untwisted}
\psfrag{tom1}[]{\small $\omega_f \to \frac13 \sqrt{3}$}
\psfrag{tom2}[]{\small $\omega_f \to - \infty $}
\psfrag{tominf}[]{\small $\omega_f \to - \infty $}
\psfrag{tom3}[]{\small $\omega_f \to \frac23 \sqrt{3}$}
\psfrag{tom4}[r]{\small $\omega_f \to 1$}
\psfrag{tr1}[]{\small $1$}
\psfrag{tr2}[]{\small $2$}
\psfrag{trv1}[]{\small $1$}
\psfrag{trvm1}[r]{\small $-1$}
\psfrag{trv112}[r]{\small $1/2$}
\psfrag{trv1n}[r]{\small $1/n$}
\psfrag{trom}[l]{\small $\omega_f(s)$}
\psfrag{trs}[]{\small $s$}
\psfrag{trc2}[]{\small $\cC_2$}
\psfrag{trcn}[]{\small $\cC_n$}
\psfrag{trfp}[]{\small $f|_{[\cC]}> 0 $}
\begin{center}
\includegraphics{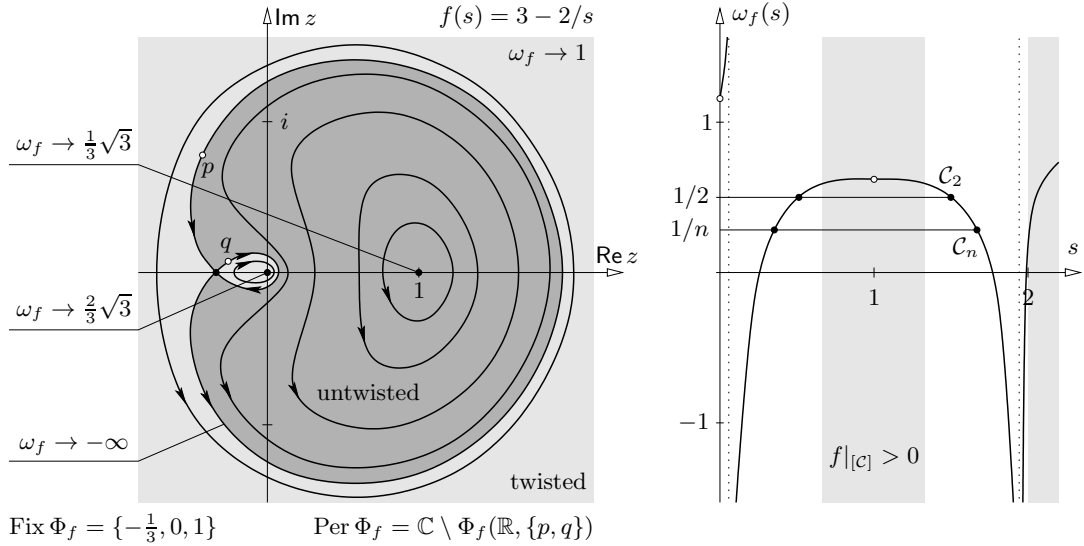}
\end{center}
\caption{For $f(s) = 3 - 2/s$ every point not on the two homoclinic
  loops associated with the saddle $-\frac13$ is periodic, with
  the centers $0$ and $1$ being the only other fixed points (left);
  see Example \ref{ex38b}. Although (\ref{eq33}) has, for
  every $n\in \N \setminus \{1\}$, an untwisted closed maximal
  solution $\cC_n$ with $\omega_f (\cC_n) = 1/n$, only finitely many
  $\cC_n$ are Jordan solutions, and none are convex (right).}\label{fign2}
\end{figure}

\begin{example}\label{ex38b}
Let $f(s) = 3 - 2/s$ for all $s\in \R^+$. Note that $f'(s) = 2/s^2>0$
and $F''(s) = 3$ for all $s\in \R^+$. The flow $\Phi_f$ has exactly
three fixed points: $-\frac13$ which is a saddle with two
associated homoclinic loops given implicitly by $6 \, \rr z = 9 |z|^2 - 12 |z|
+ 1 $, and two centers $0$ and
$1$, with $\lim_{z\to 0} \omega_f (z) = \frac23 \sqrt{3}>1$ and
$\lim_{z\to 1} \omega_f(z) = \frac13 \sqrt{3}<1$. Every point not on the two homoclinic
loops is periodic, and $\lim_{z\to \infty} \omega_f (z) = 1$. Note that $\omega_f(z)\ne 1$ for every $z\in
\mbox{\rm Per}\, \Phi_f$, by Lemma \ref{lem27}(ii). From the phase
portrait of $\Phi_f$ in Figure \ref{fign2}, it is clear that $\omega_f
(z) \to - \infty$ (respectively, $\omega_f(z) \to \infty$) as $z$
approaches a homoclinic loop from within the
untwisted or outer twisted (respectively, inner twisted)
regions. Since $\frac13 \sqrt{3}> \frac12$, for every $n\in \N \setminus \{1\}$ there exists at least one untwisted closed
maximal solution $\cC_n$ of (\ref{eq33}) with $\omega_f (\cC_n) =
1/n$. Numerical evidence strongly suggests that $\omega_f (s )>
\frac12$ for all $\frac23 \le s \le \frac43$, and hence the sign of
$f$ changes along each $\cC_n$. In other words, $\cC_n$ is not
convex for any $n$, by \cite[Thm.2.31]{kuehnel}. Moreover, it is not hard
to see that $\cC_n$ is not a simple closed curve if $n$ is large. In this
example, therefore, (\ref{eq33}) has, modulo rotations,
only a finite number of Jordan solutions, of which only the two
circles with radii $\frac13$ and $1$, centered at $0$ and
oriented clockwise and counter-clockwise, respectively, are convex.
Rigorously determining the precise number of non-circular Jordan
solutions may be a delicate task. 
\end{example}

The above examples make it clear that for the conclusion of Lemma \ref{lem37} to be
reversed in any generality, additional, possibly rather restrictive
assumptions on $f$ have to be imposed. While many assumptions are
conceivable in this regard, one clearly suggesting itself through Lemma
\ref{lem27} is that
\begin{equation}\label{eq3ast}
f(s) f'(s) \ne 0 \quad \forall s \in \R^+ \, .
\end{equation}
For example, the monomial $f_{\delta} (s) = s^{\delta}$, with $\delta
\in \R$, to be studied in detail in Section \ref{sec5} below, satisfies
(\ref{eq3ast}) for every $\delta \ne 0$. Assuming (\ref{eq3ast}), the
conclusion of Lemma \ref{lem37} can be strengthened and reversed rather neatly.

\begin{theorem}\label{thm38}
Let $f:\R^+ \to \R$ be smooth and satisfy {\rm
  (\ref{eq3ast})}. For every closed
maximal solution $\cC$ of {\rm (\ref{eq33})} the following are equivalent:
\begin{enumerate}
\item $\cC$ is untwisted, and $\omega_f (\cC)=0$ or $\omega_f (\cC) = 1/n$ for some $n\in
  \N\setminus \{1\}$;
\item $\cC$ is a Jordan solution.
\end{enumerate}
Moreover, $\omega_f (\cC) =0$ if and only if $\cC$ is a circle of
radius $|z|$ centered at $0$, with $z\in \F_f$, oriented
counter-clockwise when $z >0$, and clockwise when $z < 0$.
\end{theorem}

\begin{proof}
Throughout, let $\cC$ be a closed maximal solution, and $c\in \cC$ and
$c(t+a) = c(t)$ for some $a\in \R^+$ and all $t\in \J_c =
\R$. By Remark \ref{rem23}(iii) assume w.l.o.g.\ that $f>0$.

To prove (i)$\Rightarrow$(ii), recall from Lemma \ref{lem36} that
$z_c(0) \in \mbox{\rm Per}\, \Phi_f$. (This part of the lemma does not
require $0 \not \in [\cC]$; see Remark \ref{rem36a}.) By Lemma
\ref{lem27}, $\omega_f (\cC)>0$ unless $z_c (0) \in
\mbox{\rm Fix}\, \Phi_f$, in which case $z_c (0) \ne 0 $, and $\cC$ 
is a circle with radius $|z_c(0)|$ centered at
$0$, hence obviously a Jordan solution; since
$\kappa_c = f(|z_c(0)|) = 1/z_c(0)$, this circle is oriented
counter-clockwise when $z_c(0)>0$, and clockwise when $z_c(0)< 0$. Thus it
only remains to consider the case $\omega_f (\cC) = 1/n$ where
necessarily $z_c (0) \in \mbox{\rm Per}\, \Phi_f \setminus
\mbox{\rm Fix}\, \Phi_f$. As in the proof of Lemma
\ref{lem37}, let $s_1<s_2$ be the intersection points of
$\Phi_f\bigl(\R , z_c(0)\bigr)$ with the real axis, and assume
w.l.o.g.\ that $z_c(0)=s_2$. Noticing that $s_1s_2>0$ since $\cC$ is
untwisted, assume, for instance, that $0<s_1<s_2$. Since $\rr \dot z_c
= - \ii z_c f(|z_c|)$ is negative (respectively, positive) whenever
$\ii z_c >0$ (respectively, $\ii z_c < 0$), clearly $ \rr \Phi_f(t,s_2) \in
[s_1,s_2]$ for all $t\in \R$. Also, writing $c(t) = |c(t)|
e^{\textstyle i \alpha (t)}$ with the smooth function $\alpha$
satisfying $\alpha (0) = 0$ yields
$$
\dot \alpha (t) = \frac{\ii \bigl( \dot c(t)
  \overline{c(t)}\bigr)}{|c(t)|^2} = \frac{\rr
  \Phi_f(t,s_2)}{|\Phi_f(t,s_2)|^2} = \rr \frac1{\Phi_f(t,s_2)} >0 \quad
\forall t \in \R \, ,
$$
in accordance with (\ref{eq38b}). Thus $\alpha$ is increasing, and
with $b:= \frac12 T_f (s_2)$ for convenience,
$$
\alpha (t + 2b) - \alpha (t) = \rr \int_t^{t + 2b} \frac{{\rm
    d}u}{\Phi_f (u,s_2)} =
\int_0^{T_f (s_2)} \frac{{\rm
    d}u}{\Phi_f (u,s_2)} = 2\pi \omega_f (s_2) = \frac{2\pi}{n} \quad
\forall t \in \R \, ,
$$
as well as $\alpha (b) = \pi/n$. It follows that $c$ is
one-to-one on $[0,2nb[$, but also
$$
c(t+2nb) = |z_c (t+2nb)|e^{\textstyle i \alpha (t+2nb)} = |z_c(t)|
e^{\textstyle i\alpha (t)} = c(t) \quad \forall t \in \R \, .
$$
Thus $\cC$ is a Jordan solution. An analogous argument for
the case $s_1<s_2<0$ completes the proof of (i)$\Rightarrow$(ii).

To see that (ii)$\Rightarrow$(i) simply recall from Lemma \ref{lem27}
that (\ref{eq3ast}) implies $\omega_f (\cC)\ne 1$, and hence the claim
immediately follows from Lemma \ref{lem37}.
\end{proof}

\begin{rem}\label{rem312}
The reader may have noticed that the proof of Lemma \ref{lem37}
presented above does contain a possible, if rather
unpractical {\em necessary and sufficient condition\/} for a closed maximal
solution to be Jordan. To illustrate the simple idea, let
$\cC$ be, for instance, a (non-circular) untwisted closed
maximal solution, with associated numbers $0< s_1<s_2$ as in the proof
of Lemma \ref{lem37}. Define
$\omega_f^* : \: \bigr] 0, \frac12 T_f (s_2) \bigl[ \: \to \R$ by 
$$
\omega_f^* (t) = \pm \frac1{\pi} \, \rr \!\! \int_0^t \frac{{\rm
    d}u}{\Phi_f(u,s_2)} \quad \forall 0<t<
{\textstyle \frac12} T_f (s_2)  \, ,
$$
with the sign as in (\ref{eq29}) and (\ref{eq29aa}). Notice that
$\omega_f^*(0+ ) = 0$ and $\omega_f^*\bigl( \frac12 T_f (s_2) - \bigr) =
\omega_f (\cC)$; in fact, by (\ref{eq38b}) simply $\omega_f^*(t) = \pm
\alpha(t)/\pi$. With this, it is clear that $\cC$ is a Jordan
solution if and only if $\omega_f(\cC)
= 1/n$ for some $n\in \N$, and
\begin{equation}\label{eq3star}
0 < \omega_f^* (t) < \omega_f(\cC) \quad \forall 0< t < {\textstyle
  \frac12} T_f (s_2)   \, .
\end{equation}
Notice that in the setting of Theorem \ref{thm38}, i.e., with $f$
satisfying (\ref{eq3ast}), the function $\omega_f^*$ is increasing, as
seen in the proof, and hence the otherwise unwieldy condition (\ref{eq3star}) holds automatically.
\end{rem}

\section{Jordan solutions and monotone net winding}\label{sec4}

By Lemma \ref{lem36}, closed (maximal) solutions of (\ref{eq33}) are
plentiful unless $\omega_f$ is constant on $\mbox{\rm Per}\,
\Phi_f \setminus \mbox{\rm Fix}\, \Phi_f$. By contrast, Lemma
\ref{lem37} and Example \ref{ex38b} suggest that Jordan solutions
are much rarer. In fact, as detailed in this section and the next,
Jordan solutions of (\ref{eq33}) are exceedingly rare for many
$f$. Though much of the analysis could be carried out, at least
locally, for far more general $f$, assume from now on that the
smooth function $f:\R^+ \to \R$ satisfies (\ref{eq3ast}). This allows
Lemma \ref{lem27} and Theorem \ref{thm38} to be applied together. By
Theorem \ref{thm38}, clearly (\ref{eq33}) has at least as many different
(circular) Jordan solutions as $\F_f$ has elements. For instance,
$\F_f = \R^+ $ for $f(s) = 1/s$, and correspondingly every
(counter-clockwise oriented) circle centered at $0$ is a Jordan
solution; see Example \ref{exa210}(ii). To rule out degenerate situations
like this, and thus to make the ultimate results particularly complete
and transparent, assume in addition that $F''(s) \ne 0$ for all
$s\in \R^+$, so $\F_f$ is either empty or a singleton. For convenience,
then, let
$$
\cF = \bigl\{ f:\R^+ \to \R \enspace \mbox{\rm is smooth with } f(s)
f'(s) F''(s) \ne 0 \: \forall s \in \R^+ \bigr\} \, .
$$
Plainly, $af\in \cF$ for every $a\in \R \setminus \{0\}$ and $f\in \cF$.
Given $f\in \cF$ let $\epsilon_f = \pm 1$ be such that $\epsilon_f
F'' >0$, and notice that the open interval
$$
\I_f = \: \bigl]  \lim\nolimits_{s\to 0} \epsilon_f s f(s) , \lim\nolimits_{s\to
\infty} \epsilon_f s f(s) \bigr [
$$
is well-defined, non-empty, and does not contain $0$. For example,
with $f$ from Examples \ref{exa28}, \ref{exa29},
and \ref{exa210}(ii), clearly $f\in \cF$, and $\I_f$ equals $]0,1[$, $\R^+$, and
$-\R^+$, respectively. Every open interval contained in
$\R \setminus \{0\}$ equals $\I_f$ for an
appropriate $f\in \cF$. To analyze (\ref{eq33}) with $f\in \cF$, it is
convenient to distinguish four cases, depending on the position of
$\I_f$ relative to the two-point set $\{-1,1\}$. Three of the four cases are straightforward,
as recorded in Propositions \ref{prop41a} to \ref{prop42} below.

\begin{prop}\label{prop41a}
Let $f\in \cF$, and assume that $\I_f \cap [-1,1]=\varnothing$. Then
every closed maximal solution of {\rm (\ref{eq33})} is twisted; in
particular, {\rm (\ref{eq33})} has no Jordan solution.
\end{prop}

\begin{figure}[ht]
\psfrag{tx1}[]{\small $\rr z$}
\psfrag{tx2}[]{\small $\ii z$}
\psfrag{th1}[]{\small $1$}
\psfrag{tee1}[]{\small $1$}
\psfrag{tv1}[]{\small $i$}
\psfrag{tvm1}[]{\small $-i$}
\psfrag{torig}[]{\small $0$}
\psfrag{tf1}[l]{\small $f(s)=(s+2)/s^2$}
\psfrag{tf2}[l]{\small $f(s)=(s^2+2)/s^3$}
\psfrag{tfint1}[r]{\small $\I_f = \, ]-\infty, -1[$}
\psfrag{tfint2}[r]{\small $\I_f = \, ]-\infty, -1[$}
\psfrag{tfix1}[l]{\small $\mbox{\rm Fix} \, \Phi_f =\{ 0 \}$}
\psfrag{tfix2}[l]{\small $\mbox{\rm Fix} \, \Phi_f =\{ 0 \}$}
\psfrag{tper1}[r]{\small $\mbox{\rm Per} \, \Phi_f =\C$}
\psfrag{tper2}[r]{\small $\mbox{\rm Per} \, \Phi_f \ne \C$}
\psfrag{twis}[]{\small twisted}
\psfrag{tuntwis}[l]{\small untwisted}
\psfrag{tom2}[]{\small $\omega_f =\frac34 $}
\psfrag{tom1}[]{\small $\omega_f \to 1$}
\begin{center}
\includegraphics{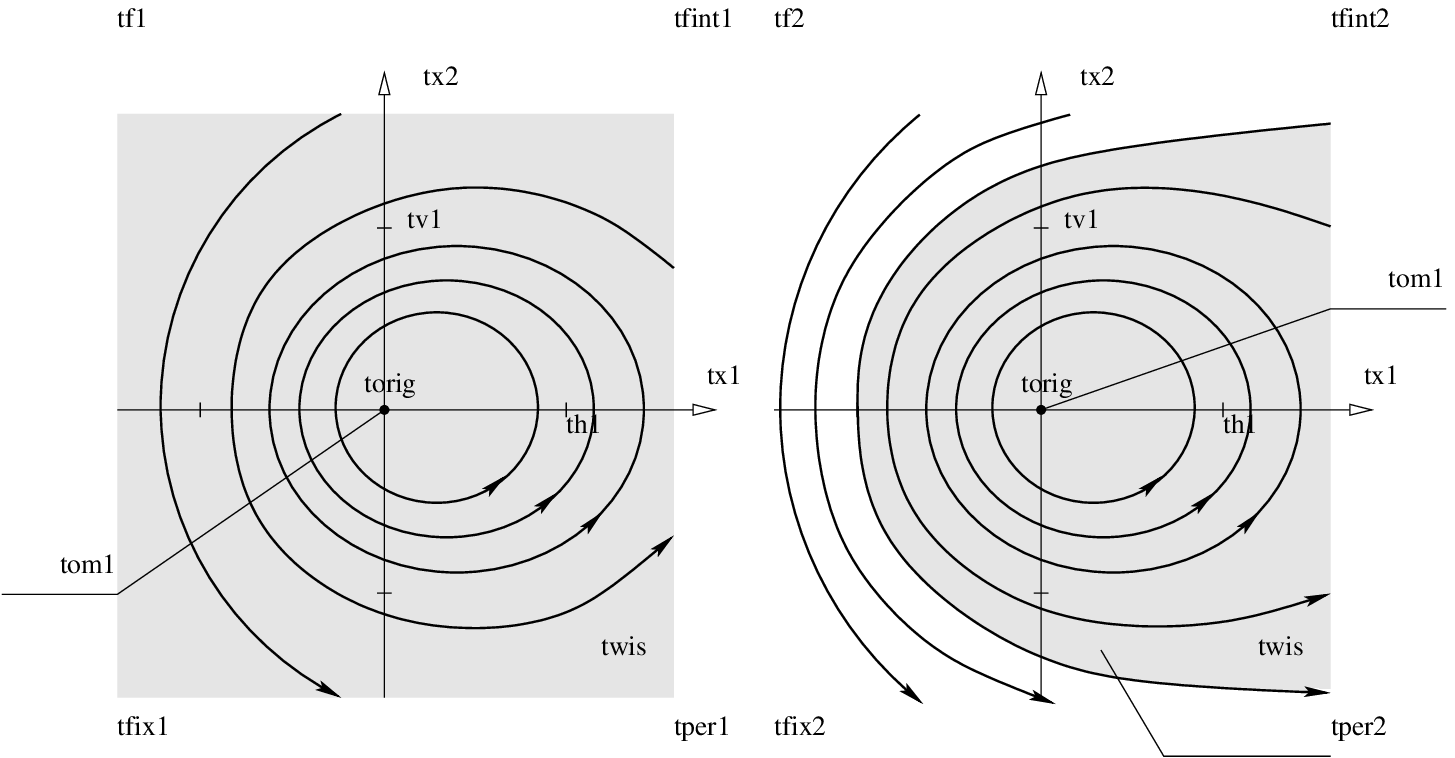}
\end{center}
\caption{If $f\in \cF$ and $\I_f \cap
  [-1,1]= \varnothing $ then every periodic orbit of $\Phi_f$ is twisted
  (left); only when $\sup \I_f = -1$ may non-periodic orbits exist (right).}\label{fig5}
\end{figure}

It is easy to see that with $f$ as in Proposition
\ref{prop41a}, every periodic orbit of $\Phi_f$ is twisted, and the
center $0$ is the only fixed point. Moreover, $\mbox{\rm Per}\,
\Phi_f = \C$, except perhaps when $\sup \I_f = -1$, in which case
$\Phi_f$ may have non-periodic points; see Figure \ref{fig5}.

\begin{prop}\label{prop41b}
Let $f\in \cF$, and assume that $\I_f \subset \: ]-\! 1,1[$. Then
every maximal solution of {\rm (\ref{eq33})} is unbounded; in
particular, {\rm (\ref{eq33})} has no Jordan solution.
\end{prop}

\begin{figure}[ht]
\psfrag{tx1}[]{\small $\rr z$}
\psfrag{tx2}[]{\small $\ii z$}
\psfrag{th1}[]{\small $1$}
\psfrag{tee1}[]{\small $1$}
\psfrag{tv1}[]{\small $i $}
\psfrag{tvm1}[]{\small $-i$}
\psfrag{torig}[]{\small $0$}
\psfrag{tf1}[l]{\small $f(s)=e^{-s}/s$}
\psfrag{tf2}[l]{\small $f(s)=1/(s+s^3)$}
\psfrag{tfint1}[r]{\small $\I_f = \, ]-1, 0[$}
\psfrag{tfint2}[r]{\small $\I_f = \, ]-1, 0 [$}
\psfrag{tfix1}[l]{\small $\mbox{\rm Fix} \, \Phi_f =\{ 0 \}$}
\psfrag{tfix2}[r]{\small $\mbox{\rm Fix} \, \Phi_f = \mbox{\rm Per} \,
  \Phi_f  = \{ 0 \}$}
\psfrag{tper1}[r]{\small $\mbox{\rm Per} \, \Phi_f =\varnothing $}
\psfrag{tper2}[r]{\small $\mbox{\rm Per} \, \Phi_f \ne \R^2$}
\psfrag{twis}[]{\small twisted}
\psfrag{tuntwis}[l]{\small untwisted}
\psfrag{tom2}[]{\small $\omega_f =\frac34 $}
\psfrag{tom1}[]{\small $\omega_f \to 1$}
\begin{center}
\includegraphics{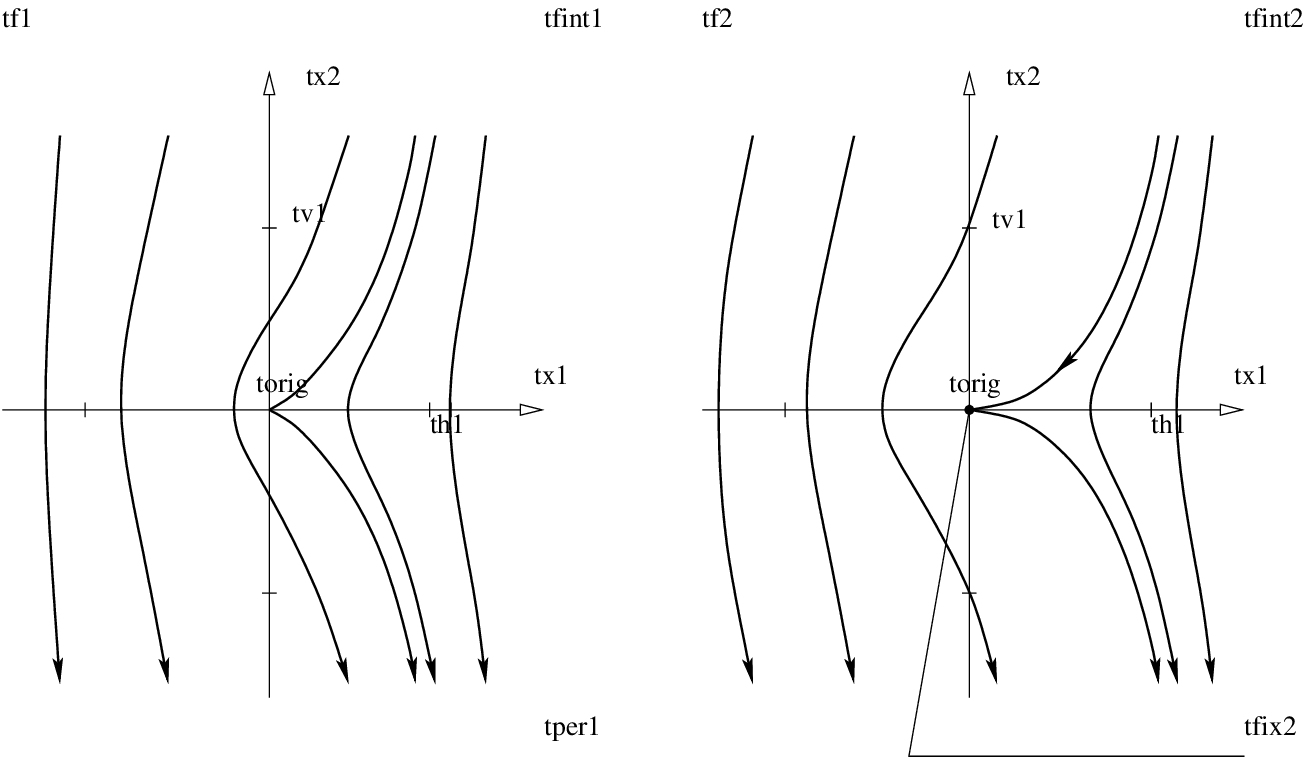}
\end{center}
\caption{If $f\in \cF$ and $\I_f \subset \:
  ]-\! 1,1[$ then $\mbox{\rm Per}\, \Phi_f = \varnothing$ (left),
  except when $\inf
  \I_f = -1$ where $0$ may be a fixed point (right).}\label{fig6}
\end{figure}

Again, with $f$ as in Proposition \ref{prop41b}, it is easy to see
that all orbits of $\Phi_f$ are unbounded, and $\mbox{\rm Per}\,
\Phi_f = \varnothing$, except perhaps when $\inf \I_f = -1$, in which
case $0$ may be a fixed point; see Figure \ref{fig6} and also
Example \ref{exa28}.

\begin{prop}\label{prop42}
Let $f\in \cF$, and assume that $-1 \in \I_f$. Then $\F_f =
\{ s \}$, with $s\in \R \setminus \{0\}$ uniquely determined by $sf(|s|) = 1$.
Every closed maximal solution of {\rm (\ref{eq33})} either is twisted,
or else equals the circle with radius $|s|$ centered
at $0$, oriented counter-clockwise when $s>0$, and clockwise
when $s<0$; in particular, {\rm (\ref{eq33})} has exactly one Jordan
solution (namely, that circle).
\end{prop}

Note that under the assumptions of Proposition \ref{prop42} the fixed
point $s $ is a saddle, its associated homoclinic loop containing
the center $0$ and all periodic orbits, each of which is
twisted; see Example \ref{exa210}(i) for a typical phase portrait.

By Propositions \ref{prop41a} to \ref{prop42},
non-circular Jordan solutions of (\ref{eq33}) with $f\in \cF$, if at
all existant, can be found only in the remaining (fourth) case, that
is, when $1\in \I_f$. Again, for convenience let
$$
\cF^* = \{ f \in \cF : 1 \in \I_f \} \, .
$$
For example, $f_{\delta} \in \cF$ for every $\delta \in \R \setminus
\{-1,0\}$, but $f_{\delta}\in \cF^*$ only when $\delta  > -1$. For
every $f\in \cF^*$, note that $\epsilon_f f > 0$, i.e., $\epsilon_f$
simply is the (constant) sign of $f$. Also, 
$\F_f = \{ \epsilon_f s_f \} = \mbox{\rm Fix}\, \Phi_f$, where $s_f \in \R^+$ is the unique
solution of $sf(s) = \epsilon_f$. The sole fixed point of $\Phi_f$ is
a center, and $\mbox{\rm Per}\, \Phi_f = \C$; see Example \ref{exa29}
for a typical phase portrait.
Observe that every (periodic) orbit not intersecting the line segment
$\epsilon_f \: ]0, s_f]$ is twisted. With Proposition
\ref{prop32} and Theorem \ref{thm38}, therefore, modulo rotations all
maximal solutions of (\ref{eq33}) with $f\in
\cF^*$ that could potentially be Jordan solutions are parametrized by
that segment. More formally, given $f\in \cF^*$ and $0< s\le s_f$, let
$\cC_{f,s}$ be a maximal solution of (\ref{eq33}) with $c(0)
= \epsilon_f s  $ and $\vartheta_c = \frac12 \pi$ for some
$c\in \cC_{f,s}$. For instance, $[\cC_{f,s_f}]$ simply is a circle with
radius $s_f$ centered at $0$. Note that $s\mapsto [\cC_{f,s}]$
is one-to-one, and clearly $[\cC_{f,s}]$ is {\em not\/} a circle when
$s<s_f$. The following properties of the curve $\cC_{f,s}$ are immediate
consequences of Lemma \ref{lem36} and Theorem \ref{thm38}.

\begin{prop}\label{prop43}
Let $f\in \cF^*$ and $0<s\le s_f$.
\begin{enumerate}
\item $\cC_{f,s}$ is closed if and only if $\omega_f (\epsilon_f s )\in \Q$.
\item $\cC_{f,s}$ is simple closed if and only if either $s=s_f$, or
  else $s<s_f$ and $\omega_f (\epsilon_f s ) = 1/n$ for some $n\in
  \N \setminus \{1\}$.
\end{enumerate}
\end{prop}

\noindent
It is clear that, as mentioned earlier and modulo rotations, the family
$(\cC_{f,s})_{0<s\le s_f}$ contains all closed maximal solutions that
are untwisted, so in particular all Jordan solutions.

\begin{prop}\label{prop44}
Let $f\in \cF^*$, and assume that $\cC$ is a closed maximal solution of {\rm
  ({\ref{eq33}})}. If $\cC$ is untwisted then $[\cC] = e^{\textstyle i \vartheta} [\cC_{f,s}]$
for some $\vartheta \in \R$ and a unique $0<s\le s_f$.
\end{prop}

Combining Propositions \ref{prop43} and \ref{prop44}, it is a simple
task to find, at least formally, all Jordan solutions. For
convenience, let
$$
\OO_f = \bigcup\nolimits_{n\ge 2} \left\{0<s<s_f : \omega_f ( \epsilon_f
s )  = \frac1{n} \right\} \, .
$$
Though a next-to-trivial consequence of the above, the following theorem may nevertheless be regarded the
main result of this section as it completely describes all Jordan solutions of
(\ref{eq33}) when $f\in \cF^*$, and hence in fact even when merely $f\in \cF$.

\begin{theorem}\label{thm45}
Let $f\in \cF^*$, and assume that $\cC$ is a closed maximal solution of {\rm
  (\ref{eq33})}. Then $\cC$ is oriented counter-clockwise when
$f>0$, and clockwise when $f<0$. Moreover, the
following are equivalent:
\begin{enumerate}
\item $\cC$ is a Jordan solution;
\item $[\cC]$ either is a circle with radius $s_f$ centered at $0$, or
  else $[\cC] = e^{\textstyle i \vartheta} [\cC_{f,s}]$ for some
  $\vartheta \in \R$ and a unique $s\in \OO_f$.
\end{enumerate}
\end{theorem}

\begin{proof}
Recall that either $f>0$ or $f<0$ since $f\in \cF^*$, and correspondingly
$\epsilon_f = 1$ or $\epsilon_f = -1$. Clearly, $\cC$ is
oriented counter-clockwise in the former case, and clockwise in the
latter \cite[Ch.1]{klingenberg}. That (i)$\Leftrightarrow$(ii) is immediate
from Theorem \ref{thm38} together with Propositions \ref{prop43} and \ref{prop44}.
\end{proof}

For every $f\in \cF^*$, Theorem \ref{thm45} establishes a bijection
between $\OO_f \subset \: \,  ]0,s_f[$ and the non-circular Jordan solutions
of (\ref{eq33}) modulo rotations. Thus, to find all
such solutions one only has to determine the set $\OO_f$, most
basically its cardinality. The remainder of this section aims at determining $\# \OO_{f}$, on the one hand by
establishing a practicable lower bound for $\omega_f$, and on the
other hand by devising a condition that
ensures $ \omega_f $ is monotone. For both
tasks, it is convenient to assume henceforth that $F(s_f) = 0$ for
every $f\in \cF^*$. (The function $F$ has been determined
only up to an additive constant so far.) Note that $|f|\in \cF^*$ whenever
$f\in \cF^*$, and $s_{|f|} = s_f$. For convenience, therefore, assume
that $f>0$ from now on. (If $f<0$ then simply replace $f$ by $-f$ in all
that follows.) Then $F$ is non-negative and (strictly) convex. Let
$F_k = F^{(k)}(s_f)$ for $k\in \N_0$, so in particular $F_0=F_1=0$ and
$F_2>0$. With a view towards Theorem
\ref{thm45}, for the further analysis it is helpful to derive an
explicit formula for $\omega_f$ on $]0,s_f[$ as follows:
For every $0<s<s_f$ there exists a unique
$s^*>s_f$ such that $F (s^*) = F (s)$, or more geometrically, $s, s^*$
are the two intersection points of $\Phi_f
(\R,  s )$ with the real axis. Note that $s\mapsto
s^*$ is smooth and decreasing, with $0_f^*:= \lim_{s\to 0} s^* < 
\infty$ and $\lim_{s\to s_f} s^* = s_f$. With this,
\begin{equation}\label{eq43}
\omega_ f(s  )= \frac{1}{\pi} \int_s^{s^*} \frac{u f(u) \, {\rm
    d}u}{\sqrt{F (s) - F(u) }
\sqrt{2 u + F(u)  - F (s)}}
\quad \forall 0< s < s_f  \, .
\end{equation}
Utilizing (\ref{eq43}), the proof of the following is a straightforward calculus
exercise.

\begin{prop}\label{prop46}
Let $f\in \cF^*$. Then $\omega_f$ is smooth and 
positive on $]0,s_f[$, with
$$
\omega_f (s_f - ) =
\frac1{\sqrt{s_f F_2 }} \, , \quad \omega'_f(s_f-) = 0 \, , \quad
\omega_f''(s_f-)  = \frac{9 F_2^2 - 3 s_f F_2 (3 F_2^2 - 2
  F_3 )  - s_f^2 ( 3 F_2 F_4 - 5F_3^2 )}{24
  \sqrt{ s_f F_2  }^5}   \, . 
$$
\end{prop}

Since $F$ is non-negative and convex, with $F (s_f) =
F' (s_f) = 0$, the ratios $F'' F/(F')^2$ and
$(F')^2/F$ define positive smooth functions on $\R^+$, with
their value for $s=s_f$ equal to $\frac12$ and $2 F_2$,
respectively. The following simple observations are useful when
establishing a lower bound for $\omega_f$.

\begin{prop}\label{lem47}
Let $f\in \cF^*$, and assume that $F''F/(F')^2$ is
increasing (respectively, decreasing) on $\R^+$. Then:
\begin{enumerate}
\item $s\mapsto s^*$ is concave (respectively, convex) on $]0,s_f[$;
\item $(F')^2/F$ is increasing (respectively, decreasing) on $\R^+$;
\item $F'$ is convex (respectively, concave) on $\R^+$.
\end{enumerate}
\end{prop}

\noindent
Using Proposition \ref{lem47}, it is now possible to establish a
reasonably tight lower bound for $\omega_f $ on $]0, s_f[$, provided
that $f$ is {\em increasing}.

\begin{lem}\label{lem48}
Let $f\in \cF^*$ be increasing. If $F''F/(F')^2$ is increasing on $\R^+$ then
\begin{equation}\label{eq48a}
\omega_f (s) >  \frac{\sqrt{s_f}}{\pi} 
\sqrt{\frac{s^* - s_f}{F'(s^*)}} \int_0^{\pi} \!\!  f \! \left( \sqrt{
\frac{ (s^*)^2 + s^2 }{2} - \frac{ (s^*)^2 - s^2}{2} \cos u}
\right) {\rm d}u
\quad \forall  0<s<s_f  \, .
\end{equation}
If $F''F/(F')^2$ is decreasing on $\R^+$ then {\rm
  (\ref{eq48a})} holds with $\displaystyle \frac{s^* - s_f}{F'(s^*)}$
replaced by $\displaystyle \frac{s-s_f}{F'(s)}$.
\end{lem}

\begin{proof}
Assume for the time being that $F''F/(F')^2$ is increasing. It will first be shown that
\begin{equation}\label{eql49a}
-\frac{F' (s)}{s^* - s} \le \frac{F(s) - F(u)}{(s^*- u) (u - s)} \le
\frac{F'(s^*)}{s^* - s} \quad \forall s< u <s^*  \, .
\end{equation}
To establish the left inequality in (\ref{eql49a}), notice that $g:[s,s^*] \to \R$ given by
$$
g(u) = F(s) - F(u) + \frac{F'(s)}{s^*-s} (s^* - u) (u- s) \quad
\forall s\le u \le s^* 
$$
satisfies $g(s) = g(s^*) = 0$, $g'(s) = 0$, and
$$
g'(s^*) = - F'(s^*) - F'(s) = \sqrt{F(s)} \left( 
\sqrt{\frac{F'(s)^2}{F(s)}} - \sqrt{\frac{F'(s^*)^2}{F(s^*)}} \right) \le 0 \, ,
$$
by Proposition \ref{lem47}(ii); also, by (iii) the function $g'$ is
concave. Consequently, if $g'(u_0) = 0$ for some $s<u_0<s^*$ then
$g'(u) \ge 0$ for all $s\le u \le u_0$, and $g(u_0) \ge 0$. Thus
$g(u)\ge 0$ for all $s\le u \le s^*$, which proves the left inequality
in (\ref{eql49a}). A completely analogous argument establishes the
right inequality.

Next it will be shown that
\begin{equation}\label{eql49b}
\frac{(s^* + u)(u+s)}{2u + F(u) - F(s)} \ge s^* + s \quad \forall s\le
u \le s^* \, .
\end{equation}
To see (\ref{eql49b}), similarly to before observe that $h:[s,s^*] \to \R$ given by
$$
h(u) = 2u + F(u) - F(s) - \frac{(s^*+ u)(u+s)}{s + s^*} \quad \forall
s\le u \le s^*
$$
satisfies $h(s) = h(s^*) = 0$, and $h'$ is convex. Since $f$ is
increasing, $h'(u_0)=0$ for a unique $s\le u_0 \le s^*$, and $h(u_0)$
is  a {\em minimal\/} value $\le 0$. In other words, $h(u) \le 0$ for all $s\le u \le s^*$,
i.e., (\ref{eql49b}) holds.

Lastly, deduce from differentiating $F(s^*) = F(s)$ twice that ${\rm
  d}s^* /{\rm d}s |_{s=s_f} = -1$, and hence $s^*\le 2s_f - s$ by
Proposition \ref{lem47}(i),
so
\begin{equation}\label{eq4nnast}
(s^*)^2 - s^2 \ge (s^*)^2 - (2s_f - s^*)^2 = 4s_f (s^* - s_f) \, .
\end{equation}
With these preparations, for every $0<s<s_f$ deduce from (\ref{eq43}),
together with (\ref{eql49a}), (\ref{eql49b}), and (\ref{eq4nnast}) that
\begin{align*}
\omega_f (s) & = \frac1{\pi} \int_s^{s^*} \!\!\!\! \frac{u
  f(u)}{\sqrt{(s^* - u)(u-s)} \sqrt{(s^* + u) (u+s)}} \sqrt{\frac{(s^*
    - u)(u-s)}{F(s) - F(u)}} \sqrt{\frac{(s^* + u)(u+s)}{2u + F(u) - F(s)}}
\, {\rm d}u \\
& \ge \frac1{\pi} \sqrt{\frac{s^* - s}{F'(s^*)}} \sqrt{s^* + s}
\int_s^{s^*} \frac{uf(u) \, {\rm d}u}{\sqrt{\bigl( (s^*)^2 - u^2 \bigr)
  (u^2 - s^2)}} \\
& > \frac1{\pi} \sqrt{\frac{s^* - s_f}{F'(s^*)}} \int_0^{\pi}
\!\! f \! \left( \sqrt{ \frac{(s^*)^2 + s^2 }{2} -\frac{(s^*)^2 - s^2}{2} \cos u}
\right) {\rm d}u \, ,
\end{align*}
which is precisely (\ref{eq48a}). Finally, if $F''F/(F')^2$ is decreasing then it is readily checked
that (\ref{eql49a}) holds with both inequalities reversed, whereas
(\ref{eql49b}) remains valid unchanged, and $s^* \ge 2s_f - s$ since
$s\mapsto s^*$ is convex, so (\ref{eq4nnast}) now reads
$$
(s^*)^2 - s^2 \ge (2s_f - s)^2 - s^2 = 4s_f (s_f - s) \, .
$$
This shows that (\ref{eq48a}) remains valid, provided that $\displaystyle
\frac{s^* - s_f}{F'(s^*)}$ is replaced by $\displaystyle \frac{s-s_f}{F'(s)}$.
\end{proof}

\begin{rem}\label{rem49a}
(i) The right-hand side in (\ref{eq48a}) tends to $1/\sqrt{s_f F_2}$
as $s\to s_f$. Thus the lower bound in Lemma
\ref{lem48} is sharp  at the right end of $]0,s_f[$. Moreover, equality
holds in (\ref{eq48a}) for every $s$ in case $f$ is constant ---
although, strictly speaking, the lemma does not apply in this case
because $f\not \in \cF$.

(ii) A variant of Lemma \ref{lem48} holds for {\em decreasing\/} $f\in
\cF^*$ as well: While (\ref{eql49a}) remains valid in this case also,
the right-hand side in (\ref{eql49b}) has to be replaced by the
trivial lower bound $\frac12 (s^* + s)$. As a consequence, the
right-hand side in (\ref{eq48a}) has to be divided by $\sqrt{2}$,
resulting in a lower bound for $\omega_f(s)$ that typically is not
sharp anywhere.
\end{rem}

This section concludes with a discussion of the monotonicity of
$ \omega_f $. Clearly, if one assumes $\omega_f$ to be,
say, decreasing on $]0, s_f[$, then Theorem \ref{thm45}  together with
Proposition \ref{prop46} immediately yields the cardinality of
$\OO_f$, and thus the number of non-circular Jordan solutions of
(\ref{eq33}).

\begin{theorem}\label{thm46d}
Let $f\in \cF^*$, and assume that $\omega_f$ is decreasing on $]0,
s_f[$. Then {\rm (\ref{eq33})} has precisely $\# \Bigl( \bigl( \N
\setminus \{1\} \bigr) \, \cap
\:\bigr]1/\omega_f (0+ ), \sqrt{s_f F_2 }\bigl[ \, \Bigr)$ different non-circular
Jordan solutions, modulo rotations.
\end{theorem}

\begin{proof}
Since $\omega_{f}$ is continuous and
decreasing on $]0, s_f[$,
$$
\bigl\{ \omega_f (s) : 0 < s < s_f\bigr\} = \left]
  \frac1{\sqrt{s_f F_2 }} , \omega_f (0+)\right[ \, ,
$$
and hence $\# \OO_f $ equals the number of integers $n\ge 2$ with
$1/\omega_f (0+) < n < \sqrt{s_f F_2}$.
\end{proof}

To establish a condition that ensures $\omega_f$ is decreasing on $]0,
s_f[$, notice that for every $f\in \cF^*$ the ratio
$F/F'$ defines a smooth function on $\R^+$, with its value
and derivative for $s= s_f$ equal to $0$ and $\frac12$, respectively. 

\begin{lem}\label{lem46a}
Let $f\in \cF^*$, and assume there exists $a\in \R$
such that
\begin{equation}\label{eq43a}
2\sqrt{s} \frac{{\rm d}}{{\rm d}s} \left( \sqrt{s} f(s)
  \frac{F(s)}{F'(s)}\right) - sf(s) + \frac{f(s) F(s)}{4} \ge a F'(s)
\sqrt{s} \quad \forall 0< s <  0_f^*  \, .
\end{equation}
Then $\omega_f'(s) < 0$ for all $0<s<s_f$.
\end{lem}

\begin{rem}\label{rem410}
The left- and right-hand sides in (\ref{eq43a}) both vanish for
$s=s_f$, with derivatives equal to $F_2 - 1/(2 s_f) - F_3/(3
F_2)$ and $a F_2 \sqrt{s_f}$, respectively. Thus, if $a$
as in Lemma \ref{lem46a} exists at all then necessarily $a = (s_f
(F_2^2 - \frac13 F_3)  - \frac12 F_2 )/( s_f^{3/2} F_2^2)$.
\end{rem}

The proof of Lemma \ref{lem46a} makes use of a simple tailor-made calculus fact
the verification of which once more is left to the interested reader;
see also \cite[Thm.2.1]{CW}.

\begin{prop}\label{prop46b}
Let $f\in \cF^*$, and assume $g:\R^+ \to \R$ is smooth. Then, for
every $0<s<s_f$,
\begin{equation*}
\frac{{\rm d}}{{\rm d}s} \int_s^{s^*} \!\!\! \frac{ g(u) \,  {\rm
    d}u}{\sqrt{F(s) - F(u) }} = \frac{F'(s)}{
  F(s)} \int_s^{s^*} \!\!\! \frac1{\sqrt{F(s) - F(u)}}
\left( 
\frac{{\rm d}}{{\rm d} u} \left( \frac{g(u) F(u)}{F'(u)}\right)
-\frac{g(u)}{2}
\right) {\rm d}u \, .
\end{equation*}
\end{prop}

\begin{proof}[Proof of Lemma \ref{lem46a}]
For every $0<s<s_f$,
\begin{equation}\label{eq46a}
0 \le \frac{F(s) - F(u)}{2u} \le 1 - \frac{s}{s^*} < 1  \quad
\forall s< u <  s^*  \, ,
\end{equation}
and consequently, with (\ref{eq43}) and the binomial formula,
\begin{align*}
\omega_f (s) & = \frac1{\pi \sqrt{2}} \int_s^{s^*}
\frac{\sqrt{u} f(u)}{\sqrt{F(s) - F(u)}} \left( 1 - \frac{F(s) -
    F(u)}{2u}\right)^{-1/2} {\rm d}u \\
& = \frac1{\pi \sqrt{2}} \int_s^{s^*} \sum\nolimits_{n=0}^{\infty}
\left( \!\!  \begin{array}{c} -1/2 \\ n \end{array}\!\!\right) (-1)^n
2^{-n} u^{1/2 - n} f(u) \bigl(F(s) - F(u)\bigr)^{n-1/2}\,  {\rm d}u \\
& = \frac1{\pi \sqrt{2}} \sum\nolimits_{n=0}^{\infty} \left(
  \!\!\!  \begin{array}{c} 2n \\ n \end{array}\!\!\! \right) 2^{-3n} g_n
(s)\\
& = \frac1{\pi \sqrt{2}} \left(  g_0 (s) + \frac{g_1(s)}{4} + \sum\nolimits_{n=2}^{\infty} \left(
  \!\!\! \begin{array}{c} 2n \\ n \end{array} \!\!\!\right) 2^{-3n} g_n
(s)\right) \, ,
\end{align*}
where the second-to-last equality is due to uniform convergence, and
for every integer $n\ge 0$,
$$
g_n (s) = \int_s^{s^*} u^{1/2 - n} f(u) \bigl(F(s) - F(u)\bigr)^{n-1/2} \, {\rm d} u
\quad \forall 0< s < s_f  \, .
$$
The derivative of $g_0$ can be computed using Proposition
\ref{prop46b}, for every $0<s<s_f$, 
$$
g_0'(s) =  \frac{F'(s)}{2F(s)} \int_s^{s^*} \!\!\! \frac1{\sqrt{u \bigl( F(s) -
    F(u)\bigr)}} \left( 
2 \sqrt{u} \frac{{\rm d}}{{\rm d}u} \left( \sqrt{u} f(u)
  \frac{F(u)}{F'(u)}\right) - uf(u)
\right) {\rm d}u \, ,
$$
whereas for $n\ge 1$, the derivative of $g_n$ simply is obtained by
formal differentiation,
$$
g_n'(s) = \bigl( n-{\textstyle \frac12} \bigr) F'(s) \int_s^{s^*} u^{1/2 - n} f(u) \bigl(F(s) -
F(u)\bigr)^{n-3/2} \, {\rm d}u \quad \forall 0<s < s_f  \, .
$$
Note that $g_n'<0$ for all $n\ge 1$ since $F'<0$ on $]0,s_f[$. By means of (\ref{eq46a}), it is
easily seen that the termwise differentiated series $\sum_{n=2}^{\infty} \left(
  \!\!\!  \begin{array}{c} 2n \\ n \end{array}\!\!\! \right) 2^{-3n} g_n'
(s)$ converges locally uniformly on $]0, s_f[$, and hence
\begin{align*}
\omega'_f(s)  & = \frac1{\pi \sqrt{2}} \left(  g_0'(s) + \frac{g_1(s)'}{4} + \sum\nolimits_{n=2}^{\infty} \left(
  \!\! \! \begin{array}{c} 2n \\ n \end{array} \!\!\!\right) 2^{-3n} g_n'
(s)\right) \\
& < \frac1{\pi \sqrt{2}} \left(  g_0'(s) + \frac{g_1'(s)}{4}\right)
\\
& = \frac{F'(s)}{2\pi \sqrt{2} F(s)} \int_s^{s^*} \frac1{\sqrt{u \bigl(F(s)
    - F(u)\bigr)}} \left( 
2 \sqrt{u} \frac{{\rm d}}{{\rm d}u} \left( \sqrt{u} f(u)
  \frac{F(u)}{F'(u)}\right) - uf(u) + \frac{f(u) F(s)}{4}
\right) {\rm d}u \\
& \le \frac{a F'(s)}{2\pi \sqrt{2} F(s)}\int_s^{s^*} \frac{F' (u)\, {\rm
    d} u}{
  \sqrt{F(s) - F(u)}}  = 0 \, ,
\end{align*}
where the last inequality is a consequence of (\ref{eq43a}) and (\ref{eq46a}).
\end{proof}

\begin{rem}\label{rem415}
With a view towards Proposition \ref{prop46b} and the proof of Lemma
\ref{lem46a}, it is tempting to write down an exact formula (rather
than an upper bound) for $\omega_f'(s)$, namely
\begin{align}\label{eq4oo1}
\omega_f'(s) = \frac{F'(s)}{2\pi \sqrt{2} F(s)} \int_s^{s^*} 
\frac1{\sqrt{u \bigl( F(s) -
    F(u)\bigr)}} & \left( 
2 \sqrt{u} \frac{{\rm d}}{{\rm d}u} \left( \sqrt{u} f(u)
  \frac{F(u)}{F'(u)}\right) - uf(u) \right. \\
& \quad \left. + \: \frac{f(u) F(s)}{4} \, \psi \! \left(
  \frac{F(s) -F(u)}{2u}\right)
\right)
\, {\rm
  d}u \, , \nonumber
\end{align}
with the (real-analytic) function
$$
\psi(t) = \frac{2t}{\sqrt{1-t}^3} + \frac2{1+\sqrt{1-t}} \quad \forall
t < 1 \, .
$$
Note that $\psi$ is convex on $[0,1[$, with $\psi(0) = 1$ and $\psi'(0)=\frac94$.
Due to its ``non-local'' nature, (\ref{eq4oo1})
appears to be rather unwieldy. In particular, the integrand typically
changes sign in $]s,s^*[$. As a consequence, any general statement about the sign
of $\omega_f'$ is bound to be a delicate affair, a fact for which the
next section is going to provide ample evidence.
\end{rem}

\section{An example: The monomial family}\label{sec5}

This final section applies the results of the preceding sections
to the monomials $f_{\delta}(s) = s^{\delta}$, with ${\delta}\in \R$. Naturally, the analysis
is quite specific to that particular family of
functions, though the techniques applied here likely are
useful also when dealing with other classes of
functions. Moreover, the section illustrates how applying the results of
the present article, notably Theorem \ref{thm46d}, though quite trivial in
theory, may nonetheless pose a considerable challenge in
practice. This is not an uncommon situation: The reader likely is familiar
with similar, seemingly simple problems in non-linear analysis that
also require for their resolution lengthy, potentially delicate and unenlightening
computations; see, for instance, \cite{AL, benguria, MY,y}. 

Recall from the previous section that $f_{\delta}\in \cF$ for every
$\delta \in \R \setminus \{-1,0\}$, and first consider the case
$\delta< -1$, where $\epsilon_{f_{\delta}} = -1$ and $\I_{f_{\delta}} = -
\R^+$. Consequently, Proposition \ref{prop42} shows that  the only
Jordan solution of $\kappa = r^{\delta}$ is the (counter-clockwise
oriented) unit circle. Next, the case $\delta = -1$ has 
been considered already in Example \ref{exa210}(ii): Every
(counter-clockwise oriented) circle centered at $0$ is a Jordan
solution of $\kappa = 1/r$, and there are no other maximal solutions
that are bounded, let alone closed or Jordan. It remains to
consider the case $\delta > - 1$, where $f_{\delta} \in \cF^*$ unless
$\delta = 0$. In this case, $\epsilon_{f_{\delta}} = 1$ and
$\I_{f_{\delta}} = \R^+$, so Theorem \ref{thm45} applies. Moreover,
$s_{f_{\delta}} = 1$, $0^*_{f_{\delta}} = (\delta + 2)^{1/(\delta +
  1)}$, and
$$
F  (s) = \frac{s^{\delta + 2} - (\delta + 2) s + \delta +
  1}{\delta + 2} \quad \forall s \in \R^+ \, ,
$$
as well as $F_k = \prod_{\ell =2}^k (\delta + 3 - \ell )$ for
all $k\ge 2$. By Proposition \ref{prop46}, $\omega_{f_{\delta}}$ is a
smooth positive function on $]0,1[$ with
\begin{equation}\label{eq51}
\omega_{f_{\delta}} (1-) = \frac1{\sqrt{\delta
  + 1}} \, , \quad
\omega'_{f_{\delta}}
(1-) = 0  \, , \quad
\omega_{f_{\delta}} ''(1-)= \frac{\delta^2}{12 \sqrt{\delta
  + 1}} \, .
\end{equation}
Thus if $-1<\delta < 0$ then $\omega_{f_{\delta}} (s)>1$ for all
$0<s<1$, by Lemma \ref{lem27}, hence $\OO_{f_{\delta}} =
\varnothing$, and again the only Jordan solution of $\kappa =
r^{\delta}$ is the (counter-clockwise oriented) unit circle. By
contrast, if $\delta = 0$ then $\omega_{f_0} (s)=1$ for
all $0<s<1$. Though (\ref{eq3ast}) fails in this case, as $f_0
\not \in \cF$, it is clear
that every (counter-clockwise oriented) circle with radius $1$ is a
Jordan solution, and there are no other maximal solutions
whatsoever. Correspondingly, $s^* = 2-s$ and
$$
\omega_{f_0} (s) = \frac2{\pi} \int_s^{2-s} \frac{u \, {\rm
    d}u}{\sqrt{u^2 - s^2} \sqrt{(2-s)^2 - u^2}} = \frac1{\pi}
\int_{s^2}^{(2-s)^2} \frac{{\rm d}u}{\sqrt{u-s^2} \sqrt{(2-s)^2 - u}}
= 1  \quad \forall 0<s<1 \, .
$$
The only case yet to be considered, therefore, is $\delta>0$. In this case, 
\begin{align}\label{eq52}
\omega_{f_{\delta}} (0+) & = \frac1{\pi} \int_0^{0_{f_{\delta}}^*}
\frac{u^{\delta + 1} {\rm d}u}{\sqrt{F(0) - F(u) } \sqrt{2u + F(u) -
    F(0)}}  \nonumber \\
& =
\frac1{\pi} \int_0^{(\delta + 2)^{1/(\delta + 1)}} \!\!\!\!
\!\!\!\!\!\! \frac{u^{\delta +
    1}{\rm d}u}{ \sqrt{u^2 - u^{2\delta + 4}/(\delta + 1)^2}} =
\frac12 + \frac1{2 (\delta + 1)}  \, ,
\end{align}
and consequently
$$
\omega_{f_{\delta}} (0+ ) -
\omega_{f_{\delta}} (1-) = \frac{\delta^2}{2(\delta + 1) (\delta + 2
  + 2 \sqrt{\delta + 1})} > 0 \, .
$$
Since $\omega_{f_{\delta}} $ attains a non-degenerate local
minimum as $s\to 1$ by (\ref{eq51}), the following certainly is a
plausible speculation.

\begin{conjecture}\label{con51}
For every $\delta >0 $ the function $\omega_{f_{\delta}}$ is decreasing on $]0,1[$.
\end{conjecture}

At the time of this writing, the author has been able to establish
the correctness of this conjecture only for $\delta \ge \frac32$; see Lemma
\ref{lem52} below. For smaller $\delta$, a somewhat weaker substitute is presented. Concretely, observe that
$\frac12 < \omega_{f_{\delta}} (0 + )< 1$ for all $\delta >0$, and if $0< \delta < 3$ then also
$\frac12 < \omega_{f_{\delta}} (1-) < 1$. If $\delta $
is not too large then these bounds are valid for all intermediate
values as well. 

\begin{lem}\label{lem53}
If $0< \delta \le \frac32$ then $\frac12 < \omega_{f_{\delta}} (s) <
1$ for all $0<s<1$.
\end{lem}

The proof of Lemma \ref{lem53} presented below makes use of several
inequalities, two of which may be of independent interest: On the one
hand, elementary calculus shows that
$$
a\le \left( \frac{(a+1)^{b+1} -1}{a}\right)^{1/b} - (b+1)^{1/b} \le a
\frac{(b+1)^{1/b}}{2} \quad \forall a \in \R^+, 0 < b \le 1 \, ,
$$
whereas for $b\ge 1$ both inequalities are reversed. In other words,
since $\max \left\{ \frac12 (b+1)^{1/b}, 1 \right\}$ equals $\frac12 (b+1)^{1/b}$ if $0<b\le
1$, and equals $1$ if $b\ge 1$,
\begin{equation}\label{eq5mm1}
a \min \left\{  \frac{(b+1)^{1/b}}{2} , 1 \right\} \le \left(
  \frac{(a+1)^{b+1} -1}{a}\right)^{1/b} \!\! - (b+1)^{1/b} \le a \max
\left\{ \frac{(b+1)^{1/b}}{2}, 1 \right\} 
\quad \forall a,b \in \R^+ \, .
\end{equation}
On the other hand, as a special case of an optimal Gautschi
inequality established in \cite{EGP},
\begin{equation}\label{eq5mm2}
a + \frac14 < \frac{\Gamma (a+1)^2}{\Gamma (a+1/2)^2} < a + \frac1{\pi}
\quad \forall a \in \R^+ \, ,
\end{equation}
and both bounds are sharp, as the left (respectively, right) inequality
becomes an equality as $a\to \infty$ (respectively, $a\to 0$).

\begin{proof}[Proof of Lemma \ref{lem53}]
Since $0<\omega_{f_{\delta}}(s)<1$ for all $0<s<1$ whenever
$\delta>0$, by Lemma \ref{lem27}(iii), it only needs to be shown that
$\omega_{f_{\delta}}(s)>\frac12 $. Clearly, $f_{\delta}$ is 
increasing, and it is readily checked that $F''F/(F')^2$ is increasing on
$\R^+$ as well. By (\ref{eq48a}),
$$
\omega_{f_{\delta}} (s) \ge \frac1{2^{\delta/2}\pi}
\sqrt{\frac{s^* - 1}{(s^*)^{\delta + 1} -1}} \bigl(
(s^*)^2+ s^2 \bigr)^{\delta / 2} \int_0^{\pi} \left( 1 - \frac{(s^*)^2 -
    s^2}{(s^*)^2 + s^2} \cos u \right)^{\delta/2} \!\! {\rm d}u \quad
\forall 0<s< 1 \, ,
$$
and utilizing the elementary estimate 
$$
\int_0^{\pi} (1 - a \cos u)^b {\rm d}u \ge \frac{2^b \sqrt{\pi} \, \Gamma
(b + 1/2 )}{ \Gamma (b+1)} \quad \forall 0\le a,b \le 1  \, ,
$$
it follows that, for every $0< \delta \le 2$,
$$
\omega_{f_{\delta}} (s)^2  \ge \frac1{\pi} \cdot \frac{s^* -
  1}{(s^*)^{\delta + 1} -1} \bigl(
(s^*)^2+ s^2 \bigr)^{\delta } \frac{\Gamma \bigl( (\delta +
  1)/2 \bigr)^2}{ \Gamma \bigl( (\delta +2)/2 \bigr)^2} \quad \forall
0<s< 1\, .
$$
Recall that $s^*>1$ for every $0<s<1$. Thus, applying (\ref{eq5mm1})
and (\ref{eq5mm2}) yields the lower bound, valid whenever $0<\delta
\le 2$,
\begin{equation}\label{eq5mm3}
\omega_{f_{\delta}}(s)^2 > \frac{\bigl( (s^*)^2 + s^2
  \bigr)^{\delta}}{\bigl( (\delta+1)^{1/\delta} + (s^* - 1) \max \{ 
  (\delta + 1)^{1/\delta}/2 , 1 \}\bigr)^{\delta}} \cdot \frac{2}{\pi
  \delta + 2} \quad \forall 0<s<1 \, .
\end{equation}
Also, recall from Proposition \ref{lem47}(i) that $s\mapsto s^*$ is
concave, with $s^*|_{s=0} =0_{f_{\delta}}^* = (\delta + 2)^{1/(\delta + 1)}$ and
$s^*|_{s=1}=1$, and hence
\begin{equation}\label{eq5mm4}
s^* \ge (\delta + 2)^{1/(\delta + 1)} (1-s) + s \quad \forall 0<s<1 \, .
\end{equation}
Replacing $s$ on the right in (\ref{eq5mm3}) by the lower bound in
terms of $s^*$ provided by (\ref{eq5mm4}), and requiring that the
resulting expression still be $>\frac14$ for all $0<s<1$ is equivalent to requiring that
$$
p(\delta, s^* -1) > 0 \quad \forall 1<s^* < (\delta+2)^{1/(\delta +
  1)} \, ,
$$
with the continuous function $p:\R^+ \times \R \to \R$ given by
$$
p(\delta , t) = (t+1)^2 + \left( \frac{t}{(\delta + 2)^{1/(\delta + 1)}
    -1} - 1\right)^2 - \left( \frac{\pi \delta + 2}{8}
\right)^{1/\delta} \left( (\delta + 1)^{1/\delta} + t \max \left\{ 
    \frac{(\delta + 1)^{1/\delta}}{2}, 1 \right\} \right) \, .
$$
Thus the assertion of the lemma immediately follows, as soon as
it is shown that in fact
\begin{equation}\label{eq5mm5}
p(\delta , t ) > 0 \quad \forall 0< \delta \le \frac32 , t \in \R \, .
\end{equation}
In other words, to prove the lemma it suffices to establish
(\ref{eq5mm5}), and this will now be done. To this end, notice that $p(\delta, \, \cdot \, )$ is a
quadratic polynomial,
$$
p(\delta, t) = p_2(\delta)t^2 - 2p_1(\delta) t + p_0(\delta) \quad \forall
\delta>0, t \in \R \, ,
$$
with continuous, positive coefficients $p_0, p_1, p_2: \R^+ \to \R^+$
given by
\begin{align*}
p_0(\delta) & = 2 - \left( \frac{\pi \delta + 2}{8}
\right)^{1/\delta} (\delta + 1)^{1/\delta} \, , \\
p_1(\delta) & = \frac1{(\delta + 2)^{1/(\delta + 1)}-1} -1 + \frac12 \left( \frac{\pi \delta + 2}{8}
\right)^{1/\delta} \max \left\{    \frac{(\delta +
    1)^{1/\delta}}{2} , 1 \right\} \, , \\
p_2(\delta) & = \frac1{\bigl( (\delta + 2)^{1/(\delta + 1)}-1\bigr)^2}
+ 1 \, , 
\end{align*}
and hence (\ref{eq5mm5}) holds, provided that
\begin{equation}\label{eq5mm6}
p_0(\delta) p_2(\delta) > p_1(\delta)^2 \quad \forall 0 < \delta \le
\frac32 \, .
\end{equation}
Now, it is readily checked that $p_0$ and $p_1$ are decreasing and
increasing on $\left]0, \frac32 \right]$, respectively, and hence to establish
(\ref{eq5mm6}), it suffices to verify that
\begin{equation}\label{eq5mm7}
\sqrt{p_0(1)} \sqrt{p_2(\delta)} > p_1(1) \quad \forall 0 < \delta \le
1 \quad \mbox{\rm and} \quad \sqrt{p_0\left(\frac32 \right)} \sqrt{p_2 (\delta)} >
p_1 \left( \frac32 \right) \quad \forall 1 \le \delta \le \frac32  \, .
\end{equation}
Notice that
$$
\sqrt{p_2(\delta)} \ge \frac1{\sqrt{2}}\cdot \frac{(\delta +
  2)^{1/(\delta + 1)}}{(\delta +
  2)^{1/(\delta + 1)} -  1}\quad \forall \delta \in \R^+ \, ,
$$
and since the lower bound on the right is increasing in $\delta$, clearly
(\ref{eq5mm7}) holds, provided that
\begin{equation}\label{eq5mm8}
\sqrt{p_0(1)} \frac{2}{\sqrt{2}} > p_1 (1) \quad \mbox{\rm and}
\quad \sqrt{p_0 \left(\frac32 \right) } \frac1{\sqrt{2}}\cdot
\frac{\sqrt{3}}{\sqrt{3}-1} > p_1 \left(\frac32 \right) \, .
\end{equation}
Utilizing the rough (rational) estimates
\begin{align*}
p_0(1) & =  \frac{6 - \pi}{4} > \frac23  \, , \quad 
p_0\left(\frac32 \right)  = 2 - \frac1{16} (40 + 30 \pi)^{2/3} >\left(
  \frac35 \right)^2 \,
, \\
p_1(1) & = \frac{8 \sqrt{3} - 6 + \pi}{16} < \frac34  \, , \quad
p_1\left(\frac32 \right)  = \frac{2^{7/5} -7^{2/5}}{7^{2/5} - 2^{2/5}}  + \frac1{32} (16 +
12 \pi)^{2/3} < 1 \, ,
\end{align*}
it is readily seen that (\ref{eq5mm8}) indeed is correct. This proves
(\ref{eq5mm6}), which in turn implies (\ref{eq5mm5}). As detailed
earlier, the latter proves the lemma.
\end{proof}

As a consequence of Lemma \ref{lem53}, $\OO_{f_{\delta}} =
\varnothing$ also when $0< \delta \le \frac32$, and again the only
Jordan solution of $\kappa = r^{\delta}$ is the (counter-clockwise
oriented) unit circle.

Lastly consider the case $\delta \ge \frac32$. In this case, the
conclusion of Conjecture \ref{con51} definitely holds.

\begin{lem}\label{lem52}
If $\delta \ge \frac32$ then the  function $\omega_{f_{\delta}} $ is decreasing on $]0,1[$.
\end{lem}

The proof of Lemma \ref{lem52} presented below makes use of a
simple calculus fact: Given $n\in \N$, non-zero real numbers $a_1,
\ldots , a_n$, and real numbers $b_1> \ldots > b_n$, consider
the real-analytic function $g:\R \to \R$ given by
\begin{equation}\label{eq5nn0}
g(t) = \sum\nolimits_{\ell =1}^n a_{\ell} e^{\textstyle b_{\ell}  t} \, ,
\end{equation}
and let $\sigma(g)$ be the number of sign changes in the finite
sequence $(a_1, \ldots , a_n)$, more formally, $\sigma(g) = \# \{ 1
\le \ell \le n : a_{\ell-1} a_{\ell} < 0\}$ where $a_0:= a_1$. Plainly, $0\le
\sigma(g) \le n-1$. Since $\lim_{t\to -\infty} g(t) e^{\textstyle -b_n t} =
a_n \ne 0$ and $\lim_{t\to \infty} g(t) e^{\textstyle -b_1 t} = a_1 \ne 0$,
the equation $g(t)=0$ has only finitely many real roots, each of which
has finite multiplicity. The following variant of Descartes' rule
\cite{curtiss} makes this more precise.

\begin{prop}\label{prop5nn1}
Let $n\in \N$, $a_1, \ldots , a_n \in \R \setminus
\{0\}$, and $b_1, \ldots, b_n \in \R$ with $b_1 >
\ldots > b_n$. Then the total number of real roots (counted with
multiplicities) of $g(t)=0$, with $g$ as in {\rm (\ref{eq5nn0})}, equals
$\sigma(g) - 2k $ for some $k \in \N_0$.
\end{prop}

\begin{proof}[Proof of Lemma \ref{lem52}]
Though computationally intense in its details, the following argument
has a simple basic strategy: Intending to utilize Lemma
\ref{lem46a} with $f=f_{\delta}$, notice first that (\ref{eq43a}),
with $a= \frac16  (4\delta + 3)/(\delta + 1)$ as per Remark \ref{rem410},
can be written equivalently but more concisely as
$$
p\left( 2(\delta + 1), \frac{\log s}{2} \right) \ge 0 \quad \forall 0<s< (\delta +
2)^{1/(\delta + 1)} \, ,
$$
with the real-analytic function $p:\R^2 \to \R$ given by
\begin{align}\label{eq5nn1}
p(\varepsilon , t)  =  & \enspace 6 \varepsilon  e^{(4 \varepsilon  -
  1)t}- 4 ( \varepsilon+ 2)  ( 2 \varepsilon - 1 ) e^{3 \varepsilon t}
+9 \varepsilon (\varepsilon-2) e^{(3\varepsilon - 1)t } + 3
\varepsilon^2 t^{3(\varepsilon- 1)t} \nonumber \\ & + 12
(\varepsilon+2) (2\varepsilon-1) e^{2 \varepsilon t}  
  - 6\varepsilon (5\varepsilon-3) e^{(2 \varepsilon - 1)t} - 18\varepsilon^2 e^{(2\varepsilon-3)t}   \\
&  - 12 (\varepsilon+2) (2\varepsilon -1 )
e^{\varepsilon t} + 3\varepsilon(\varepsilon+2) (4\varepsilon-1) e^{(\varepsilon-1 )t}
-3  \varepsilon^2 (4\varepsilon - 5) e^{(\varepsilon-3) t }+ 4 (\varepsilon+2) (2\varepsilon-1)
\, . \nonumber 
\end{align}
Since $\delta \ge \frac32$ precisely if $\varepsilon = 2(\delta + 1) \ge
5$, the assertion of the lemma immediately follows from Lemma
\ref{lem46a}, as soon as it is shown that in fact
\begin{equation}\label{eq5nn2}
p(\varepsilon, t) \ge 0 \quad \forall \varepsilon \ge 5 , t \in \R \, .
\end{equation}
Thus, to prove the lemma it suffices to establish
(\ref{eq5nn2}), and this will now be done in several steps. Usage of the
same symbols as in the proof of Lemma \ref{lem53} will hopefully not
confuse the reader but rather highlight the parallels between both proofs.

Henceforth assume $\varepsilon \ge 5$, and notice at
the outset that $\sigma \bigl(  p(\varepsilon , \, \cdot \,
)\bigr)=6$, as well as $\partial^k p /\partial t^k (\varepsilon, 0)
=0$ for $k=0,1,2,3$, whereas $\partial^4 p /\partial t^4 (\varepsilon, 0)
= 24 \varepsilon^3 (\varepsilon + 2) (2\varepsilon^2 - 11 \varepsilon + 8)>0$. For every $\varepsilon >
0$, therefore, $p(\varepsilon , \, \cdot \, )$ has $t=0$ as a $4$-fold
root, and so by Proposition \ref{prop5nn1} has either two or zero
additional real roots. In the latter case, clearly (\ref{eq5nn2}) is
correct.

First it will be shown that $p(\varepsilon, t)\ge 0$ for all $t\in
\R$, provided that $\varepsilon$ is large enough. To this end,
consider the real-analytic function $\widehat{p}:\R^2 \to \R$ given
by
$$
\widehat{p} (\varepsilon , t) = 2 \varepsilon e^{(4\varepsilon - 1)t} - (2
\varepsilon - 1) (3 \varepsilon - 1) e^{3 \varepsilon t} + 2 \varepsilon
(3 \varepsilon-1) e^{(3 \varepsilon - 1)t}- 2(3 \varepsilon-1) e^{2
  \varepsilon t} + (\varepsilon-1) e^{\varepsilon t} \, .
$$
Notice that $\sigma \bigl(\widehat{p} (\varepsilon , \, \cdot \, )
\bigr)=4$, with $t=0$ being a $4$-fold root, and
consequently $\widehat{p} (\varepsilon, t) \ge 0$ for all $t\in \R$,
by Proposition \ref{prop5nn1}. Thus
\begin{equation}\label{eq5nn3}
p(\varepsilon , t) \ge p(\varepsilon,  t) - 4\frac{\varepsilon + 2}{3
  \varepsilon - 1} \widehat{p} (\varepsilon , t) = 2\varepsilon
\frac{5\varepsilon - 11}{3\varepsilon - 1}  \varepsilon
e^{(4\varepsilon - 1)t} + 0 \cdot e^{3 \varepsilon t}  + \varepsilon
(\varepsilon - 34) e^{(3 \varepsilon - 1)t} + 3\varepsilon^2 e^{3(
  \varepsilon - 1)t} + \ldots \, ,
\end{equation}
where $\ldots$ indicates the remaining terms of $p(\varepsilon , t)$
from (\ref{eq5nn1}), of which only the coefficients of
$e^{2\varepsilon t}$ and $e^{\varepsilon t}$ differ from those in
(\ref{eq5nn1}), both being larger now in absolute value but having
kept their respective signs. Now, if $\varepsilon \ge 34$ then $\sigma
\bigl( p(\varepsilon, \, \cdot \, ) - 4(\varepsilon +2)/(3\varepsilon
- 1) \widehat{p} (\varepsilon , \, \cdot \, ) \bigr) = 4$, and since
$t=0$ is a $4$-fold root of both $p(\varepsilon, \, \cdot \, ) $ and
$\widehat{p} (\varepsilon , \, \cdot \, ) $, the right-hand side of
(\ref{eq5nn3}) is $\ge 0$ for all $t\in \R$, by Proposition
\ref{prop5nn1}. In particular, $p(\varepsilon, t)\ge 0$ for all
$\varepsilon \ge 34$ and $t\in \R$.

To consider the remaining cases in (\ref{eq5nn2}), for the remainder
of this proof assume
$5\le \varepsilon \le 34$, and for computational convenience let
$\varepsilon = \nu + 5$, so $0 \le \nu \le 29$. Notice that
$$
p(\varepsilon, t) = \sum\nolimits_{n=0}^{\infty} \frac{\partial ^n
  p}{\partial t^n} (\varepsilon , 0) \frac{t^n}{n!} = t^4
\sum\nolimits_{n=0}^{\infty} p_n(\nu) \frac{t^n}{(n+4)!} \, ,
$$
where for every $n\in \N_0$,
\begin{align*}
p_n (\nu)  =  & \enspace \frac{\partial ^{n+4}
  p}{\partial t^{n+4}} (\nu + 5, 0)  \\
 =  & \enspace 6(\nu + 5) (4\nu + 19)^{n+4} -
4(\nu + 7) (2\nu + 9) (3\nu + 15)^{n+4} \\
& + 9(\nu + 5) (\nu + 3) (3\nu +
14)^{n+4} + 3(\nu + 5)^2 (3\nu + 12)^{n+4} + 12 (\nu + 7) (2\nu + 9)
(2\nu + 10)^{n+4} \\
& - 6 (\nu + 5) (5\nu + 22) (2\nu + 9)^{n+4} - 18 (\nu
+ 5)^2 (2\nu + 7)^{n+4} - 12 (\nu + 7) (2\nu + 9) (\nu +5)^{n+4} \\
& +
3(\nu + 5) (\nu + 7) (4\nu + 19) (\nu + 4)^{n+ 4} - 3(\nu + 5)^2 (4\nu
+ 15) (\nu + 2)^{n+4}
\end{align*}
is a polynomial of degree $n+6$; for example,
\begin{align*}
p_0 (\nu) & = 48\nu^6 + 1272 \nu^5 + 13464 \nu^4 + 71664 \nu^3 +
195360 \nu^2 + 235800 \nu + 63000 \\
& = 24 (\nu+5)^3 (\nu + 7) (2\nu^2 + 9\nu + 3) \, ,\\
p_1(\nu) & = 336 \nu^7 + 10440 \nu^6 + 135240 \nu^5 + 939840 \nu^4 +
3734784 \nu^3 + 8270760 \nu^2 + 8893800 \nu + 2898000 \\
& = 24 (\nu+5)^3 (\nu + 7) (14\nu^3 + 127 \nu^2 + 321 \nu + 138) \, .
\end{align*}
To establish (\ref{eq5nn2}) for all $5\le \varepsilon \le 34$ and
$t\ge 0$, clearly it is enough to demonstrate that
\begin{equation}\label{eq5nn4}
p_n (\nu) \ge 0 \quad \forall 0 \le \nu \le 29, n \in \N_0 \, .
\end{equation}
In order to do this, notice that for all $0\le \nu \le 29$ the three
inequalities
\begin{align*}
3(\nu + 5) (4\nu + 19)^{n+4} & \ge 2 (\nu + 7) (2\nu + 9) (3\nu +
15)^{n+4} \, , \\
3(\nu + 3) (3\nu + 14)^{n+4} & \ge 2 (5\nu + 22) (2\nu + 9)^{n+4} \, ,
\\
(3\nu + 12)^{n+4} & \ge 6(2\nu + 7)^{n+4} \, ,
\end{align*}
hold simultaneously, and hence $p_n(\nu) \ge 0$, provided that
$n\ge 13$. Thus (\ref{eq5nn4}) is correct for all $n\ge 13$. For the
remaining cases $2\le n\le 12$, a straightforward albeit tedious
calculation (involving only integers, and much aided by symbolic
mathematical software) shows that, just as for $p_0$ and $p_1$, all
coefficients of $p_n$ are positive (integers). Hence $p_n(\nu)\ge 0$
for all $0\le \nu\le 29$ and $n\in \N_0$, i.e., (\ref{eq5nn4}) indeed
is correct. As seen earlier, this yields $p(\varepsilon,t)\ge 0$ for
all $5\le \varepsilon \le 34$ and $t\ge 0$.

Finally, it remains to establish (\ref{eq5nn2}) for $5\le \varepsilon
\le 34$ and $t< 0$. To this end, consider the real-analytic function 
$q:\R^2 \to \R$ given by
\begin{align*}
q(\varepsilon, t) = & \enspace e^{(4\varepsilon - 1) t} p(\varepsilon
, -t) \\
= & \enspace 4 (\varepsilon+2) (2\varepsilon-1) e^{(4 \varepsilon  -
  1)t} - 3  \varepsilon^2 (4\varepsilon - 5) e^{(3 \varepsilon+2) t }
+ 3\varepsilon(\varepsilon+2) (4\varepsilon-1) e^{3 \varepsilon t}
\\
& - 12 (\varepsilon+2) (2\varepsilon -1 ) e^{(3 \varepsilon - 1) t} 
- 18\varepsilon^2 e^{ 2(\varepsilon+1) t}  - 6\varepsilon (5\varepsilon-3) e^{2 \varepsilon t}
+ 12 (\varepsilon+2) (2\varepsilon-1) e^{(2 \varepsilon -1)  t}  \\
& + 3
\varepsilon^2 t^{(\varepsilon + 2)t} +9 \varepsilon (\varepsilon-2)
e^{\varepsilon  t }- 4 ( \varepsilon+ 2)  ( 2 \varepsilon - 1 ) e^{(
  \varepsilon - 1)  t} + 6 \varepsilon  \, ,
\end{align*}
and observe that the validity of $p(\varepsilon, t)\ge 0$ for $5\le
\varepsilon \le 34$ and $t<0$ follows from
\begin{equation}\label{eq5nn5}
q(\varepsilon  , t) \ge 0 \quad \forall 5 \le \varepsilon \le 34 ,
t\ge 0 \, ,
\end{equation}
so it is sufficient to establish (\ref{eq5nn5}). To do this, it is
natural to imitate the earlier argument: Write
$$
q(\varepsilon, t) = \sum\nolimits_{n=0}^{\infty} \frac{\partial ^n
  q}{\partial t^n} (\varepsilon , 0) \frac{t^n}{n!} = t^4
\sum\nolimits_{n=0}^{\infty} q_n(\nu) \frac{t^n}{(n+4)!} \, ,
$$
where for every $n\in \N_0$,
\begin{align*}
q_n (\nu)  =  & \enspace \frac{\partial ^{n+4}
  q}{\partial t^{n+4}} (\nu + 5, 0)  \\
 =  & \enspace 4(\nu + 7) (2\nu + 9) (4\nu + 19)^{n+4} -
3(\nu + 5)^2 (4\nu + 15) (3\nu + 17)^{n+4} \\
& + (\nu + 7) (4\nu + 19) (3\nu +
15)^{n+5} -12 (\nu + 7) (2\nu + 9) (3\nu + 14)^{n+4} \\
& - 18  (\nu + 5)^2 
(2\nu + 12)^{n+4} - 3  (5\nu + 22) (2\nu + 10)^{n+5} + 12  (\nu
+ 7) (2\nu + 9)^{n+5} \\
& +3  (\nu + 5)^2  (\nu +7)^{n+4} +
9(\nu + 3)  (\nu + 5)^{n+ 5} - 4(\nu + 7) (2 \nu + 9)  (\nu + 4)^{n+4}
\end{align*}
again is a polynomial of degree $n+6$; for example,
\begin{align*}
q_0 (\nu) & =  p_0 (\nu) = 24 (\nu+5)^3 (\nu + 7) (2\nu^2 + 9\nu + 3) \, ,\\
q_1(\nu) & = 624 \nu^7 + 19560 \nu^6 + 254880 \nu^5 + 1772520 \nu^4 +
6980496 \nu^3 + 15004440 \nu^2 \\
 & \quad + 14767200 \nu + 3087000 \\
& = 24 (\nu+5)^3 (\nu + 7) (26\nu^3 + 243 \nu^2 + 594 \nu + 147) \, .
\end{align*}
In complete analogy to (\ref{eq5nn4}), it suffices to show that
\begin{equation}\label{eq5nn6}
q_n(\nu) \ge 0 \quad \forall 0\le \nu \le 29, n \in \N_0 \, .
\end{equation}
For all $0\le \nu\le 29$, it is readily checked that
\begin{align*}
(4\nu + 19)^{n+4} \ge \max &  \left\{ 
\frac{3(\nu+5)^2 (4\nu + 15)}{(\nu + 7)(2\nu + 9)} (3\nu + 17)^{n+4},
12(3\nu +14)^{n+4}, \right. \\
& \enspace \left. \frac{18(\nu+5)^2}{(\nu+7)(2\nu + 9)} (2\nu +
12)^{n+4}, \frac{3(5\nu + 22)}{(\nu+7)(2\nu + 9)} (2\nu + 10)^{n+5}
\right\} \, , 
\end{align*}
and hence $q_n(\nu)\ge 0$, provided that $n\ge 44$. Similarly to
before, it can be confirmed by direct calculation that for $q_2,
\ldots , q_{43}$ all coefficients are positive (integers), just as for
$q_0$ and $q_1$. In other words, (\ref{eq5nn6}) is correct, which in
turn yields $p(\varepsilon, t) \ge 0$ for all $5\le \varepsilon \le
34$ and $t\le 0$. At long last, therefore, (\ref{eq5nn2}) has been
established. As detailed earlier, an application of Lemma
\ref{lem46a} now completes the proof.
\end{proof}

\begin{figure}[ht]
\psfrag{tpe}[]{\small $p(\varepsilon, t)>0$}
\psfrag{tt1}[]{\small $t$}
\psfrag{taa}[]{\small $\A$}
\psfrag{te1}[]{\small $\varepsilon$}
\psfrag{tlem52}[l]{\small Lemma \ref{lem52}}
\psfrag{tangle}[]{\small $\pi/6$}
\psfrag{tnm2}[]{\small $-2$}
\psfrag{tn0}[]{\small $0$}
\psfrag{tn2}[]{\small $2$}
\psfrag{tn5}[]{\small $5$}
\psfrag{tvn5}[]{\small $5$}
\psfrag{tvnm5}[]{\small $-5$}
\psfrag{tt0}[]{\tiny $0$}
\psfrag{tt2}[]{\tiny $2$}
\psfrag{ttm2}[]{\tiny $-2$}
\psfrag{tt10}[]{\tiny $10$}
\psfrag{tt20}[]{\tiny $20$}
\psfrag{tt30}[]{\tiny $30$}
\psfrag{tzv0}[]{\tiny $10^0$}
\psfrag{tzv20}[]{\tiny $10^{20}$}
\psfrag{tzv40}[]{\tiny $10^{40}$}
\psfrag{tzv60}[]{\tiny $10^{60}$}
\psfrag{ttt}[]{\small $t$}
\psfrag{tte}[]{\small $\varepsilon$}
\psfrag{tpp}[]{\small $\displaystyle \frac{p(\varepsilon ,
    t)}{e^{2\varepsilon t} t^4}$}
\vspace*{-4mm}
\begin{center}
\includegraphics{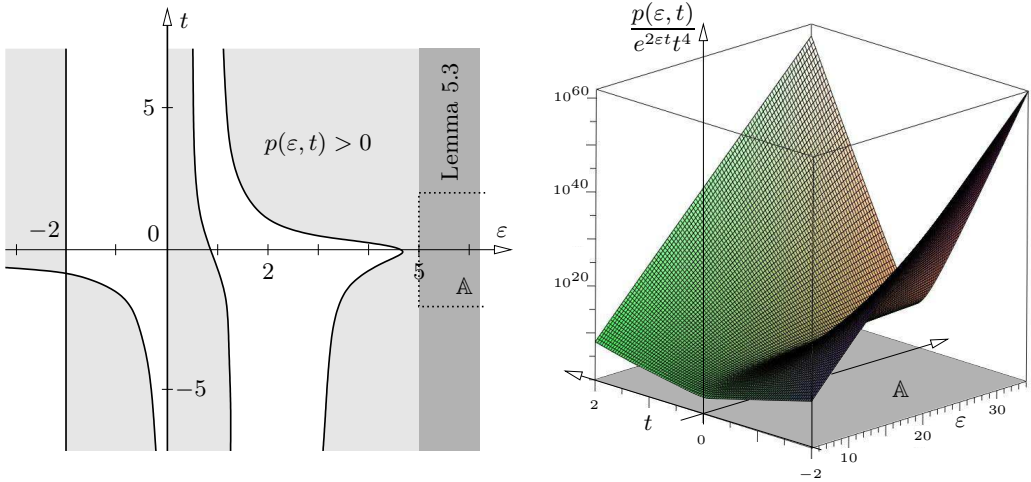}
\end{center}
\vspace*{-6mm}
\caption{Solid black curves indicate the zero locus of
  $p=p(\varepsilon , t)$ given by
  (\ref{eq5nn1}), with the set $\{ p >
  0 \}$ shown in grey (left). Plotting the real-analytic function
  $p(\varepsilon,t)e^{-2\varepsilon t}/t^4$ suggests that $p(\varepsilon, t) \ge a
  t^4e^{2\varepsilon t}$ for all $\varepsilon\ge 5$, $t\in \R^+$,
  where $a$ may be as large as $10^{3.416}= 2606$ (right).}\label{figa1}
\end{figure}

\begin{rem}\label{rem5nn1}
(i) With the same symbols as in the proof of Lemma \ref{lem52}, for every $n\in \N_0$ let
$\nu_n$ be the largest real root of $p_n q_n$. In essence, the above  proof
hinges on the fact that $\nu_n < 0$ for every $n$. A careful analysis
reveals that $\sup_{n\in \N_0} \nu_n = \nu_1 = -0.2781$, and hence the
same argument could be used to establish the
monotonicity of $\omega_{f_{\delta}}$, i.e., the conclusion of
Conjecture \ref{con51}, whenever $\delta \ge \frac12 (5+\nu_1) - 1 =
1.360$. It can be checked numerically, however, that $\inf_{t\in \R}
p(\varepsilon , t) <0$ whenever $2<\varepsilon \le 4.660$. Lemma
\ref{lem46a} therefore is incapable of establishing Conjecture
\ref{con51} for $0<\delta \le 1.330$.

(ii) Numerical evidence strongly suggests that $p(\varepsilon , t) \ge
a t^4 e^{2\varepsilon t}$ for an appropriate $a >0$ and all $\varepsilon
\ge 5$, $t\in \R$. Obviously, such a lower
bound on $p$, if indeed correct, implies (\ref{eq5nn2}). In the light of this, an alternative proof of Lemma \ref{lem52}
might be provided using rigorous (or validated) numerics; see, e.g.,
\cite{AH, tucker} for context. For
$\varepsilon \ge 5$, it is easy to see that $p(\varepsilon , t)/t^4  >0$
whenever $\varepsilon \ge 34$ or $|t| \ge 2$. Thus in order to prove
Lemma \ref{lem52}, one only has to rigorously verify that $\min_{\A}
p/t^4 >0$ for the compact rectangle $\A=[5,34]\times [-2,2]$; see Figure
\ref{figa1}. The reader may want to notice that usage of rigorous
numerics or other forms of computer assistance is not uncommon for problems of a
similar flavour in non-linear analysis; see, e.g., \cite{AL,AKP}.

(iii) Given any specific {\em rational\/} $\delta \ge \frac32$, or $\varepsilon \ge 5$, the conclusion of Lemma
\ref{lem52} may be arrived at in yet another way: Letting $\varepsilon
=m/n\ge 5$ with coprime $m,n\in \N$, notice that
$$
\frac{n^3 p(m/n, nt)}{(e^t -1)^4} = p_{m,n} (e^t) \quad \forall t\in
\R \, ,
$$
with the appropriate polynomial $p_{m,n}$ with integer coefficients
and degree $4m-n-4$. Using classical algebra tools, it may be
straightforward to see directly that $p_{m,n} (s) >0$ for all $s\in
\R^+$. For example, take $\delta=\frac32 $, hence $\varepsilon = 5$, so
$m=5$, $n=1$, and a short calculation yields
\begin{align}\label{eq518}
p_{5,1} (s) = & \enspace 30 s^{15} + 120 s^{14} + 300 s^{13} + 600
s^{12} + 798 s^{11} + 807s^{10} + 540 s^9 - 15 s^8 \\
& \enspace  - 870 s^7 -1281 s^6
-1164 s^5 - 435 s^4 + 540s^3 + 1395 s^2 + 1008 s +252 \, . \nonumber
\end{align}
By Descartes' rule, $p_{5,1}=0$ has precisely zero or two real roots
(counted with multiplicities) on $\R^+$, and the rough estimate
implied by (\ref{eq518}),
$$
p_{5,1}(s) \ge 15 \left( 213 s^9 \min \{1,s\}^6 - 251 s^4
  \max\{1,s\}^4 + 213 \min \{1,s\}^3
\right) \quad \forall s \in \R^+ \, , 
$$
can be used to show that in fact $p_{5,1}(s)>0$ for all $s\in \R^+$. Thus $p(5,t)\ge
0$ for all $t\in \R$, and hence $\omega_{f_{3/2}}$ is decreasing on $]0,1[$.
\end{rem}

As a consequence of Lemma \ref{lem52}, Theorem \ref{thm46d} together with
(\ref{eq51}) and (\ref{eq52}) yields
$$
\# \left( \bigl( 
\N \setminus \{1\} \bigr) \: \cap \: \Bigr] 2 \frac{\delta + 1}{\delta + 2} ,
\sqrt{\delta + 1} \Bigl[
\right) = \big\lceil \sqrt{\delta + 1} \, \big\rceil  - 2
$$
as the number of different non-circular Jordan solutions of $\kappa =
r^{\delta}$ modulo rotations, whenever $\delta \ge
\frac32$; here $\lceil a \rceil$ denotes the smallest
integer not smaller than $a\in \R$. By means of an obvious rescaling,
the results of this section so far can be summarized and slightly
extended as follows.

\begin{theorem}\label{thm57}
Let $a\in \R^+$, $\delta \in \R$, and assume that $\cC$ is a Jordan solution of $\kappa =
a r^{\delta}$. Then $\cC$ is oriented counter-clockwise, and the
following hold:
\begin{enumerate}
\item if $\delta \le 3$ and $\delta \ne -1, 0$ then $[\cC]$ is the
  circle with radius $a^{-1/(\delta + 1)}$ centered at $0$;
\item if $\delta = -1$ then $a=1$, and $[\cC]$ is a circle centered at $0$;
\item if $\delta = 0$ then $[\cC]$ is a circle with radius $a^{-1}$;
\item if $\delta > 3$ then either $[\cC]$ is the circle with radius
  $a^{-1/(\delta + 1)}$ centered at $0$, or else $[\cC]$ is
  non-circular, $[\cC] = a^{-1/(\delta + 1)}e^{ \textstyle i \vartheta} [\cC_{f_{\delta},s}]$ for some $\vartheta \in
  \R$ and a unique $s\in \OO_{f_{\delta}}$, with $\OO_{f_{\delta}}\subset \:\:
  ]0,1[$ containing precisely $\big\lceil \sqrt{\delta + 1} \, \big\rceil - 2$ elements. 
\end{enumerate}
\end{theorem}

\begin{proof}
Replacing $c\in \cC$ by $a^{1/(\delta + 1)}c$, it can be assumed that
$a=1$, provided that $\delta \ne -1$. Thus only the assertion
regarding $\delta = -1$ requires further justification. For $f(s) =
a/s$ with $a\in \R^+$, it is readily seen that $\mbox{\rm Per}\,
\Phi_f = \varnothing$ if $a>1$. By contrast, if $a<1$ then $0$
is a center, by Proposition \ref{prop22}, and every other orbit is
periodic and twisted. In either case, $\cC$ is not Jordan, by Theorem
\ref{thm38}, so necessarily $a=1$.
\end{proof}

\begin{figure}[ht]
\psfrag{tx1}[]{\small $\rr c$}
\psfrag{tx2}[]{\small $\ii c$}
\psfrag{th1}[]{\small $1$}
\psfrag{tee1}[]{\small $1$}
\psfrag{tv1}[]{\small $i$}
\psfrag{tvm1}[]{\small $-i $}
\psfrag{torig}[]{\small $0$}
\psfrag{tf1}[l]{\small $f(s)=e^{-s}/s$}
\psfrag{tf2}[l]{\small $f(s)=1/(s+s^3)$}
\psfrag{tfint1}[r]{\small $\kappa = r^4$}
\psfrag{tfint2}[r]{\small $\kappa = r^9$}
\psfrag{tfix1}[l]{\small $\mbox{\rm Fix} \, \Phi_f =\{ 0 \}$}
\psfrag{tfix2}[r]{\small $\mbox{\rm Fix} \, \Phi_f = \mbox{\rm Per} \,
  \Phi_f  = \{ 0 \}$}
\psfrag{tper0}[r]{\small $\omega_{f_4} \to  \frac1{\sqrt{5}}$}
\psfrag{tper0a}[l]{\small $\omega_{f_9} \to  \frac1{\sqrt{10}}$}
\psfrag{tper1}[r]{\small $[\cC_{f_4, 0.4819}]$}
\psfrag{tper1a}[r]{\small $\omega_{f_4} =  \frac12$}
\psfrag{tper2}[l]{\small $[\cC_{f_9, 0.1955}]$}
\psfrag{tper2a}[l]{\small $\omega_{f_9} =  \frac12$}
\psfrag{tper3}[r]{\small $[\cC_{f_9, 0.8477}]$}
\psfrag{tper3a}[r]{\small $\omega_{f_9} =  \frac13$}
\begin{center}
\includegraphics{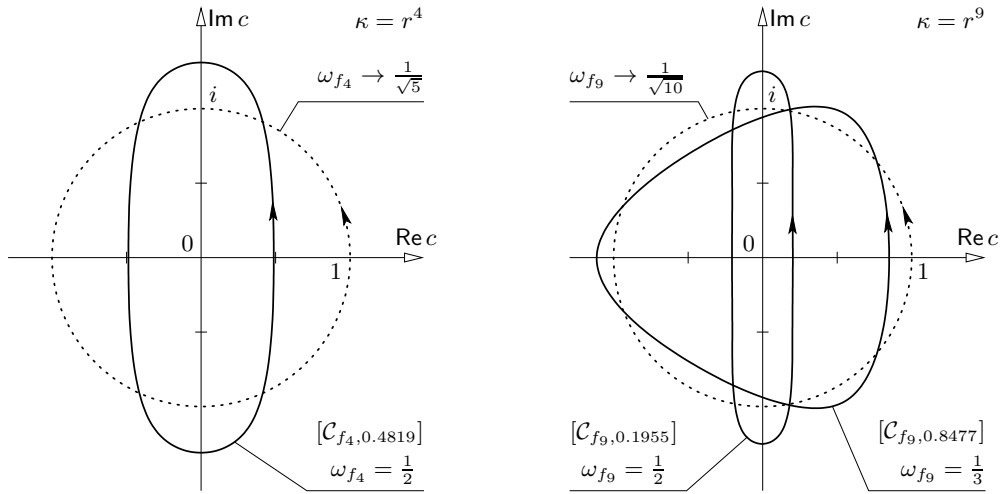}
\end{center}
\caption{Apart from the unit circle (dotted), and modulo rotations, $\kappa = r^{\delta}$ has precisely one
  (non-circular) Jordan solution when $\delta = 4$ (left; see also
  Figure \ref{fig0}), and has precisely two such solutions when $\delta = 9$ (right).}\label{figom}
\end{figure}

Notice that Theorem \ref{thm57} asserts in particular that $\kappa = a
r^{\delta}$ with $a\in\R^+$ has no non-circular Jordan solution when $\delta \le 3$,
but has, modulo rotations, precisely $n\in \N$ different
such solutions when $n^2 + 2n < \delta \le n^2 + 4n + 3$.
For instance, $\kappa= r^4$ and $\kappa
= r^9$ have precisely one and two non-circular Jordan solutions,
respectively, modulo rotations; see Figures \ref{fig0} and \ref{figom}.

\begin{rem}
As mentioned already in the Introduction, it is a well-documented
empirical observation that the oval shapes of worn stones never seem to be exact
ellipses, but rather appear to be a bit bulkier. In this regard, the reader may find it interesting to note
that none of the non-circular ovals $[\cC_{f_{\delta},s}]$ of Theorem
\ref{thm57}(iv) is an ellipse either. Indeed, suppose that for a
Jordan solution $\cC$ of $\kappa = a r^{\delta}$ with $a\in\R^+$ the set $[\cC]$ was an
ellipse with semi-axes $A\ge B>0$. If $\delta \ge 0$ then $A/B^2 = a
A^{\delta}$ as well as $B/A^2 = aB^{\delta}$, thus $a^2 A^4 =
B^{2(1-\delta)} = a^{\delta - 1} A^{(1-\delta)^2}$, and hence
$a^{3-\delta} = A^{(\delta -3)(\delta + 1)}$. It follows that either
$a=A^{-\delta - 1}$ and consequently $A=B$, or else $\delta = 3$. By Theorem \ref{thm57},
$A=B$ in any case, i.e., $[\cC]$ is a circle. (Usage of Theorem \ref{thm57}, though convenient, is not
essential here: A simple calculation shows directly that any ellipse
solving $\kappa = a r^3$ necessarily is a circle.) A similar argument
applies when $\delta < 0$.
\end{rem}

\subsection*{Supplement: The limiting shapes in \cite{andrews} revisited}

As noted earlier, the results of this section are quite specific to
the monomial family. However, the basic tools developed earlier in
order to obtain these results, notably the auxiliary planar flow $\Phi_f$ (Section
\ref{sec2}), geometric correspondances (Section \ref{sec3}), and
analytic estimates (Section \ref{sec4}), all may be useful in other
contexts as well. This supplementary section briefly describes one such
context, motivated by the classification of limiting shapes for isotropic curve
flows in \cite{andrews}. As the author intends to give a detailed account
elsewhere, only an outline is presented here that highlights
the similarity to the main results of the present article; as such, a
few non-essential assumptions are made to simplify the
exposition.

Let $g:\R^+\to \R$ be smooth, and assume for convenience that $g$ is
increasing and $g(0+)\ge 0$. Fix any $G:\R\setminus \{0\}\to
\R$ with $G'(s) = 2s - 2/g(|s|)$ for all $s\in \R \setminus
\{0\}$. Note that $G$ is convex on $\R^+$, with a global
non-degenerate minimum at $s_g$, where $s_g >0$ is the unique solution
of $sg(s)=1$. This usage of the symbol $s_g$ is consistent with
earlier usage of $s_f$. Assume henceforth that $G(s_g)=0$. On $\C
\setminus i \R$ consider an ODE for $z=z(t)$ very similar to (\ref{eq21}),
\begin{equation}\label{eq5A1}
\dot z = i z g(| \rr z |)  - i \, .
\end{equation}
Note that $G(\, \rr z ) + (\, \ii z)^2$ is a first integral of (\ref{eq5A1}). While
(\ref{eq5A1}) does not in general generate a (global, topological) flow on $\C$,
unlike (\ref{eq21}), it does generate a flow, henceforth denoted
$\Psi_g$, on the open convex set $\A_g:=\{ z\in \C : \rr z > 0 , 
G(\, \rr z) + (\, \ii z )^2 < G(0+)\}$. For instance, if $g$ can be extended
smoothly (or merely as a Lipschitz function) to $s=0$ then $\A_g =
\{z\in \C : \, \rr z > 0 \}$. In any case, every $z \in \A_g$ is a periodic point of
$\Psi_g$, and $\mbox{\rm Fix}\, \Psi_g = \{s_g\}$. Similarly to
(\ref{eq29}), define
$$
\nu_g(z) = \frac1{2\pi} \int_0^{T_g(z)} g(| \rr \Psi_g(t, z) |) \, {\rm d}t
\quad \forall z\in \A_g \, .
$$
Thus $\nu_g(s_g ) = 0$ whereas $\nu_g(z)>0$ for every $z\in \A_g
\setminus \{s_g\}$. In analogy to (\ref{eq241}) it is readily
confirmed that
$$
\lim\nolimits_{z \to s_g } \nu_g (z) = \sqrt{\frac{2}{G''(s_g)}} =
\frac1{\sqrt{1 + s_g^2 g'(s_g)}} \in \; ]0,1[ \, .
$$
Next, similarly to Section \ref{sec3} say that an oriented smooth
curve $\cC$ is a solution of 
\begin{equation}\label{eq5A4}
\kappa = g\bigl( \big|  \rr (z \overline{n})  \big| \bigr) \, ,
\end{equation}
if $\kappa_c(t) = g \bigl( \big| \rr  \bigl( c(t) (-i)
\overline{\dot c (t)}\bigr)  \big|\bigr) = g \bigl( \big| \ii \bigl(
\overline{c(t)} \dot c (t)
\bigr) \big| \bigr)$ for all $t\in \J_c$ with $\overline{c(t)} \dot c
(t) \not \in \R$, where $c$ is some (and hence any) element of
$\cC$. By the same calculation as in Section \ref{sec3}, if $\cC$ is a
solution of (\ref{eq5A4}) and $c\in \cC$ then $z_c$ given by
(\ref{eq31}) solves (\ref{eq5A1}). Thus, just as for
(\ref{eq33}), a correspondance can be established between maximal
solutions of (\ref{eq5A4}) modulo ${\rm O}(2)$-congruence on the one
hand and orbits of $\Psi_g$ on the other hand. (A careful analysis needs to pay attention to
the possibility of $\ii (\overline{c} \dot c)=0$, a situation reflected for $g(0+)=0$ by the
invariance in (\ref{eq5A1}) of the imaginary axis.) In
particular, it makes sense to define $\nu_g (\cC) = \nu_g\bigl( z_c
(0)\bigr)$ for any $c\in \cC$. With this, if $\cC$ is a Jordan solution of (\ref{eq5A4})
other than the (counter-clockwise oriented) circle with radius $s_g$
centred at $0$ then $\nu_g (\cC) = 1/n$ for some $n\in \N$. As
in the case of (\ref{eq33}), it is natural to find all Jordan
solutions of (\ref{eq5A4}), and to do so by determining the range of
$\nu_g$ as accurately as possible. With $s^*>s_g$ defined uniquely
by $G(s^*) = G(s)$ for every $0<s<s_g$, as in Section \ref{sec4},
\begin{equation}\label{eq5A5}
\nu_g(s) = \frac1{\pi} \int_s^{s^*} \frac{{\rm
    d}u}{\sqrt{G(s) - G(u)}} \quad \forall 0<s<s_g \, .
\end{equation}
This usage of the symbol $s^*$ is slightly inconsistent
with its earlier usage. However, rather than confusing the reader,
it will hopefully highlight the parallels between the analyses
of (\ref{eq33}) and (\ref{eq5A4}), respectively. Note that (\ref{eq5A5}), though very similar in spirit to
(\ref{eq43}), is considerably simpler, in at least two respects: On the one
hand, utilizing the smooth positive function $K:= G/(G')^2$, together
with (a trivially adjusted version of) Proposition \ref{prop46b}, yields
\begin{equation}\label{eq5A6}
\nu_g'(s) = \frac{G'(s)}{2 \pi G(s)} \int_s^{s^*} \frac{G'(u)
  K'(u)}{\sqrt{G(s) - G(u)}} \, {\rm d}u = \frac{G'(s)}{\pi G(s)}
\int_s^{s^*} \sqrt{G(s) - G(u)} K''(u) \, {\rm d}u \quad \forall
0<s<s_g \, ,
\end{equation}
which is far less unwieldy than (\ref{eq4oo1}). For instance, (\ref{eq5A6})
makes obvious the known fact that $s\mapsto \nu_g(s)$
is monotone whenever $K$ is convex or concave \cite{benguria, chicone,
  CW, MY, y}. Applying Proposition \ref{prop46b} once more yields, for every $0<s<s_g$,
$$
\nu_g''(s) = \frac1{2\pi G(s) K(s)} \int_s^{s^*} \frac{G(u)
  K''(u)}{\sqrt{G(s) - G(u)}} \, {\rm d}u - \frac{G'(s)K'(s)}{2\pi
  G(s)K(s)} \int_s^{s^*} \sqrt{G(s) - G(u)} K''(u) \, {\rm d}u \, .
$$
In analogy to Proposition \ref{prop46}, therefore, $\nu_g$ is smooth
on $]0,s_g[$, with
\begin{equation}\label{eq5A7}
\nu_g (s_g -) = \sqrt{\frac2{G''(s_g)}} \, , \quad \nu_g'(s_g-) = 0 \,
, \quad
\nu_g''(s_g -) = \frac{5G^{(3)}(s_g)^2 - 3 G''(s_g) G^{(4)}(s_g)}{12
  \sqrt{2} G''(s_g)^{5/2}} \, .
\end{equation}
On the other hand, notice that $\pi \nu_g (s)$ can be interpreted as the
{\em true\/} minimal period of the point $s\in \A_g$ in the planar Hamiltonian flow on
$\A_g$ generated by
\begin{equation}\label{eq5A8}
\dot w = 2 \, \ii w - i G'(\, \rr w ) \, ,
\end{equation}
which has $s_g$ as a non-degenerate center; cf.\ \cite[Sec.\ 2]{andrews}. This interpretation
makes $\nu_g$ directly amenable to the substantial literature on $1$-DOF
Hamiltonian systems, notably on the periods of such systems; see,
e.g., \cite{benguria, chicone, CD, CW, cima1, cima2,
  Rothe, villa, walter, ZP}. By contrast, much of the delicate analysis from
Proposition \ref{prop46} onward has been necessitated by the fact  that no similar interpretation
seems to exist for $\omega_f$; see also Remark \ref{rem26a}.

Finally, to illustrate the above for a familiar example, consider once
again the monomial family, that is, let $g=f_{\delta}$ with $\delta
\in \R^+$. Here, for all $s\in \R^+$,
$$
G(s) = \left\{
\begin{array}{ll}
\displaystyle s^2 + \frac2{\delta - 1} s^{1 -\delta} - \frac{\delta+
  1}{\delta - 1} & \mbox{\rm if } \delta \ne 1
\, , \\[2mm]
s^2 - 2 \log s - 1 & \mbox{\rm if } \delta = 1 \, ,
\end{array}
\right.
$$
and $\A_{f_{\delta}} = \{z\in \C : \, \rr z > 0\}$ precisely if $\delta \ge
1$. Also, $s_{f_{\delta}} = 1$, and (\ref{eq5A7}) yields
$$
\nu_{f_{\delta}} (1-) = \frac1{\sqrt{\delta + 1}} \, , \quad
\nu_{f_{\delta}}'(1-) =0 \, , \quad \nu_{f_{\delta}}''(1-) =
\frac{\delta (\delta - 3)}{12 \sqrt{\delta + 1}} \, . 
$$
Thus, as $s\to 1$ the function $\nu_{f_{\delta}}$ attains a
non-degenerate maximum (respectively, minimum) if $0<\delta < 3$
(respectively, if $\delta > 3$).  In addition, a straightforward calculation
shows that
$$
\nu_{f_{\delta}} (0+) = \frac1{1 + \min \{\delta, 1\}} \quad \forall
\delta \in \R^+ \, .
$$
As a consequence, $(\delta-3) \bigl( \nu_{f_{\delta}} (0+) -
\nu_{f_{\delta}} (1-) \bigr)> 0$ for every $\delta \in \R^+\setminus \{3\}$, which
certainly makes it plausible to speculate that $\nu_{f_{\delta}}$ is
increasing (respectively, decreasing) on $]0,1[$ if $0<\delta < 3$
(respectively, if $\delta > 3$). Notice how this is the precise analogue
of Conjecture \ref{con51}. Unlike for the latter, it is not 
hard to prove this speculation correct in its entirety. Given that
(\ref{eq5A5}) is considerably simpler than (\ref{eq43}), as noted
above, this may not come as a complete surprise \cite{andrews}.

\begin{prop}\label{prop5o1}
On $]0,1[$, the function $\nu_{f_{\delta}}$ is increasing for $0<\delta
< 3$, and decreasing for $\delta > 3$.
\end{prop}

As Proposition \ref{prop5o1} suggests, the case $\delta =3$ is
somewhat special: Indeed, here $s^*=1/s$ for every $0<s<1$, and hence
$$
\nu_{f_3} (s) = \frac1{\pi} \int_s^{1/s} \frac{{\rm d}u}{
  \sqrt{(s-1/s)^2 - (u-1/u)^2}} = \frac1{2\pi} \int_{s^2}^{1/s^2}
\frac{{\rm d}u}{\sqrt{(u-s^2)(1/s^2 - u)}} = \frac12 \quad \forall 0<s<1
\, .
$$
In other words, the center $1$ of (\ref{eq5A8}) is {\em isochronous\/} for
$g=f_3$. This ``surprising affine invariance property''\cite{andrews}
reflects the fact that the phase portrait of (\ref{eq5A8}) for
$g=f_3$ is invariant under the diffeomorphism $z \mapsto (\, \rr
z)^{-1} - i \, \ii z$ of $\A_{f_3}$. Further analysis shows that every orbit
$\Psi_{f_3} (\R, s)$ with $s\in \R^+$ corresponds to a
(counter-clockwise oriented) ellipse with semi-axes $s,1/s$. Thus,
every maximal solution of $\kappa = | \rr ( z \overline{ n}) |^3$ is
an ellipse centered at $0$ with interior area $\pi$.

\begin{figure}[ht]
\psfrag{tx1}[]{\small $s$}
\psfrag{tx2}[l]{\small $\omega_{f_{\delta}}(s)$}
\psfrag{tx3}[l]{\small $\nu_{f_{\delta}}(s)$}
\psfrag{tvl1}[]{\small $1$}
\psfrag{thl1}[]{\small $1$}
\psfrag{thl0}[]{\small $0$}
\psfrag{tvl12}[]{\small $\frac12$}
\psfrag{tvl13}[]{\small $\frac13$}
\psfrag{tvl1n}[]{\small $\frac1{n}$}
\psfrag{tdell05}[l]{\small $\delta = \frac12$}
\psfrag{tdell0}[l]{\small $(\delta = 0)$}
\psfrag{tdell1}[l]{\small $\delta = 1$}
\psfrag{tdell15}[l]{\small $\delta = \frac32$}
\psfrag{tdell3}[l]{\small $\delta = 3$}
\psfrag{tdell4}[l]{\small $\delta = 4$}
\psfrag{tdell9}[l]{\small $\delta = 9$}
\psfrag{tdell8}[l]{\small $\delta = 8$}
\psfrag{tdelln2m1}[l]{\small $\delta = n^2 -1$}
\begin{center}
\includegraphics{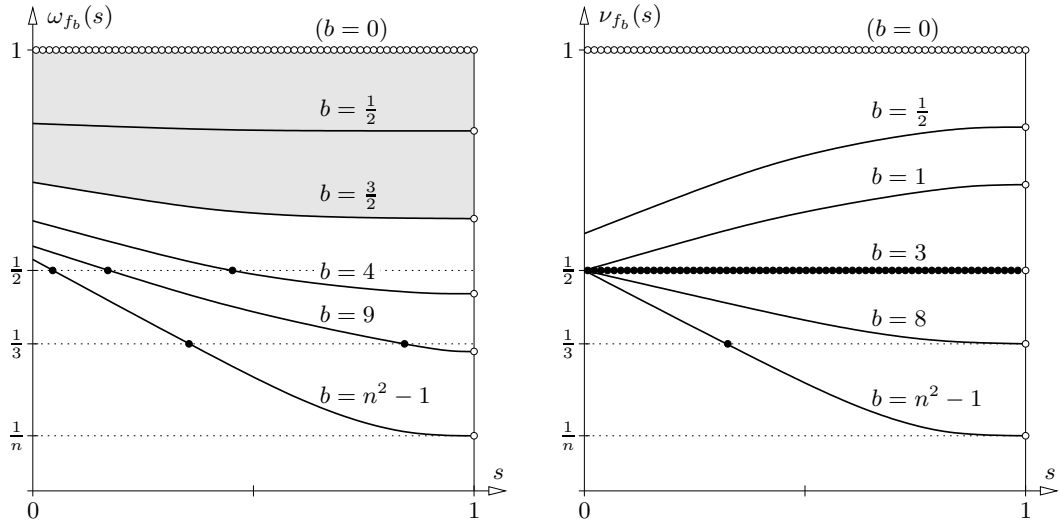}
\end{center}
\caption{Finding all Jordan solutions of $\kappa = r^{\delta}$ (left)
  and $\kappa = | \rr (z \overline{n})|^{\delta}$ (right) with $\delta\in \R^+$, by
  qualitatively graphing $\omega_{f_{\delta}}$ and
  $\nu_{f_{\delta}}$, respectively. At the time of this writing, monotonicity (in $s$) of
  $\omega_{f_{\delta}}$ for $0<\delta < \frac32$ is conjectural only (grey
  region). Solid black dots indicate non-circular Jordan solutions,
  while circles represent circular solutions with radius $1$.}\label{fig9}
\end{figure}

As a consequence of Proposition \ref{prop5o1}, for
every $0<\delta \le 8$ with $\delta \ne 3$ the only Jordan solution of
\begin{equation}\label{eq5A9}
\kappa = | \rr ( z \overline{  n}) |^{\delta}
\end{equation}
is the (counter-clockwise oriented) unit circle, whereas for $\delta > 8$ there exist precisely
$\bigl\lceil \sqrt{\delta +1} \bigr\rceil - 3$ different non-circular Jordan
solutions of (\ref{eq5A9}), modulo rotations; see Figure \ref{fig9} and \cite[Thm.\ 5.1]{andrews}.

\subsubsection*{Acknowledgements}
The author was partially supported by an {\sc Nserc} Discovery
Grant. He owes deep gratitude to T.P.\ Hill and K.E.\ Morrison who in 2005
conjectured (correctly, as it turned out) that the only Jordan
solution of $\kappa = r$ is the (counter-clockwise oriented) unit
circle, and who greatly helped this work come to fruition
through continued interest and advice. Insightful comments by an anonymous
referee led to a much improved presentation. Thanks also to J.\
Muldowney, M.\ Niksirat, T.\ Schmah, and C.\ Xu for several
enlightening conversations over the years. Parts of this work were
completed while the author was a visitor at the Universit\"{a}t
Wien. He is much indebted to R.\ Zweim\"{u}ller for many kind acts of
hospitality.

\end{document}